\documentclass[3p,11pt]{elsarticle}

\usepackage{url}
\usepackage{amssymb}
\usepackage{amsmath}
\usepackage{amsthm}
\usepackage{stmaryrd}
\usepackage{siunitx}
\usepackage{commath}
\usepackage{algorithm}
\usepackage{algorithmic}
\usepackage{subfig}
\usepackage{enumerate}
\usepackage{multirow}
\usepackage{makecell}
\usepackage{lipsum}
\usepackage{amsfonts}
\usepackage{graphicx}
\usepackage{bm}
\usepackage{color}
\usepackage{mathrsfs}
\usepackage{caption}
\usepackage{fontspec}

\usepackage[hidelinks]{hyperref}
\hypersetup{
colorlinks   = true, 
urlcolor     = red, 
linkcolor    = red, 
citecolor   = red 
}

\biboptions{sort&compress}
\usepackage[misc]{ifsym}
\usepackage{color}
\usepackage{cleveref}

\theoremstyle{definition}
\newtheorem{definition}{Definition}[section]

\newtheorem{remark}{Remark}[section]

\makeatletter
\def\ps@pprintTitle{%
\let\@oddhead\@empty
\let\@evenhead\@empty
\def\@oddfoot{}%
\let\@evenfoot\@oddfoot}

\newcommand{\tnorm}{\@ifstar\@tnorms\@tnorm}
\newcommand{\@tnorms}[1]{%
\left|\mkern-1.5mu\left|\mkern-1.5mu\left|
#1
\right|\mkern-1.5mu\right|\mkern-1.5mu\right|
}
\newcommand{\@tnorm}[2][]{%
\mathopen{#1|\mkern-1.5mu#1|\mkern-1.5mu#1|}
#2
\mathclose{#1|\mkern-1.5mu#1|\mkern-1.5mu#1|}
}

\renewenvironment{proof}[1][\proofname]{\par
\vspace{-\topsep}
\pushQED{\qed}%
\normalfont
\topsep0pt\partopsep0pt 
\trivlist
\item[\hskip\labelsep
\itshape
#1\@addpunct{.}]\ignorespaces
}{%
\popQED\endtrivlist\@endpefalse
\addvspace{6pt plus 6pt} 
}
\makeatother

\newtheoremstyle{mystyle}
{3pt} 
{1pt} 
{} 
{} 
{\bfseries} 
{.} 
{.5em} 
{} 

\theoremstyle{mystyle} 
\newtheorem{theorem}{Theorem}[section]
\newtheorem{lemma}{Lemma}[section]

\allowdisplaybreaks[4]

\theoremstyle{definition}
\newtheorem{example}[theorem]{Example}

\theoremstyle{remark}

\numberwithin{equation}{section}

\usepackage{float}
\usepackage{comment}

\makeatletter

\@addtoreset{equation}{section}

\makeatother

\begin{document}
\abovedisplayskip=4.2pt plus 4.2pt
\abovedisplayshortskip=0pt plus 4.2pt
\belowdisplayskip=4.2pt plus 4.2pt
\belowdisplayshortskip=2.1pt plus 4.2pt
\begin{frontmatter}
\title{Functional normalizing flow for statistical inverse problems of partial differential equations}

\author[1]{Yang Zhao}

\author[1]{Haoyu Lu}

\author[1]{Junxiong Jia \texorpdfstring{\corref{first}}{}}
\cortext[first]{Corresponding author: Junxiong Jia (jjx323@xjtu.edu.cn)}

\author[2]{Tao Zhou}

\address[1]{School of Mathematics and Statistics, Xi'an Jiaotong University,
	Xi'an, Shaanxi 710049, China}
\address[2]{ Academy of Mathematics and Systems Sciences, Chinese Academy of Sciences,
    Beijing, 100190, China}

\begin{abstract}
Inverse problems of partial differential equations are ubiquitous across various scientific disciplines and can be formulated as statistical inference problems using Bayes' theorem. To address large-scale problems, it is crucial to develop discretization-invariant algorithms, which can be achieved by formulating methods directly in infinite-dimensional space. We propose a novel normalizing flow based infinite-dimensional variational inference method (NF-iVI) to extract posterior information efficiently. Specifically, by introducing well-defined transformations, the prior in Bayes' formula is transformed into post-transformed measures that approximate the posterior. To circumvent the issue of mutually singular probability measures, we formulate general conditions for the employed transformations. As guiding principles, these conditions yield four concrete transformations. Additionally, to minimize computational demands, we have developed a conditional normalizing flow variant, termed CNF-iVI, which is adapt at processing measurement data of varying dimensions while requiring minimal computational resources. We apply the proposed algorithms to three typical inverse problems governed by the simple smooth equation, the steady-state Darcy flow equation, and the electric impedance tomography. Numerical results confirm our theoretical findings, illustrate the efficiency of our algorithms, and verify the discretization-invariant property.
\end{abstract}

\begin{keyword}
inverse problems, infinite-dimensional variational inference, functional normalizing flow, Bayesian analysis for functions, partial differential equations.
\end{keyword}
\end{frontmatter}

\section{Introduction}
Driven by their widespread applications in seismic exploration, radar imaging, and other fields, inverse problems involving partial differential equations (PDEs) have witnessed significant advancements in recent decades \cite{Kirsch2011Book,Benning2018ActaNum}. As computational power continues to grow, researchers are increasingly focused on not only obtaining estimated solutions but also conducting statistical analyses to quantify uncertainties, which is crucial for tasks such as artifact detection \cite{Zhou2020SJIS}. The Bayesian inverse approach offers a robust framework for addressing inverse problems within the context of PDEs by transforming them into statistical inference problems, thereby enabling the analysis of parameter uncertainties \cite{dashti2013bayesian}. 

Typically, inverse problems of PDEs are posed in infinite-dimensional spaces \cite{Kirsch2011Book}, which presents challenges for applying well-established finite-dimensional Bayesian inference techniques \cite{bishop2006pattern,kaipio2006statistical}. To bridge this gap, two predominant strategies have been developed: 
\begin{itemize} 
\item Discretize-then-Bayesianize: PDEs are initially discretized to approximate the original problem in some finite-dimensional space, and the reduced, approximated problem is then solved using Bayes' method \cite{kaipio2006statistical}. 
\item Bayesianize-then-discretize: Bayes' formula and algorithms are initially constructed in infinite-dimensional space, and after the infinite-dimensional algorithm is built, some finite-dimensional approximation is carried out \cite{dashti2013bayesian}. 
\end{itemize}

Both approaches offer distinct advantages and disadvantages. The \emph{Discretize-then-Bayesianize} approach allows us to employ all the Bayesian inference methods developed in the statistical literature \cite{kaipio2006statistical,Calvetti2023Book} to solve the inverse problems. However, due to the infinite-dimensional nature of the original problems, two critical challenges arise: 
\begin{itemize} 
\item Model consistency: Finite-dimensional models, such as those representing prior probability measures, exhibit distinct properties compared to their infinite-dimensional counterparts. A notable example is the total variation prior, which has been extensively studied in \cite{lassas2004can,Dashti2012IPI}. 
\item Algorithm applicability: To preserve the intrinsic structure of original infinite-dimensional problems, algorithms initially designed for finite-dimensional spaces must be carefully adapted to the infinite-dimensional setting. This reformulation ensures consistent algorithmic behavior across various discretizations, as demonstrated in \cite{Jia2021IPI,jia2021variational,cotter2013,Jia2022VINet,sui2024noncentered}. 
\end{itemize} 
To address these challenges, the \emph{Bayesianize-then-discretize} approach has attracted significant attention from researchers in recent years \cite{bui2013computational,Cotter2009IP,dashti2013bayesian}.

A primary challenge in Bayesian inference is the efficient extraction of posterior information. From the \emph{Bayesianize-then-discretize} perspective, an infinite-dimensional Markov chain Monte Carlo (MCMC) algorithm known as preconditioned Crank-Nicolson (pCN) has been proposed and analyzed in detail \cite{cotter2013,Pillai2014SPDE}. This algorithm maintains consistent sampling efficiency across different discretizations. In addition to the pCN algorithm, other types of infinite-dimensional sampling algorithms have been proposed, such as the infinite-dimensional sequential Monte Carlo algorithm \cite{Beskos2015SC} and the infinite-dimensional importance sampling algorithm \cite{Agapiou2017SS}. To enhance sampling efficiency, infinite-dimensional MCMC algorithms with gradient and geometric informative proposals have been developed. Examples include the infinite-dimensional Metropolis-adjusted Langevin algorithm \cite{Thanh2016IPI} and the geometric pCN algorithm \cite{Beskos2017JCP}. Although these algorithms possess mesh independence properties, they are challenging to apply to some large-scale inverse problems of PDEs \cite{Fichtner2011Book}.

In the finite-dimensional setting, variational inference (VI) methods have been extensively studied in machine learning to mitigate the computational demands of MCMC sampling algorithms \cite{zhang2018advances}. For example, a mean-field assumption based VI approach was employed to solve finite-dimensional inverse problems with hyper-parameters in prior and noise distributions \cite{Guha2015JCP, Jin2012JCP, jin2010hierarchical}.
Projected Stein variational gradient descent methods were constructed to solve inverse problems with low-intrinsic dimensions \cite{Chen2021SISC}, and normalizing flow in Euclidean space have been extensively studied as a means to approximate target posterior distributions by transforming simple distributions  \cite{planarflow,sylvesterflow,householderflow1,MR4749392}.

Here, we briefly outline the key distinctions between MCMC sampling and variational inference (see \cite{blei2017variational} for details). While MCMC methods typically require higher computational resources, they guarantee asymptotically exact samples from the target distribution \cite{robert2004monte}. In contrast, variational inference seeks an approximate distribution that is computationally tractable, leading to significant computational gains despite some loss of accuracy compared to MCMC. Variational inference's compatibility with optimization techniques such as stochastic gradient descent and distributed optimization \cite{robbins1951stochastic, kushner2003stochastic, welling2011bayesian} makes it well-suited for large datasets and rapid model exploration. MCMC excels in scenarios with smaller datasets where precise sampling is prioritized. Another critical factor influencing the choice between MCMC and variational inference is the geometry of the posterior distribution. For instance, mixture models often exhibit multiple modes corresponding to different label permutations. In such cases, Gibbs sampling excels due to its ability to efficiently explore these modes \cite{bishop2006pattern}. However, when Gibbs sampling is not feasible, variational inference can outperform general MCMC methods, such as Hamiltonian Monte Carlo, even for smaller datasets \cite{kucukelbir2017automatic}.

In contrast to the extensive body of research in finite-dimensional settings, the application of variational inference methods to infinite-dimensional settings remains relatively unexplored. Generally speaking, there are two main research directions for VI methods in infinite-dimensional settings: 
\begin{itemize} 
\item Parameterize measures directly: This involves directly parameterizing the approximate measure and optimizing these parameters to approximate the posterior. For instance, when restricting approximations to Gaussian measures, a novel Robbins-Monro algorithm was developed using a calculus-of-variations perspective \cite{pinski2015algorithms, Pinski2015SIAMMA}. Under the classical mean-field assumption, a general VI framework within separable Hilbert spaces was recently proposed \cite{jia2021variational, Jia2022VINet, sui2024noncentered}. 
\item Parameterize transformation: By parameterizing a transformation, a simple measure can be mapped into a more complex measure to approximate the posterior. Optimization of these transformation parameters drives the approximation process. For example, the infinite-dimensional Stein variational gradient descent introduced in \cite{jia2021stein} represents a function space particle optimization method with rigorous mathematical underpinnings in separable Hilbert spaces. 
\end{itemize}

In this work, we focus on defining normalizing flow in infinite-dimensional spaces, aligning with the second research direction. Specifically, we derive general transformation forms capable of mapping simple measures to complex distributions within the function space. 
{\color{black}{It is noteworthy that the definition of normalizing flow in function space has been considered before in \cite{shiuniversal}. However, the objectives of this work differ significantly from ours. In \cite{shiuniversal}, the authors focus on constructing a generative model to fit a collection of functions $\mathbb{D}=\{u_i\}_{i=1}^N$ resides in some function spaces, rather than approximating the posterior of Bayes' formula. More importantly, our model establishes a rigorous theory guaranteeing the equivalence of the measure before and after the transformation, which provides a solid theoretical foundation for posterior approximation.}} Precisely, it is crucial to ensure that the post-transformed measure is absolutely continuous with respect to the prior measure when constructing the normalizing flow-based infinite-dimensional variational inference (NF-iVI) algorithm. Failure to ensure this condition would render the Radon-Nikodym (RN) derivative of the post-transformed measure with respect to the prior measure meaningless, thus contradicting the entire theoretical framework of the infinite-dimensional variational inference algorithm \cite{Pinski2015SIAMMA, Jia2022VINet}. Leveraging the properties of Gaussian measures \cite{gaussianmeasure}, we develop a rigorous theoretical framework for transformations within functional normalizing flow. By formulating specific conditions, our approach guarantees the equivalence between the post-transformed and pre-transformed measures (detailed in Subsection \ref{subsec:normalizing flow in function space}), ensuring the model's well-definedness in function space.

In addition, we have devised a strategy to alleviate the computational demands of NF-iVI. Variational inference algorithms require analyzing the Kullback-Leibler (KL) divergence between the approximate and target posterior measures, which involves taking an expectation over the approximate measure. This typically necessitates computationally intensive Monte Carlo estimation. In our algorithm, the derivative of the KL divergence with respect to the parameters of the approximate measure must be calculated at each iteration. This necessitates computing the KL divergence and its gradient at each iteration, significantly increasing computational cost. While the NF-iVI algorithm offers advantages over traditional methods like pCN, it still presents considerable computational challenges (see Section \ref{sec4} for details).
We observe that for a fixed inverse problem model, distinct measurement data induce distinct posterior distributions, necessitating computationally intensive model retraining for each new set of measurement data. To enhance the computational efficiency of our algorithm when processing new measurement data, we propose to incorporate a conditional component into functional normalizing flow. By introducing measurement data as a conditional input, we enable the model to adapt to different measurement data without requiring extensive retraining. When the conditional network is trained, we can directly obtain an acceptable estimate of the posterior for any measurement data. We refer to this approach \textbf{c}onditional \textbf{n}ormalizing \textbf{f}low based \textbf{i}nfinite-dimensional \textbf{v}ariational \textbf{i}nference (CNF-iVI). Furthermore, if we want to obtain a more accurate approximation, we can consider further training for specific measurement data based on the pre-trained conditional functional normalizing flow.

In summary, this work mainly contains four contributions:
 \begin{itemize}
 \item We introduce an infinite-dimensional model called functional normalizing flow, which is capable of transforming a simple measure into a complex approximation of the target measure. To ensure the validity of these transformations, we establish a rigorous theoretical framework that guarantees the equivalence between the transformed and pre-transformed measures.
 \item We introduce two linear flows: functional Householder flow and functional projected transformation flow, and two nonlinear flows: functional planar flow and functional Sylvester flow, as exact examples. Rigorous proofs confirm that all four proposed models satisfy the conditions of the theoretical framework we established for functional normalizing flow.
 \item We have developed a model called conditional functional normalizing flow by incorporating a conditional component. This component enables the direct entry of measurement data of any length, yielding an acceptable posterior estimate. To achieve greater accuracy, the model can be further refined through targeted training on specific measurement data.
 \item We have applied NF-iVI and CNF-iVI to typical linear and nonlinear inverse problems, specifically elliptic inverse problems, steady-state Darcy flow inverse permeability estimations and electric impedance tomography, thereby validating their efficiency.
\end{itemize}

The outline of this paper is as follows: In Section \ref{sec:functional normalizing flow}, we briefly introduce the theoretical framework and algorithmic process of NF-iVI. In Subsection \ref{subsec:bayes and vi}, we present the infinite-dimensional Bayesian theory proposed in \cite{pinski2015algorithms,Pinski2015SIAMMA} and provide a brief overview of infinite-dimensional variational inference algorithms. In Subsection \ref{subsec:normalizing flow in function space}, we introduce the NF-iVI algorithm in function spaces. In Subsection \ref{subsec:neural_op_connection}, we provide a detailed comparison between the conventional neural operator and the functional normalizing flow. In Subsections \ref{subsec:linear transformation} and \ref{subsec:nonleaner transformation}, we provide four examples of functional normalizing flow. In Subsection \ref{subsec:discrete invariance}, we prove the discretization-invariance of the proposed functional normalizing flows. In Section \ref{sec3}, based on NF-iVI, we introduce CNF-iVI, adding an auxiliary conditional input to enable the functional normalizing flow parameters to be determined by the measurement information. The detailed training procedure for CNF-iVI is also provided. Furthermore, we describe a method to improve the results of conditional functional normalizing flow. In Sections \ref{sec4} and \ref{sec5}, we apply the algorithms to three typical inverse problems, verifying the feasibility of both the NF-iVI and CNF-iVI algorithms. Additionally, in each numerical simulation, we demonstrate the mesh independence as expected for the \emph{Bayesianize-then-discretize} approach. In Section \ref{sec6}, we summarize our achievements, acknowledge some deficiencies, and explore further research directions.

\section{Functional Normalizing Flow}\label{sec:functional normalizing flow}
{\color{black} In this section, we establish the theoretical framework of the functional normalizing flow and demonstrate its specific application to the PDE inverse problem involving infinite-dimensional parameters and finite-dimensional measurements.} To circumvent the obstacle posed by the singularity of measures in infinite-dimensional spaces, we introduce a novel theoretical framework. Building upon this foundation, four specific flow models are proposed and shown to satisfy the established theoretical criteria.
\subsection{Bayesian Approach and Variational Inference}\label{subsec:bayes and vi}
For a large class of inverse problems involving partial differential equations, sparse measurement data are practically adopted \cite{Cotter2009IP, nickl2020bernstein}, as such data are more easily acquired. {\color{black} In this subsection, we introduce the foundational concepts of infinite-dimensional Bayesian inverse problems with finite-dimensional sparse data, as well as the underlying principles of infinite-dimensional variational inference. }

Let $\mathcal{H}_u$, $\mathcal{H}_w$ be separable Hilbert spaces representing the parameter space and solution space, respectively, and $N_d$ be a positive integer. 
Denote $\mathcal{N}(u, \mathcal{C}_0)$ as a Gaussian measure with mean $u$ and covariance operator $\mathcal{C}_0$.
The inverse problem can be described as
\begin{align}\label{eq:ca}
 \bm{d} = \mathcal{S}\mathcal{G}(u) + \bm{\epsilon},
\end{align}
where $\bm{d}\in \mathbb{R}^{N_d}$ is the measurement data, $u\in \mathcal{H}_u$ is the parameter of interest, $\mathcal{G}$ is the PDE solution operator from $\mathcal{H}_u$ to $\mathcal{H}_{w}$, $\mathcal{S}$ is the measurement operator from $\mathcal{H}_{w}$ to $\mathbb{R}^{N_d}$, and $\bm{\epsilon}$ is a Gaussian random vector with zero mean and covariance matrix  $\bm{\Gamma}_{\text{noise}} := \tau^{-1}\textbf{I}$ ($\tau$ is a fixed positive number, and $\textbf{I}$ denotes the matrix of identity), which means
\begin{align}\label{eq:noise}
 \bm{\epsilon} \sim \mathcal{N}(0, \bm{\Gamma}_{\text{noise}}).
\end{align}

Based on the framework of infinite-dimensional Bayesian inference \cite{weglein2003inverse}, we are able to preserve the fundamental Bayes' formula for the inverse problem:
\begin{align}\label{equ:bayes infinite}
	\qquad\,\,\frac{d\mu}{d\mu_0}(u) = \frac{1}{Z_{\mu}}\exp (-\Phi(u) ),
\end{align}
where $\Phi : \mathcal{H}_u \rightarrow \mathbb{R}$ is defined as
\begin{align}\label{eq:likelihood}
 \Phi(u) = \frac{1}{2}\lVert \bm{d} - \mathcal{S}\mathcal{G}(u)\rVert ^{2}_{\bm{\Gamma}_{\text{noise}}}
\end{align}
with 
$\lVert  \cdot  \rVert_{\bm{\Gamma}_{\text{noise}}}:=\lVert \bm{\Gamma}_{\text{noise}}^{-1/2}   \cdot\rVert$, 
and $Z_{\mu}$ is a positive finite constant given by 
\begin{align*}
Z_{\mu} = \int_{\mathcal{H}_u} \exp (-\Phi(u))\mu_0(du).
\end{align*}

Directly extracting quantitative information from the posterior distribution in Bayes' formula (\ref{equ:bayes infinite}) is computationally prohibitive \cite{cotter2013}. To circumvent this challenge, variational inference approximates the intractable posterior with a more tractable measure. 
We select $\mathcal{M}(\mathcal{H}_u)$, a set of measures on the space $\mathcal{H}_u$, as the approximating measure set. To ensure that the Radon-Nikodym (RN) derivative is well-defined, every measure in $\mathcal{M}(\mathcal{H}_u)$ is required to be equivalent to the prior measure $\mu_0$ \cite{Pinski2015SIAMMA}.
 
For any $\nu \in \mathcal{M}(\mathcal{H}_u)$,  the KL divergence between $\nu$ and posterior $\mu$ is given by
\begin{align*}
	D_{\text{KL}}(\nu||\mu)=\int_{\mathcal{H}_u} \ln\left(\frac{d\nu}{d\mu}(u)\right)\nu(dx)=\int_{\mathcal{H}_u} \ln\left(\frac{d\nu}{d\mu_0}(u)\right)-\ln\left(\frac{d\mu}{d\mu_0}(u)\right)\nu(du).
\end{align*}
The primary objective of variational inference is to identify a measure, denoted $\nu^*$, that minimizes the KL divergence between an approximate distribution $\nu$ and the target posterior distribution $\mu$:
\begin{align}\label{equ:KLmin}
\nu^*=\mathop{\arg\min}\limits_{\nu\in \mathcal{M}(\mathcal{H}_u)}D_{\text{KL}}(\nu||\mu).
\end{align}
Solving the minimization problem (\ref{equ:KLmin}) typically requires a suitable parameterization of the set $\mathcal{M}(\mathcal{H}_u)$: Overly intricate measures can hinder optimization, while overly simplistic ones may compromise approximate accuracy. To address this trade-off, \cite{Pinski2015SIAMMA,pinski2015algorithms} employed carefully designed Gaussian measures, whereas \cite{jia2021variational, sui2024noncentered, Jia2022VINet} opted for the mean-field approach.

\subsection{Normalizing Flow in Function Space}\label{subsec:normalizing flow in function space}
In this paper, we draw inspiration from the ideas presented in Euclidean space \cite{planarflow}. A series of parametric transformations will be applied to a simple measure to generate a flexible set of complex measures. This collection of measures forms the approximating measures set $\mathcal{M}(\mathcal{H}_u)$ for variational inference as outlined in (\ref{equ:KLmin}). Consistent with normalizing flow in Euclidean space \cite{planarflow,sylvesterflow,householderflow1}, we refer to this set of models as functional normalizing flow.

The core of the model consists of two components: the pre-transformed measure and the transformations. The pre-transformed measure should be analytically tractable to enable efficient sampling, while the transformations must be flexible to convert this simple measure into a complex one that accurately approximates the target measure. We select the prior $\mu_0$ of Bayes' formula (\ref{equ:bayes infinite}) as the pre-transformed measure. The transformations are a series of parameterized operators, denoted as 
\( f_{\theta_n}^{(n)} \), mapping from \( \mathcal{H}_u \) to \( \mathcal{H}_u \). Here, \( n = 1, 2, \ldots, N \), and \( \theta_n \in \Theta_n \) represents the parameters of \( f_{\theta_n}^{(n)} \). The space \(\Theta_n\) is a parameter space, whose specific form depends on the parameterization of \( f_{\theta_n}^{(n)} \).
For instance, in the context of functional planar flow illustrated in Subsection \ref{subsec:nonleaner transformation}, the parameter space $\Theta_n$ is defined as follows: 
\begin{align*}
    \Theta_n=\{(w_n, u_n, b_n) \in \mathcal{H} \times \mathcal{H} \times \mathbb{R} \,\, | \,\, \langle u_n, w_n \rangle_{\mathcal{H}_u}>-1\},
\end{align*}
where $\mathcal{H}$ is the Cameron-Martin space of the prior $\mu_0$. 
Let $\theta = \{\theta_1, \theta_2, \ldots, \theta_N\}$ denote the collection of parameters. The post-transformed measures are constructed by composing a series of transformations, which convert the pre-transformed measure into a complex one. 
We write that
\begin{align*}
    f_{\theta}(u)=f_{\theta_{N}}^{(N)} \circ \cdots \circ f_{\theta_2}^{(2)} \circ f_{\theta_1}^{(1)}(u),
\end{align*}
and denote $u_{\theta}^{N}:=f_{\theta}(u_0)$ where $u_0$ is a sample from $\mu_0$. Then $u_{\theta}^N$ will be a sample from the measure $\mu_{f_{\theta}}=\mu_0 \circ f_{\theta}^{-1}$, where $\mu_0 \circ f_{\theta}^{-1}$ is the law of $f_{\theta}$ with respect to $\mu_0$ \cite{gaussianmeasure}. 
By varying the parameters $\theta$, the transformed measures $\mu_{f_{\theta}}$ form a set of measures on 
$\mathcal{H}_u$, denoted by
\begin{align*}
\mathcal{M}(\mathcal{H}_u)=\left\{ \mu_{f_{\theta}} \, |  \, \theta \in \Theta \right\},
\end{align*}
where $\Theta=\Theta_1 \times \Theta_2 \times \cdots \times \Theta_N$ denotes the space of all possible values of  $\theta$. We let the set $\mathcal{M}(\mathcal{H}_u)$ constitute the approximating measures for the variational inference defined in (\ref{equ:KLmin}). 

Functional normalizing flow gives us a way to specify the approximate measures required for variational inference. However, since all the measures in $\mathcal{M}(\mathcal{H}_u)$ should be equivalent to the prior $\mu_0$ \cite{Pinski2015SIAMMA}, the choice of the parametric transformations $\{f_{\theta_n}^{(n)}\}_{n=1}^{N}$ may pose significant challenges. The measures within the infinite-dimensional function space are inherently prone to be singular with each others. To elucidate this concept, we present two pertinent examples:
\begin{example}
Let \(\mu_0\) be a Gaussian measure on \(L^2(\mathbb{R})\), and let \(\mathcal{H}(\mu_0)\) be its Cameron-Martin space. Assume \(m \notin \mathcal{H}(\mu_0)\), define \( f_m(u) = u + m \), and let \(\mu_m = \mu_0 \circ f_m^{-1}\). We know that \(\mu_m\) is singular with respect to \(\mu_0\).
\end{example}
\begin{example}
Let \(\mu_0\) be a Gaussian measure on \(L^2(\mathbb{R})\), \( f_2(u) = 2u \), and \(\mu_2 = \mu_0 \circ f_2^{-1}\), then we know that \(\mu_2\) is singular with respect to \(\mu_0\).
\end{example}
It is noteworthy that although the transformations $f_{m}$ and $f_{2}$ are simple, they will convert $\mu_0$ into measures that are singular with respect to $\mu_0$, which contradicts the theoretical framework of infinite-dimensional variational inference. Consequently, the judicious selection of transformations is essential for normalizing flow in function space, requiring both model flexibility and the preservation of measure equivalence. Here, we introduce a general theorem characterizing the transformations of functional normalizing flow. Assuming $\mathcal{F}_{\theta_{n}}^{(n)}$ is an operator mapping $\mathcal{H}_u$ to itself for $n = 1, 2, \ldots, N$, parameterized by $\theta_{n}$, $D\mathcal{F}_{\theta_n}^{(n)}(u)$ represents the Fréchet derivative of $\mathcal{F}_{\theta_n}^{(n)}$, and $\theta = \{\theta_1 , \theta_2, \ldots, \theta_N\}$. We define transformations as follows:
\begin{equation}\label{equ:transformation}
    f_{\theta_n}^{(n)}(u)=u+\mathcal{F}_{\theta_n}^{(n)}(u).
\end{equation}
We now present the following theorem, which establishes conditions ensuring the equivalence of measures before and after the transformations defined in equation (\ref{equ:transformation}), and provides their corresponding RN derivative.

{\color{black}
\begin{theorem}\label{thm:RN}
Let $\mathcal{H}_u$ be a separable Hilbert space equipped with the Gaussian measure  $\mu_0 = \mathcal{N}(0,\mathcal{C}_0)$, and $\mathcal{H}=\mathcal{H}(\mu_0)$ be the Cameron-Martin space of the measure $\mu_0$.  For each $n = 1, 2, \ldots ,  N$, let $f_{\theta_n}^{(n)}(u)=u+\mathcal{F}_{\theta_n}^{(n)}(u)$, where $\mathcal{F}_{\theta_n}^{(n)}$ is an operator from $\mathcal{H}_u$ to $\mathcal{H}_u$.  The composite transformation $f_{\theta}$ is given by  $f_{\theta}(u)=f_{\theta_N}^{(N)} \circ f_{\theta_{N-1}}^{(N-1)} \circ \cdots \circ f_{\theta_1}^{(1)}(u)$, and $\mu_{f_{\theta}} $ denotes its corresponding push-forward measure by  $\mu_{f_{\theta}}=\mu_0 \circ  f_{\theta}^{-1}$. For any $\theta \in \Theta$, assume that for all $n = 1,2,\ldots,N$ the following conditions hold:
\begin{itemize}
\item The space $\text{Im}(\mathcal{F}_{\theta_n}^{(n)}) \subset \mathcal{H}$, where $\text{Im}(\mathcal{F}_{\theta_n}^{(n)})$ denotes the image of $\mathcal{F}_{\theta_n}^{(n)}$.
\item The operator $\mathcal{F}_{\theta_n}^{(n)}$ has finite dimensional range.
\item The operator $f_{\theta_n}^{(n)}$ is bijective.
\item For any $u \in \mathcal{H}_u$,  all point spectrum of $D\mathcal{F}_{\theta_n}^{(n)}(u)$ are not in $(-\infty,-1]$.
\end{itemize}
Then $\mu_{f_{\theta}}$ is equivelent with $\mu_0$, and the RN derivative of $\mu_{f_{\theta}}$ with respect to $\mu_0$ is given by:
\begin{align*}
    \frac{d\mu_{f_{\theta}}}{d\mu_0}(f_{\theta}(u))=\prod\limits_{n=1}^{N}\left | {\det}_{1}(Df_{\theta_n}^{(n)}(u_{n-1}))   \right |^{-1} \exp\left(\frac{1}{2}\langle f_{\theta}(u)-u,f_{\theta}(u)-u\rangle_{\mathcal{H}}+  \langle u,u-f_{\theta}(u)\rangle_{\mathcal{H}}\right),
\end{align*}
where $u_n=f_{\theta_n}^{(n)} \circ f_{\theta_{n-1}}^{(n-1)} \circ \cdots \circ f_{\theta_1}^{(1)}(u)$ for each  $n = 1, 2, \ldots ,  N$, and ${\det}_{1}(Df_{\theta_n}^{(n)}(u_{n-1}))$ is the Fredholm-Carleman determinant \cite{gaussianmeasure} of the linear operator $Df_{\theta_n}^{(n)}(u_{n-1})$ (A rigorous definition can be found in the Appendix).
\end{theorem}
}

Theorem \ref{thm:RN} establishes the equivalence between transformed measures $\mu_{f_{\theta}}$ and the prior $\mu_0$. Consequently, the collection of transformed measures, denoted $\mathcal{M}(\mathcal{H}_u)$, constitutes a suitable set of approximating measures for variational inference. Leveraging Theorem \ref{thm:RN}, we can compute the necessary KL divergence
\begin{align*}
    D_{\text{KL}}(\mu_{f_{\theta}}||\mu)&=\int_{\mathcal{H}_u}\ln\left(\frac{d\mu_{f_{\theta}}}{d\mu}(u)\right)d\mu_{f_{\theta}}(u)\\
    &=\int_{\mathcal{H}_u}\ln\left(\frac{d\mu_{f_{\theta}}}{d\mu_0}(u)\right)d\mu_{f_{\theta}}(u)-\int_{\mathcal{H}_u}\ln\left(\frac{d\mu}{d\mu_0}(u)\right)d\mu_{f_{\theta}}(u)\\
    &=E_{\mu_{f_{\theta}}}\ln\left(\frac{d\mu_{f_{\theta}}}{d\mu_0}(u)\right)-E_{\mu_{f_{\theta}}}\ln\left(\frac{d\mu}{d\mu_0}(u)\right),
\end{align*}
where $E_{\mu_{f_{\theta}}}$ and $E_{\mu_0}$ denote the expectation with respect to $\mu_{f_{\theta}}$ and $\mu_0$, respectively.  Furthermore, 
\begin{align*}
    &E_{\mu_{f_{\theta}}}\ln\left(\frac{d\mu_{f_{\theta}}}{d\mu_0}(u)\right)=E_{\mu_0}\ln\left(\frac{d\mu_{f_{\theta}}}{d\mu_0}(f_{\theta}(u))\right)\\
    &=E_{\mu_0}\ln\left(\prod\limits_{n=1}^{N}\left | {\det}_{1}(Df_{\theta_n}^{(n)}(u_{n-1}))   \right |^{-1} \exp\left(\frac{1}{2}\langle f_{\theta}(u)-u,f_{\theta}(u)-u\rangle_{\mathcal{H}}+  \langle u,u-f_{\theta}(u)\rangle_{\mathcal{H}}\right)\right)\\
    &= E_{\mu_0}\left(-\sum\limits_{n=1}^{N}\ln\left(\left | {\det}_{1}(Df_{\theta_n}^{(n)}(u_{n-1}))   \right |\right)+\frac{1}{2}\langle f_{\theta}(u)-u,f_{\theta}(u)-u\rangle_{\mathcal{H}}+  \langle u,u-f_{\theta}(u)\rangle_{\mathcal{H}}\right).
\end{align*}
The optimization problem (\ref{equ:KLmin}) thus yields the following reformulation: 
\begin{align*}
\theta^*=\mathop{\arg\min}\limits_{\theta\in \Theta}D_{\text{KL}}(\mu_{{f}_{\theta}}||\mu).
\end{align*}
Monte Carlo sampling allows us to estimate the KL divergence, which serves as the loss function, to optimize the transformation parameters \(\theta\). Details of the algorithm are presented in Algorithm \ref{alg A}. The estimated gradient of $\theta_k$ in Algorithm \ref{alg A} can be computed by the preconditioned stochastic gradient-based optimization methods such as RMSprop (Root Mean Square Propagation) \cite{RMS}, AdaGrad (Adaptive Gradient) \cite{Adagrad}, and Adam (Adaptive Moment Estimation) \cite{Adam}. In this article, we adopt the Adam optimizer, with implementation details provided in Section \ref{sec4}. 

\begin{algorithm}[ht!]
\caption{Functional Normalizing Flow}
\label{alg A}
\begin{algorithmic}[1]
\STATE{Initialize the hyperparameter $\lambda$ to $\lambda_0$, parameter $\theta$ to $\theta_0$, set the sampling size to $N$, the iterative number to $K$, and the learning rate schedule to $\alpha_k$. Then, initialize the iteration counter $k$ to $0$;}
\REPEAT
\STATE{Sampling $N$ functions $\{u_i\}_{i=1}^{N}$ from prior $\mu_0$;}
\STATE{Update the parameters $\theta_{k+1}=\theta_{k}-\alpha_k\nabla_{\theta_k}L(\theta_k)$ with\\$\nabla_{\theta_k} L(\theta_k)\approx \frac{1}{N}\sum\limits_{i=1}^{N}\nabla_{\theta_k}\ln\left(\frac{d\mu_{f_{\theta_k}}}{d\mu_0}\left(f_{\theta_k}(u_i)\right)\right)-\frac{1}{N}\sum\limits_{i=1}^N\nabla_{\theta_k}\ln\left(\frac{d\mu}{d\mu_0}\left(f_{\theta_k}(u_i)\right)\right).$\\Variants of stochastic gradient-based
optimization method can be employed in this step;}
\STATE{$k=k+1$;}
\UNTIL{$k = K$.}
\STATE{Return the final result $\theta=\theta_K$.}
\end{algorithmic}
\end{algorithm}
{\color{black}
\begin{remark}
The loss function involves terms expressed in the $\mathcal{H}$-norm, such as $\langle f_{\theta}(u)-u, f_{\theta}(u)-u \rangle_{\mathcal{H}}$ and $\langle u,u-f_{\theta}(u) \rangle_{\mathcal{H}}$. Note that the $\mathcal{H}$-inner product is defined by 
$$\langle u,v\rangle_{\mathcal{H}}=\langle \mathcal{C}_0^{-\frac{1}{2}}u, \mathcal{C}_0^{-\frac{1}{2}}v \rangle_{{H}_u},$$ 
where $\mathcal{C}_0$ is a trace class operator. Without appropriate regularity conditions on the arguments of the $\mathcal{H}$-norm terms, these terms are likely to be infinite. However, to satisfy the measure equivalence requirement of Theorem~\ref{thm:RN} between the approximate and prior measures, the arguments must possess sufficient regularity, which in turn necessitates the use of the $\mathcal{H}$-norm. In fact, the equivalence of the post-transformed measure $\mu_{f_{\theta}}$ with the prior $\mu_0$, combined with established results from \cite{stuart2010inverse, Cotter2009IP}, typically ensures that the approximate measure is equivalent to the posterior. This equivalence is critical because it ensures the KL divergence (the loss function) is well-defined and finite.
\end{remark}
}

\subsection{Linear Transformation}\label{subsec:linear transformation}
Subsection \ref{subsec:normalizing flow in function space} presents a general theoretical framework that outlines the conditions for flow models in function space to be well-defined. In this subsection, we introduce two specific linear flow models and provide rigorous proofs that they satisfy these theoretical conditions. Theorem \ref{thm:RN} requires $f_{\theta_n}^{(n)}$ to be bijective for all $n$ and $\theta_n$. To ensure this property for linear flow, Lemma \ref{lem:libijective} is introduced.
\begin{lemma}\label{lem:libijective}
Let $\mathcal{H}_u$ be a separable Hilbert space, $f_{\theta_n}^{(n)}:\mathcal{H}_u \rightarrow \mathcal{H}_u$, and $f_{\theta_n}^{(n)}(u)=u+\mathcal{F}_{\theta_n}^{(n)}(u)$, where $\mathcal{F}_{\theta_n}^{(n)}$ is a bounded linear operator. 
Assume that the following conditions hold for all $n = 1, 2, \ldots , N$:
\begin{itemize}
\item The operator $\mathcal{F}_{\theta_n}^{(n)}$ has finite dimensional range.
\item For any $u \in \mathcal{H}_u$, all point spectrum of $D\mathcal{F}_{\theta_n}^{(n)}(u)$ are not in $(-\infty,-1]$.
\end{itemize}
Then the operator $f_{\theta_n}^{(n)}$ is bijective.
\end{lemma}
\subsubsection{Functional Householder Flow}
The first model, termed functional Householder flow, is a generalization of the Householder flow in Euclidean space \cite{householderflow1,Householderflow2}. Each layer of functional Householder flow employs the following transformation:
\begin{align*}
    f_{\theta_n}^{(n)}(u)=u-0.5v_n(\langle v_n,u\rangle_{\mathcal{H}_u}+b_n),
\end{align*}
where $v_n \in \mathcal{H}_u$,  $b_n \in \mathbb{R}$, and the parameters set is given by  $\theta_n=\{v_n, b_n\}$. We can ensure that the model meets the conditions in Theorem \ref{thm:RN} by setting some restrictions.
{\color{black}
\begin{theorem}\label{thm:Householder}
Let $ f_{\theta_n}^{(n)}(u)=u-0.5v_n(\langle v_n,u\rangle_{\mathcal{H}_u}+b_n)$. Assume that the following conditions hold for all $n = 1, 2, \ldots , N$:
\begin{itemize}
\item  $v_n \in \mathcal{H}$.
\item $\langle v_n,v_n\rangle_{\mathcal{H}_u}=1$.
\end{itemize}
Then $ f_{\theta_n}^{(n)}$ satisfies all the conditions of Theorem \ref{thm:RN}
\end{theorem}
}
The parameters \(\{v_n\}_{n=1}^{N}\) of the flow model, being functions themselves, require an efficient parameterization. This is essential not only for computational efficiency but also for ensuring compliance with the constraints outlined in the theorem. A novel parameterization method for the flow model is proposed here.
We denote $\{\phi_i\}_{i=1}^M$ as the first $M$ eigenfunctions of the prior covariance operator. The parameterization process for $v_n$ involves two steps.  Firstly, an auxiliary function $\hat{v}_n$ is introduced and parameterized as 
\begin{align}\label{equ:vn}
    \hat{v}_n=\sum\limits_{i=1}^{M}\alpha_i^n \phi_i.
\end{align}
Secondly, to ensure that the normalization is met, that is, \(\langle v_n, v_n \rangle_{\mathcal{H}_u} = 1\), we parameterize \(v_n\) using \(\hat{v}_n\) as shown in equation (\ref{equ:vn}):
\begin{align*}
v_n=\frac{\hat{v}_n}{\lVert\hat{v}_n\rVert_{\mathcal{H}_u}}.
\end{align*}

As the image of $\mathcal{F}_{\theta_n}^{(n)}$ is one-dimensional, the model’s expressive power is limited. To accurately approximate target distributions that are far from the prior, a deep network with multiple layers is often necessary. Detailed numerical results are presented in Section \ref{sec4}.
{\color{black}
\begin{remark}
    Similar with article \cite{chen2024learning}, we first require the eigenfunctions of the covariance operator $\mathcal{C}_0$ before running the algorithm. While these eigenfunctions can sometimes be derived analytically, numerical methods are necessary when analytical solutions are infeasible. For instance, we can discretize the operator $\mathcal{C}_0$ using the method in \cite{bui2013computational} and numerically obtain the eigenfunctions at some discrete level, employing techniques such as the double-pass algorithm \cite{saibaba2016randomized} or any other classical algorithm presented in \cite{golub2013matrix}.
\end{remark}
}

\subsubsection{Functional Projected Transformation Flow}
The second model is termed functional projected transformation flow. Let $\{\phi_i\}_{i=1}^{M}$ be the first $M$ eigenfunctions of the prior covariance operator, and we define two operators $\mathcal{P}$, $\mathcal{Q}$:
\begin{align}\label{equ:P}
    \mathcal{P}u=(\langle u,\phi_1 \rangle_{\mathcal{H}_u}, \langle u,\phi_2 \rangle_{\mathcal{H}_u},\ldots,\langle u,\phi_M \rangle_{\mathcal{H}_u})^{T},
\end{align}
\begin{align}\label{equ:Q}
\mathcal{Q}d=\sum\limits_{i=1}^{M}d_i \phi_i.
\end{align}
The transformation applied to each layer is defined as follows:
\begin{align*}
f_{\theta_{n}}^{(n)}(u)=u+\mathcal{Q}R_n(\mathcal{P}u+b_n),
\end{align*}
where $\mathcal{P}$ is defined in (\ref{equ:P}), $\mathcal{Q}$ is defined in (\ref{equ:Q}), 
$R_n$ is a $M \times M$ matrix, $b_n \in \mathbb{R}^M$, and the parameters set is given by $\theta_n=\{R_n, b_n\}$. The following theorem guarantees the model's adherence to the theoretical framework proposed in Subsection \ref{subsec:normalizing flow in function space}.
{\color{black}
\begin{theorem}\label{thm:pro}
Let $ f_{\theta_{n}}^{(n)}(u)=u+\mathcal{Q}R_n(\mathcal{P}u+b_n)$. Assume that the following conditions hold for all $n = 1, 2, \ldots , N$:
\begin{itemize}
\item The space $\text{Im}(\mathcal{Q})\subset\mathcal{H}$.
\item  All the point spectrum of $R_n$ are not in $(-\infty,-1]$.
\end{itemize}
Then $ f_{\theta_n}^{(n)}$ satisfies all the conditions of Theorem \ref{thm:RN}.
\end{theorem}
}
Note that the image of $\mathcal{F}_{\theta_n}^{(n)}$
is a $M$-dimensional space, functional projected transformation flow may exhibit greater expressive capacity than functional Householder flow. Detailed numerical results can be found in Section \ref{sec4}.

\subsection{Nonlinear Transformation}\label{subsec:nonleaner transformation}
In this subsection, we will explore two different nonlinear flow models. Theorem \ref{thm:RN} requires that the function $f_{\theta_n}^{(n)}$ be invertible for all values of $n$ and $\theta_n$. To guarantee this property for our nonlinear flow, we introduce Lemma \ref{lem:bijective}.
\begin{lemma}\label{lem:bijective}
Let $\mathcal{H}_u$ be a separable Hilbert space, $f_{\theta_n}^{(n)}:\mathcal{H}_u \rightarrow \mathcal{H}_u$, and $f_{\theta_n}^{(n)}(u)=u+\mathcal{F}_{\theta_n}^{(n)}(u)$. Assume that the following conditions hold for all $n = 1, 2, \ldots , N$:
\begin{itemize}
\item The operator $\mathcal{F}_{\theta_n}^{(n)}$ has finite dimensional range.
\item $\mathcal{F}_{\theta_n}^{(n)}(\mathcal{H}_u)$ is a bounded subset in $\mathcal{H}_u$.
\item For any $u \in \mathcal{H}_u$, all point spectrum of $D\mathcal{F}_{\theta_n}^{(n)}(u)$ are not in $(-\infty,-1]$.
\end{itemize}
Then the operator $f_{\theta_n}^{(n)}$ is bijective.
\end{lemma}
\subsubsection{Functional Planar Flow}
The first nonlinear model, termed functional planar flow, is a generalization of the planar flow in Euclidean space \cite{planarflow}. Each layer of functional planar flow employs the transformation:
\begin{align}\label{equ:functional planar}
    f_{\theta_n}^{(n)}(u)=u+u_nh(\langle w_n,u\rangle_{\mathcal{H}_u}+b_n),
\end{align}
where $u_n,w_n \in \mathcal{H}_u$, $b_n \in \mathbb{R}$, $h(x)=\tanh(x)$ and the parameters set is given by  $\theta_n=\{u_n, w_n, b_n\}$. 
The following theorem guarantees the model's adherence to the theoretical framework proposed in Subsection \ref{subsec:normalizing flow in function space}.
{\color{black}
\begin{theorem}\label{thm:planar}
Let $f_{\theta_n}^{(n)}(u)=u+u_nh(\langle w_n,u\rangle_{\mathcal{H}_u}+b_n)$. Assume that the following conditions hold for all $n = 1, 2, \ldots , N$:
\begin{itemize}
\item $w_n, u_n \in \mathcal{H}$.
\item $\langle u_n,w_n\rangle_{\mathcal{H}_u}>-1$.
\end{itemize}
Then $ f_{\theta_n}^{(n)}$ satisfies all the conditions of Theorem \ref{thm:RN}.
\end{theorem}
}
Given that the parameters $u_n$ and $w_n$ of our flow model are functions, an efficient parameterization is essential both for computational efficiency and for ensuring compliance with the assumptions of the theorem. To address this, we propose a novel parameterization method for the flow model. Let $\{\phi_i\}_{i=1}^{M}$ represent the first $M$ eigenfunctions of the prior covariance operator. Our parameterization of $u_n$ and $w_n$ involves a two-step process. Firstly, auxiliary functions $\hat{u}_n$, $\hat{w}_n$ are introduced and parameterized as
\begin{align}\label{equ:un}
    \hat{u}_n=\sum\limits_{i=1}^{M}\alpha_i^n \phi_i, \qquad
    \hat{w}_n=\sum\limits_{i=1}^{M}\beta_i^n \phi_i.
\end{align}
Secondly, considering the conditions of Theorem \ref{thm:planar}, i.e. $\langle u_n,w_n \rangle_{\mathcal{H}_u} >- 1$, we parameterize $u_n$ and $w_n$ using $\hat{u}_n$ and $\hat{w}_n$ as follows:
\begin{align}\label{equ:wnchange}
w_n=\hat{w}_n,
\end{align}
\begin{align}\label{equ:unchange}
u_n=\hat{u}_n + [q(\langle \hat{u}_n,\hat{w}_n \rangle_{\mathcal{H}_u})-\langle \hat{u}_n,\hat{w}_n \rangle_{\mathcal{H}_u}]\frac{\hat{w}_n}{\lVert\hat{w}_n\rVert_{\mathcal{H}_u}},
\end{align}
where $q(x)=\ln(1+e^x)-1$.

Similar with functional Householder flow, the image of $\mathcal{F}_{\theta_n}^{(n)}$ is constrained to a one-dimensional space, which limits the model's expressive capabilities. To effectively approximate target measures that are significantly different from the prior, a deep network architecture with multiple layers is often necessary. The effectiveness of this approach is demonstrated through numerical experiments in Section \ref{sec4}.

\subsubsection{Functional Sylvester Flow}
The second nonlinear model we are exploring, known as functional Sylvester flow, is an extension of Sylvester flow, a concept traditionally applied in Euclidean space \cite{sylvesterflow}. Each layer of functional Sylvester flow utilizes the following transformation:
\begin{align*}
f_{\theta_n}^{(n)}(u)=u+\mathcal{A}_nh(\mathcal{B}_nu+b_n),
\end{align*}
where $\mathcal{B}_n$ is a bounded linear operator mapping from $\mathcal{H}_u$ to $\mathbb{R}^M$, $\mathcal{A}_n$ is a bounded linear operator mapping from $\mathbb{R}^M$ to $\mathcal{H}_u$, $b_n \in \mathbb{R}^M$, $h(x)=\tanh(x)$, and the parameters set is given by $\theta_n=\{\mathcal{A}_n, \mathcal{B}_n, b_n\}$. The following theorem ensures that this model adheres to the theoretical framework we outlined in Subsection \ref{subsec:normalizing flow in function space}.

{\color{black}
\begin{theorem}\label{thm:Sylvester}
Let  $ f_{\theta_n}^{(n)}(u)=u+\mathcal{A}_nh(\mathcal{B}_nu+b_n)$. Assume that the following conditions hold for all $n = 1, 2, \ldots , N$:
\begin{itemize}
\item The space $\text{Im}(\mathcal{A}_n)\subset \mathcal{H}$.
\item The operator $\mathcal{A}_n$ has finite dimensional range.
\item All the point spectrum of $\mathcal{B}_n\mathcal{A}_n$ are not in $(-\infty,-1]$.
\end{itemize}
Then $ f_{\theta_n}^{(n)}$ satisfies all the conditions of Theorem \ref{thm:RN}.
\end{theorem}
}
Given that the parameters $\mathcal{A}_n$ and $\mathcal{B}_n$ of the flow model are operators, efficient parameterization is essential for both computational efficiency and ensuring compliance with the theorem's constraints. We select the first $M$ eigenfunctions of the prior covariance operator, denoted as $\{\phi_i\}_{i=1}^{M}$. The parameterization of $\mathcal{B}_n$ involves two steps:
\begin{itemize}
\item First, $\mathcal{B}_n^{(1)}$ projects $u$ onto a finite-dimensional subspace, defining $p_{\mathcal{B}}$:
\begin{align*}
 p_{\mathcal{B}}=\mathcal{B}_n^{(1)}u:=
 \begin{pmatrix}
 \langle u,\phi_1 \rangle_{\mathcal{H}_u}\\
 \langle u,\phi_2 \rangle_{\mathcal{H}_u}\\
 \vdots\\
 \langle u,\phi_M \rangle_{\mathcal{H}_u}
 \end{pmatrix}.
\end{align*}
\item Second, $\mathcal{B}_n^{(2)}$ applies a linear transformation on $p_{\mathcal{B}}$ via the square matrix $R_{\mathcal{B}}^n$:
\begin{align}\label{equ:R1}
\mathcal{B}_n^{(2)}p_{\mathcal{B}}:=R_{\mathcal{B}}^np_{\mathcal{B}},
\end{align}
where $R_{\mathcal{B}}^n$ is a square matrix that maps $p_{\mathcal{B}}$ to the target.
\end{itemize}
The overall operator $\mathcal{B}_n$ for $u\in \mathcal{H}_u$ is defined as $\mathcal{B}_n= \mathcal{B}_n^{(2)}\mathcal{B}_n^{(1)}u$, representing the composition of the two aforementioned steps.

The parameterization of $\mathcal{A}_n$ is achieved through the following two steps: 
\begin{itemize}
\item The first step, $\mathcal{A}_n^{(1)}$ maps the vector $d$ to $p_{\mathcal{A}} \in \mathbb{R}^{M}$ using the square matrix $R_{\mathcal{A}}^n$:
\begin{align}\label{equ:R2}
 p_{\mathcal{A}}=\mathcal{A}_n^{(1)}d:=R_{\mathcal{A}}^nd.
\end{align}
\item The second step, $\mathcal{A}_n^{(2)}$ maps the vector $p_{\mathcal{A}}$ back to the function space spanned by the basis functions $\{\phi_i\}$:
\begin{align*}
\mathcal{A}_n^{(2)}p_{\mathcal{A}}:=\sum\limits_{i=1}^{M}p_{\mathcal{A}}^{i}\phi_{i},
\end{align*}
where $p_{\mathcal{A}}^{i}$ is the value of the $i\text{-}th$ position of the vector $p_{\mathcal{A}}$.
\end{itemize}
The overall operator $\mathcal{A}_n$ is defined for $d\in \mathbb{R}^M$ as  $\mathcal{A}_n= \mathcal{A}_n^{(2)}\mathcal{A}_n^{(1)}d$, representing the composition of the two aforementioned steps.

To satisfy the conditions imposed by Theorem \ref{thm:Sylvester}, we establish the following result.

\begin{theorem}\label{thm:change}
Let $\mathcal{H}_u$ be a separable Hilbert space. Let $\mathcal{B} : \mathcal{H}_u \to \mathbb{R}^M$ and $\mathcal{A} : \mathbb{R}^M \to \mathcal{H}_u$ be bounded linear operators. Then the following statements hold:
\begin{itemize}
\item The point spectrum of $I+\mathcal{A}\mathcal{B}$ corresponds one-to-one with the eigenvalues of $I_M+\mathcal{B}\mathcal{A}$.
\item The Fredholm determinant of $I+\mathcal{A}\mathcal{B}$ is equal to the determinant of $I_M+\mathcal{B}\mathcal{A}$.
\end{itemize}
\end{theorem}
Based on Theorem \ref{thm:change}, we have the following equality:
\begin{align*}
{\det}_{1}(I+\mathcal{A}_n\mathrm{diag}(h'(\mathcal{B}_nu+b_n))\mathcal{B}_n)=\det(I_M+\mathrm{diag}(h'(\mathcal{B}_nu+b_n))\mathcal{B}_n\mathcal{A}_n).
\end{align*}
This allows us to efficiently compute the Fredholm determinant of the operator $Df_{\theta_n}^{(n)}(u) = I + \mathcal{A}_n\mathrm{diag}(h'(\mathcal{B}_nu+b_n))\mathcal{B}_n$. 
Theorem \ref{thm:change} establishes the equivalence between the point spectrum of  
$I + \mathcal{A}_n\mathrm{diag}(h'(\mathcal{B}_nu+b_n))\mathcal{B}_n$ and the eigenvalues of the matrix $I_M + \mathrm{diag}(h'(\mathcal{B}_nu+b_n))\mathcal{B}_n\mathcal{A}_n$.  In order to ensure that the eigenvalues of $\mathrm{diag}(h'(\mathcal{B}_nu+b_n))\mathcal{B}_n\mathcal{A}_n$ are not in $(-\infty,-1]$, we set $R_{\mathcal{B}}^n$ in (\ref{equ:R1}) as an upper triangular matrix, $R_{\mathcal{A}}^n$ in (\ref{equ:R2}) as a lower triangular matrix, the diagonal elements of $R_{\mathcal{B}}^n$ are all 1, and the diagonal elements of $R_{\mathcal{A}}^n$ are all greater than $-1$. Through this method, $\mathrm{diag}(h'(\mathcal{B}_nu+b_n))\mathcal{B}_n\mathcal{A}_n$ is an upper triangular matrix, and the diagonal elements are all greater than $-1$. Hence we can compute the determinant of $I_M + \mathrm{diag}(h'(\mathcal{B}_nu+b_n))\mathcal{B}_n\mathcal{A}_n$ quickly.

Since the image of $\mathcal{F}_{\theta_n}^{(n)}$ is a $M$-dimensional space, functional Sylvester flow may exhibits a higher degree of expressive power compared to functional planar flow. We will delve into detailed numerical results in Section \ref{sec4}.

\subsection{Discretization Invariance}\label{subsec:discrete invariance}
The model, defined within a function space, inherently implies discretization invariance. In general, it can be summarized as the following three properties \cite{kovachki2023neural}:
\begin{itemize}
\item It operates on any discretization of the input function, meaning it can accept any set of points within the input domain.
\item Its output can be evaluated at any point in the output domain.
\item As the discretization becomes finer, the model converges to a continuous operator.
\end{itemize}

In this subsection, we investigate the discretization invariance of the operator $\mathcal{F}_{\theta_n}^{(n)}$ within the flow model $f_{\theta_n}^{(n)} = I + \mathcal{F}_{\theta_n}^{(n)}$. Following the approach outlined in \cite{kovachki2023neural}, a formal definition of the discretization invariance is provided. 

\begin{definition}
We call a discrete refinement of the domain $D \subset \mathbb{R}^d$ any sequence of nested sets $D_1 \subset D_2 \subset \cdots \subset D$ with $|D_L|=L$ for any $L \in \mathbb{N}$ such that, for any $\epsilon >0$, there exists a number $L=L(\epsilon) \in \mathbb{N}$ such that
\begin{align*}
D \subset \bigcup_{x \in D_L}\{ y \in \mathbb{R}^d: \lVert y-x\rVert < \epsilon \}.
\end{align*}
\end{definition}
\begin{definition}
Given a discrete refinement ${\{D_L\}}_{L=1}^{\infty}$ of the domain $D \subset \mathbb{R}^d$, any member $D_L$ is
called a discretization of $D$.
\end{definition}
\begin{definition}
Suppose $\mathcal{H}_u$ is a Hilbert space of $\mathbb{R}^m$ valued functions on the domain $D \subset \mathbb{R}^d$. For any $\theta_n$, let $\mathcal{F}_{\theta_n}^{(n)}: \mathcal{H}_u \rightarrow \mathcal{H}_u$ be an operator, $D_L$ be an $L$-point discretization of $D$, and operator $\hat{\mathcal{F}}_{\theta_n(L)}^{(n)}:  \mathbb{R}^{Lm} \rightarrow \mathcal{H}_u$. For any $K \subset C(D)$ compact, we define the discretized uniform risk as
\begin{align*}
R_K(\mathcal{F}_{\theta_n}^{(n)},\hat{\mathcal{F}}_{\theta_n(L)}^{(n)})=\sup_{a \in K}\left\lVert\hat{\mathcal{F}}_{\theta_n(L)}^{(n)}(a|_{D_L})-\mathcal{F}_{\theta_n}^{(n)}(a)\right\rVert_{\mathcal{H}_u}.
\end{align*}
\end{definition}

\begin{definition}
Given $a$ discrete refnement ${\{D_n\}}_{n= 1}^\infty$ of the domain 
 $D\subset \mathbb{R} ^d$. For any fixed $\theta_n$, we say $\mathcal{F}_{\theta_n}^{(n)}$ is discretization-invariant if there exists a sequence of maps $\hat{\mathcal{F}}_{\theta_n(1)}^{(n)},\hat{\mathcal{F}}_{\theta_n(2)}^{(n)},\ldots$ where $\hat{\mathcal{F}}_{\theta_n(L)}^{(n)}:\mathbb{R}^{Lm}\to\mathcal{H}_u$ such that, for any compact set
$K \subset C(D)$, 
\begin{align*}
\lim_{L\to\infty}R_K\left(\mathcal{F}_{\theta_n}^{(n)},\hat{\mathcal{F}}_{\theta_n(L)}^{(n)}\right)=0.
\end{align*}
\end{definition}

Under the proposed definition, the aforementioned four flow models can be shown to exhibit discretization invariance. 
\begin{theorem}\label{thm:discrete invariance}
Let $D\subset\mathbb{R}^d$ be a domain for some $d \in \mathbb{N}$. Suppose that $\mathcal{H}_u$ can be continuously embedded in $C(D)$. Then for any $ n\in\mathbb{N}$, the layers $\mathcal{F}_{\theta_n}^{(n)}: \mathcal{H}_u \rightarrow \mathcal{H}_u$ of functional planar flow, functional Householder flow, functional Sylvester flow, and functional projected transformation flow are discretization-invariant.
\end{theorem}
Next, we compare and contrast the functional normalizing flow framework with standard normalizing flow in Euclidean space. A key property contributing to the success of normalizing flow in computer vision is the tractability of the density of the post-transformed distribution. To achieve this, the transformation $\widetilde{f}(z)$ of normalizing flow in Euclidean space must be invertible \cite{planarflow,nvp}. In this case, the post-transformed density can be expressed as
\begin{align*}
p_1(z') = p_0(z) \left| \det \left(\frac{\partial \widetilde{f}(z)}{\partial z}\right) \right|^{-1},
\end{align*}
where $z'=\widetilde{f}(z)$, $z, z'\in \mathbb{R}^{L}$ and $\widetilde{f}$ maps from $\mathbb{R}^L$ to $\mathbb{R}^L$, the densities of the post-transformed and pre-transformed distributions are denoted by $p_1(z)$ and $p_0(z)$, respectively. This tractable form allows for the computation of the Kullback-Leibler divergence between the target distribution and the post-transformed distribution, which is used as the loss function. 
 
In this paper, we demonstrate the applicability of similar ideas within the infinite-dimensional Hilbert space $\mathcal{H}_u$, aligning well with the Bayesianize-then-Discretize approach for inverse problem. Similar to the normalizing flow in Euclidean space, ensuring the invertibility of the transformation $f_{\theta}$ is crucial for the tractability of the RN derivative between the post-transformed and pre-transformed measures \cite{gaussianmeasure}. While various different approaches have been proposed to ensure the invertibility of normalizing flow transformations in Euclidean space \cite{planarflow,nvp}, we present general theorems (Lemmas \ref{lem:libijective} and \ref{lem:bijective}) that establish conditions for the invertibility of transformations within the infinite-dimensional function space.  This ensures the tractability of the RN derivative between the post-transformed measure and the prior, enabling the calculation of the KL divergence as the loss function.

Lastly, we demonstrate that, by picking a parameterization that is inconsistent in function space, the model we proposed will degenerate to the standard normalizing flow in Euclidean space. For illustrative purposes, we will consider a simple example using functional planar flow (\ref{equ:functional planar}). Let $u$ be a real-valued function defined on the domain $D \subset \mathbb{R}^m$, $\mathcal{H}_u=L^2(D)$,  $x^1, x^2, \ldots , x^L \in D$ be the points at which the input function $u$ is evaluated and denote $\bm{u} = (u(x^1), u(x^2), \ldots , u(x^L)) \in R^L$ the vector of evaluations. Let $w_n(x^i)=w_n^{(i)}$, $u_n(x^i) = u_n^{(i)}$ for $i = 1, \ldots, L$ where $w_n^{(i)}, u_n^{(i)} \in \mathbb{R}$ are some constants. Employing a Monte Carlo approximation, the functional planar flow can be discretized as

\begin{align}\label{equ:Euclid planar}
\widetilde{f}_{\theta_n}^{(n)}(\bm{u})=\bm{u}+\bm{u}_n h\left(\frac{1}{L}\sum\limits_{i=1}^{L}w_n^{(i)}u(x^{i})+b_n\right).
\end{align}
This expression can be simplified to
\begin{align}\label{equ:Euclid planar1}
\widetilde{f}_{\theta_n}^{(n)}(\bm{u})=\bm{u}+\bm{u}_n h(\langle\bm{w}_n, \bm{u}\rangle_{\mathbb{R}^L}+b_n),
\end{align}
where $\bm{u}_n=(u_n^{(1)}, \ldots ,  u_n^{(L)})$, $\bm{w}_n=(\frac{w_n^{(1)}}{L}, \ldots ,  \frac{w_n^{(L)}}{L})$, $b_n \in \mathbb{R}$ is constant,  $h(x)=\tanh(x)$. Parameters $\theta_n$ can be regard as $\theta_n = \{\bm{u}_n, \bm{w}_n, b_n\}$. The construction (\ref{equ:Euclid planar}) can be generalized to arbitrary number of layers by composing several transformations together
\begin{align}\label{equ:planar}
    \widetilde{f}_{\theta}(\widetilde{u})=\widetilde{f}_{\theta_{N}}^{(N)} \circ \widetilde{f}_{\theta_{N-1}}^{(N-1)} \circ\cdots \circ \widetilde{f}_{\theta_{1}}^{(1)}(\widetilde{u}).
\end{align}
The obtained model (\ref{equ:planar}) is the planar flow defined in the Euclidean space \cite{planarflow}. It is evident that when employing the planar flow (\ref{equ:Euclid planar1}) in function spaces, the parameterization of $w_n$ becomes contingent upon the specific discretization of the input $\bm{u}$. Likewise, the parameterization of $u_n$ depends on the desired discretization of the output. Consequently, for varying discretization requirements, the entire model necessitates retraining, precluding the ability to generalize to arbitrary discretized inputs and outputs. Consequently, standard planar flow in Euclidean space are not consistent within the infinite-dimensional function space, and it is not discretization-invariant.

\subsection{Connection with conventional neural operator}\label{subsec:neural_op_connection}
{\color{black}Our proposed functional normalizing flow shares a similar neural architecture to conventional neural operators \cite{kovachki2023neural, li2020fourier, anandkumar2020neural, kovachki2024operator}, as they both learn mappings between function spaces. However, the two models address fundamentally different problems. The neural operator aims to learn a deterministic one-to-one mapping $\mathcal{G}$ between function spaces using paired data $\{u_i, \mathcal{G}(u_i)\}_{i=1}^N$, minimizing a point-wise distance loss $L(\theta) = \sum_{i=1}^N d(\mathcal{G}(u_i), f_\theta(u_i))$. In contrast, functional normalizing flow seeks to construct a mapping $f_{\theta}$ that transforms a tractable source measure $\mu_0$ into a complex target measure $\mu$, requiring no paired data. Its objective is to ensure that the transformed ensembles $\{f_{\theta}(u_i)\}_{i=1}^N$ serve as approximate samples from $\mu$, which is achieved by minimizing the KL divergence $D_{\text{KL}}(\mu_{f_{\theta}}||\mu)$. 

Beyond these fundamental differences in objectives and loss functions, it is noteworthy that both functional normalizing flow and a wide range of neural operators exploit low-rank structures to enhance efficiency. As a representative and concise example, the low-rank neural operator (LNO) in Subsection 4.2 of \cite{kovachki2023neural} projects infinite-dimensional input functions onto a low-dimensional subspace via projection bases, and subsequently reconstructs the outputs back into the original function space using reconstruction bases. Specifically, consider a single layer of the LNO, an input function $v_t$ is first projected onto the spatial basis family $\{\psi_k\}_{k=1}^r$ via the $L^2$ inner product $\langle \psi_k, v_t \rangle$, compressing it into a finite-dimensional latent coefficient vector. These coefficients are then mapped back to the infinite-dimensional space using reconstruction bases $\{\varphi_k\}_{k=1}^r$, combined with a pointwise linear transformation. Consequently, a single layer of the LNO is defined by
\[
v_{t+1}(x)=\sigma\left(W_tv_t(x)+\sum_{k=1}^r \langle \psi_k, v_t \rangle \varphi_k(x)+b_t\right),
\]
where the trainable parameters consist of the basis families $\{\psi_k\}_{k=1}^r$ and $\{\varphi_k\}_{k=1}^r$, along with the weight matrix $W_t$ and bias $b_t$.
However, despite this architectural similarity, the role of low-rank structures within their respective theoretical frameworks is fundamentally distinct. In neural operators, low-rank structures are employed merely as an empirical technique to reduce computational costs and accelerate deterministic learning, without being a theoretical necessity. In contrast, for our functional normalizing flow model tailored for infinite-dimensional Bayesian variational inference, low-rank structures are theoretically necessitated. Since probability measures defined on infinite-dimensional function spaces are inherently prone to mutual singularity \cite{gaussianmeasure, MR2663405}, functional normalizing flow transformations must strictly preserve measure equivalence, as established in Theorem \ref{thm:RN}; otherwise, the KL-based loss function and the RN derivative would become inapplicable. Crucially, low-rank structures constitute an essential component of this theorem. Beyond the low-rank constraint, functional normalizing flow is subject to further theoretical requirements. To satisfy these conditions, our framework restricts functional normalizing flow architecture to layers exclusively based on a residual perturbation form \(f_n(u)=u+\mathcal{F}_n(u)\). We endow the network with sufficient complexity through \(\mathcal{F}_n\) while strictly satisfying all the constraints of Theorem \ref{thm:RN}.

Furthermore, another noteworthy point is the parameterization of the projection and reconstruction basis functions. In the context of neural operators, there are generally two primary approaches to parameterizing these basis functions. The first approach, exemplified by the Fourier neural operator (FNO) \cite{li2020fourier} and the principal component analysis neural operator (PCA-Net) \cite{bhattacharya2021model, kovachki2024operator}, utilizes a fixed set of bases as projection directions. The second approach, represented by the LNO \cite{kovachki2023neural}, employs trainable coordinate networks to adaptively learn the projection and reconstruction bases. While the first parameterization method yields high computational efficiency, its expressiveness may be constrained by the choice of basis. The second method allows for adaptive learning of bases and is theoretically more flexible. However, this generality entails a trade-off, as emphasized by the authors of \cite[p.~27]{kovachki2023neural}: ``the Low-rank
neural operator is finite-dimensionalized on a product space formed from the Barron space of neural
networks $\cdots$ The Barron space (Low-rank operator and DeepONet) representation is usually the most generic and flexible as it is widely applicable. However, this can lead to induced biases and reduced representation power on specific problems''.

In the numerical experiments of this paper, we adopt a parameterization method analogous to the first strategy. For instance, consider the functional planar flow \(f_n(u)=u+u_n h(\langle u,w_n \rangle_{\mathcal{H}_u}+b)\). Given a Gaussian prior \(\mu_0 = \mathcal{N}(0,\mathcal{C}_0)\), where the covariance operator \(\mathcal{C}_0\) admits eigen-pairs \(\{\lambda_i, \phi_i\}_{i=1}^{\infty}\), we parameterize the learnable components as \(u_n=\sum\limits_{i=1}^r \lambda_i\alpha_i\phi_i\) and \(w_n=\sum\limits_{i=1}^r \lambda_i\beta_i\phi_i\), where $\{\alpha_i,\beta_i\}_{i=1}^{r}$ are the trainable parameters. This utilization of the prior information from \(\mu_0\) ensures that the designed neural network satisfies the theoretical requirements established in Theorem \ref{thm:RN}. It is noteworthy that, conceptually, our functional normalizing flow architecture could also employ trainable networks to parameterize the low-rank bases. For instance, within the functional planar flow, \(u_n(x)\) and \(w_n(x)\) could be modeled as trainable neural networks. In this context, it is crucial to ensure that such parameterizations enable the resulting functional normalizing flow to satisfy the theoretical conditions of Theorem \ref{thm:RN}. In this work, we focus on the first parameterization for functional normalizing flow and do not examine the second in detail, leaving its thorough investigation to future work.
}

\section{Conditional Normalizing Flow in Function Space}\label{sec3}
While Section \ref{sec:functional normalizing flow} establishes the theoretical framework of functional normalizing flow, the practical application of Algorithm \ref{alg A} is hindered by the computational cost of repeatedly solving forward problems. To address this limitation, strategies for reducing the algorithm's iteration are essential.
In order to mitigate the computational burden, we propose conditional functional normalizing flow, which involve pre-training a model to map measurement data and their corresponding measurement points to the parameters of a functional normalizing flow, enabling efficient estimation of posterior for new measurement {\color{black} (similar techniques have also been employed for normalizing flows in Euclidean space \cite{ardizzone2019guided, winkler2019learning})}.

Our objective is to learn the map from measurement information (measurement data $\bm{d}$ and their corresponding measurement points $\bm{x}$) to the parameters of a functional normalizing flow, enabling approximation of posterior corresponding to the given measurement information.  Following model training, the conditional functional normalizing flow can efficiently generate approximate posterior for new measurement information. While this provides an initial estimate, further refinement can be achieved through fine-tuning the model on specific measurement information, leveraging the initial estimate as a starting point. This approach significantly reduces the number of iterations required compared with training from a random initialization.

When constructing conditional functional normalizing flow, the position of the measurement points $\bm{x}$ and the measurement data $\bm{d}$ both play important roles. 
Therefore, we need to reformulate the framework proposed in Section \ref{sec:functional normalizing flow}. 

Let $\bm{x}=(x^1, x^2, \ldots, x^{N_{\bm{d}}})$ represents the measurement points of the PDE solution $w$, $\bm{d}$ denotes the measurement data, and let $\mathcal{G}: \mathcal{H}_u \rightarrow \mathcal{H}_w$ represents the PDE solution operator, where $\mathcal{H}_u$ is the space of parameter, and $\mathcal{H}_w$ is the space of solution. The measurement points, denoted by $\bm{x}$, influence the measurement operator $\mathcal{S}_{\bm{x}}$:
\begin{align*}
    \mathcal{S}_{\bm{x}}: \mathcal{H}_w \rightarrow \mathbb{R}^{N_{\bm{d}}},
\end{align*}
where $N_{\bm{d}}$ is the number of observed points. Specifically, the operator $\mathcal{S}_{\bm{x}}$ acts as follows:
\begin{align*}
    \mathcal{S}_{\bm{x}}(w)=(w(x^1),w(x^2), \ldots , w(x^{N_{\bm{d}}})).
\end{align*}
The inverse problem can be formulated as
\begin{align*}
 \bm{d} = \mathcal{S}_{\bm{x}}\mathcal{G}(u) + \bm{\epsilon},
\end{align*}
where \(\bm{d} \in \mathbb{R}^{N_{\bm{d}}}\) is the measurement data, \(u \in \mathcal{H}_u\) is the parameter of interest, and \(\bm{\epsilon}\) is a Gaussian random vector with zero mean and covariance matrix \(\bm{\Gamma}_{\text{noise}} := \tau^{-1}\textbf{I}\) (\(\tau\) is a fixed positive number, \(\textbf{I}\) denotes the \(N_{\bm{d}}\)-dimensional identity matrix), which means
\begin{align*}
 \bm{\epsilon} \sim \mathcal{N}(0, \bm{\Gamma}_{\text{noise}}).
\end{align*}
Our objective is to infer the possible values of the model parameter $u$ based on the measurement data $\bm{d}$ and their corresponding measurement points $\bm{x}$ by the corresponding Bayes' formula:
\begin{align}\label{eq:post2}
	\qquad\,\,\frac{d\mu_{(\bm{x},\bm{d})}}{d\mu_0}(u) = \frac{1}{Z_{\mu}^{(\bm{x},\bm{d})}}\exp (-\Phi_{(\bm{x},\bm{d})}(u) ),
\end{align}
where $\mu_{(\bm{x},\bm{d})}$ is the posterior, and the potential function $\Phi_{(\bm{x},\bm{d})} : \mathcal{H}_u \rightarrow \mathbb{R}$ is defined as
\begin{align*}
 \Phi_{(\bm{x},\bm{d})}(u) = \frac{1}{2}\lVert \bm{d} - \mathcal{S}_{\bm{x}}\mathcal{G}(u)\rVert ^{2}_{\bm{\Gamma}_{\text{noise}}}
\end{align*}
with $\lVert  \cdot  
 \rVert_{\bm{\Gamma}_{\text{noise}}}:=\lVert \bm{\Gamma}_{\text{noise}}^{-1/2}   \cdot\rVert$, and $Z_{\mu}^{(\bm{x},\bm{d})}$ is the normalization constant  given by \begin{align*}
	Z_{\mu}^{(\bm{x},\bm{d})} = \int_{\mathcal{H}_u} \exp (-\Phi_{(\bm{x},\bm{d})}(u))\mu_0(du).
\end{align*}

The goal of conditional functional normalizing flows is to approximate the posterior distribution defined by Bayes’ formula, where the forward operator $\mathcal{G}$ is fixed while the measurement information $(\bm{x}, \bm{d})$ varies. By entering the measurement pair $(\bm{x}, \bm{d})$ into the neural network, we aim to directly construct an approximate probability measure $\mu_{(\bm{x}, \bm{d})}$ for the corresponding posterior distribution, as defined by Bayes’ formula (\ref{eq:post2}). It is important to note that the measurement information $(\bm{x},\bm{d})$ is both dynamic and heterogeneous. Variations in the number, positions, and values of measurement points lead to varying dimensions and scales of the information. To facilitate the application of neural networks, we therefore seek to transform this heterogeneous information into a fixed-length vector representation.

We develop a method to integrate two pieces of information, i.e., measurement data $\bm{d}$ and their corresponding locations $\bm{x}$. The adjoint operator of $\mathcal{S}_{\bm{x}}$,  denoted $\mathcal{S}_{\bm{x}}^*$, maps vector in $\mathbb{R}^{N_{\bm{d}}}$ to functional in the dual space $\mathcal{H}_w^*$. Formally,
 \begin{align*}
    \mathcal{S}_{\bm{x}}^*: \mathbb{R}^{N_{\bm{d}}} \rightarrow \mathcal{H}_w^*.
\end{align*}
The specific action of the adjoint operator is given by
\begin{align*}
     \mathcal{S}_{\bm{x}}^*{\bm{d}}(u)  = (\bm{d},\mathcal{S}_{\bm{x}}u)_{\mathbb{R}^{N_{\bm{d}}}} = \sum\limits_{i=1}^{N_{\bm{d}}}d_i u(x_i),
\end{align*}
where $d_i$ denotes the $i$-th value of the vector $\bm{d}$.
The adjoint operator $\mathcal{S}_{\bm{x}}^*$ can be used to integrate measurement data $\bm{d}$ and their corresponding measurement points $\bm{x}$ into a functional $\mathcal{S}_{\bm{x}}^{*}\bm{d}$. However, neural networks typically require vectors input. To address this issue, we propose calculating the values of the functional on a set of fixed functions, resulting in a vector $\bm{v}$ that can be directly input into the neural network. The entire process $\mathcal{W}(\bm{x},\bm{d})=\bm{v}$ can be divided into the following two steps:
\begin{itemize}
 	\item {\color{black} Basis selection: Choose a set of basis functions $\{\phi_i\}_{i=1}^{M}$ that span the relevant function space. Here we select the first $M$ eigenvectors of the prior covariance operator $\mathcal{C}_0$}.
	\item Functional evaluation: For each basis function $\phi_i$, compute $\mathcal{S}_{\bm{x}}^*{\bm{d}}(\phi_i)$, and collect these scalar values into a vector $\bm{v}=(\mathcal{S}_{\bm{x}}^*{\bm{d}}(\phi_1), \mathcal{S}_{\bm{x}}^*{\bm{d}}(\phi_2), \ldots , \mathcal{S}_{\bm{x}}^*{\bm{d}}(\phi_n))$. 
\end{itemize}

Under these circumstances, we can construct a neural network that maps the vector $\bm{v}$ to the parameters of the functional normalizing flow. Notably, the parameters of the functional normalizing flow are denoted as $\theta=\{\theta_1,\theta_2,\ldots,\theta_N\}$, where $\theta_n$ corresponds to the parameters of a transformation $f_{\theta_n}^{(n)}$ within the model for $n = 1, 2, \ldots, N$. We construct a neural network $\mathcal{N}_{\lambda_n}^{(n)}(\bm{v})$ for each transformation, mapping the vector $\bm{v}$ to the parameters $\theta_n$. This neural network has its own parameters, denoted by $\lambda_n$. By training the network, we can effectively control the parameters $\theta_n$ of each transformation $f_{\theta_n}^{(n)}$ based on the measurement information $\bm{x}$ and $\bm{d}$. The collection of all neural networks $\{\mathcal{N}^{(n)}_{\lambda_n}(\bm{v})\}_{n=1}^{N}$ can be compactly represented as $\mathcal{N}_{\lambda}(\bm{v})=\theta$, where $\lambda = \{\lambda_1, \lambda_2, \ldots, \lambda_N\}$ denotes the combined parameters of these neural networks.

Unlike functional normalizing flow introduced in Section \ref{sec:functional normalizing flow}, the conditional network $\mathcal{N}_{\lambda}(\bm{v})$ propose a novel approach where the model parameters for each measurement information  $(\bm{x}, \bm{d})$ are generated directly, rather than being learned through training. In other words, once the parameters $\lambda$ of the neural network $\mathcal{N}_{\lambda}(\bm{v})$ are trained, we can directly generate the parameters of the corresponding flow model for any given measurement information $(\bm{x}, \bm{d})$. Next, we detail the training process for the parameters $\lambda$ of the network $\mathcal{N}_{\lambda}(\bm{v})$.

The network $\mathcal{N}_{\lambda}(\bm{v})$ is trained on the dataset $D_{train}=\{(\bm{x}_i, \bm{d}_i)\}_{i=1}^{N_{train}}$, where $({\bm{x}}_i, {\bm{d}}_i)$ represents a pair of measurement data $\bm{d}_i$ and its corresponding measurement points $\bm{x}_i$ for $i = 1 , 2, \ldots, N_{train}$. For a specific measurement information pair $(\bm{x}^*,\bm{d}^*)$, the parameters of the corresponding functional normalizing flow are given by $\theta^*=\mathcal{N}_{\lambda}(\bm{v}^*)$, where $\bm{v}^*=\mathcal{W}(\bm{x}^*,\bm{d}^*)$. As shown in (\ref{eq:post2}), we denote the posterior corresponding to $(\bm{x}^*,\bm{d}^*)$ as $\mu_{(\bm{x}^*,\bm{d}^*)}$. We hypothesize that the flow model $f_{\theta^*}$ can transform the prior $\mu_0$ into a measure $\nu_{\lambda}(\bm{v}^*)$ that closely approximates the posterior $\mu_{(\bm{x}^*,\bm{d}^*)}$. Consequently, we expect the Kullback-Leibler divergence $D_{\text{KL}}(\nu_{\lambda}(\bm{v}^*)||\mu_{(\bm{x}^*,\bm{d}^*)})$ to be small.
Let $\{ \bm{v}_i=\mathcal{W}(\bm{x}_i,\bm{d}_i) \}_{i=1}^{N_{train}}$ denote the vectors corresponding to the dataset $D_{train}$, the overall error function can then be written as
\begin{align*}
q(\lambda)=&\mathbb{E}_{(\bm{x},\bm{d})}\{D_{\text{KL}}(\nu_{\lambda}(\bm{v})||\mu_{(\bm{x},\bm{d})})\}
\approx \frac{1}{N_{train}}\sum\limits_{i=1}^{N_{train}}D_{\text{KL}}(\nu_{\lambda}(\bm{v}_i)||\mu_{(\bm{x}_i,\bm{d}_i)}).
\end{align*}
For a detailed description of the algorithm, we refer the reader to Algorithm \ref{alg B}.

\begin{algorithm}[ht!]
\caption{Conditional Functional Normalizing Flow}
\label{alg B}
\begin{algorithmic}[1]
\STATE{Initialize the neural network parameter $\lambda$ to $\lambda_0$, the learning rate schedule to $\alpha_k$, the training steps to $K$, the sampling size to $M$, and the training dataset to $D_{train}=\{(\bm{x}_1, \bm{d}_1), (\bm{x}_2, \bm{d}_2), \ldots , (\bm{x}_{N_{train}},\bm{d}_{N_{train}})\}$. Then, initialize the iteration counter $k$ to $0$;}
\REPEAT
\STATE{Randomly select $M$ training data $\{(\bm{x}_{r_1},\bm{d}_{r_1}), \ldots , (\bm{x}_{r_M},\bm{d}_{r_M})\}$ from dataset $D_{train}$, where $r_i \in \{1, 2, \ldots, N_{train}\}$;}
\STATE{Calculate their corresponding vectors $\{\bm{v}_{r_1},\bm{v}_{r_2},\ldots,\bm{v}_{r_M}\}$ with $\bm{v}_{r_i}=\mathcal{W}(\bm{x}_{r_i},\bm{d}_{r_i})$;}
\STATE{Calculate the approximate measures $\{ \nu_{\lambda_k}(\bm{v}_{r_1}),\nu_{\lambda_k}(\bm{v}_{r_2}),\ldots,\nu_{\lambda_k}(\bm{v}_{r_M}) \}$ with the neural network $\mathcal{N}_{\lambda_k}(\bm{v})$;}
\STATE{For each $i=r_1, r_2,\ldots,r_M$, generate $N_u$ functions from the measure $\nu_{\lambda_k}(\bm{v}_{i})$, and denote them as $\{ u_{i1}^{(\lambda_k)}, u_{i2}^{(\lambda_k)}, \ldots, u_{iN_u}^{(\lambda_k)} \}$;}
\STATE{Updating the parameters $\lambda_{k+1}=\lambda_k-\alpha_k\nabla_{\lambda_k}q(\lambda_k)$ with \begin{align*}
\nabla_{\lambda_k}q(\lambda_k)&=\nabla_{\lambda_k}E_{(\bm{x},\bm{d})}\{D_{\text{KL}}(\nu_{\lambda_k}(\bm{v})||\mu_{(\bm{x},\bm{d})})\}\\ 
&\approx \frac{1}{MN_u}\nabla_{\lambda_k}\sum\limits_{i=1}^{M}\left(\sum\limits_{j=1}^{N_u}\text{ln}\left(\frac{d\nu_{\lambda_k}(\bm{v}_{r_i})}{d\mu_0}(u_{ij}^{(\lambda_k)})\right)- \text{ln}\left(\frac{\mu_{(\bm{x}_{r_i},\bm{d}_{r_i})}}{d\mu_0}(u_{ij}^{(\lambda_k)})\right)\right).
\end{align*}
Variants of stochastic gradient-based
optimization method can be employed in this step;}
\STATE{$k=k+1$;}
\UNTIL{$k = K$;}
\STATE{Return the final result $\lambda_{K}$;}
\end{algorithmic}
\end{algorithm}
{\color{black}
The neural network $\mathcal{N}_{\lambda}(\bm{v})$ serves as a computationally efficient estimator for the initial posterior approximation. For a specific measurement pair $(\bm{x}^*,\bm{d}^*)$, the input vector $\bm{v}^*=\mathcal{W}(\bm{x}^*,\bm{d}^*)$ produces the initial flow parameters $\theta^*=\mathcal{N}_{\lambda}(\bm{v}^*)$. The corresponding normalizing flow $f_{\theta^*}$ approximates the true posterior $\mu_{(\bm{x}^*,\bm{d}^*)}$ by pushing forward the prior  to the measure $\nu_{\lambda}(\bm{v}^*)$. However, due to the inherent capacity limits of $\mathcal{N}_{\lambda}(\bm{v})$, this initial estimate may deviate from the target posterior, potentially leading to inaccurate uncertainty quantification. To refine this approximation, we propose a sample-specific fine-tuning strategy: the inverse problem associated with $(\bm{x}^*,\bm{d}^*)$ is further optimized via Algorithm \ref{alg A} with the initialization $\theta_0=\theta^*$. As demonstrated in Section \ref{sec5}, this informed initialization significantly accelerates convergence compared to random starts.
}

\section{Numerical Examples of Functional Normalizing Flow}\label{sec4}
To demonstrate the effectiveness of our novel framework for solving inverse problems, we present three illustrative examples, i.e., the simple elliptic equation, the steady-state Darcy flow equation and the electrical impedance tomography (EIT) problem. These examples enable a straightforward visualization of the approximate posterior distribution. {\color{black} The code used for simulations is available on GitHub at \url{https://github.com/jjx323/FunctionalNormalizingFlow}. All programs ran on a system with a 13th Gen Intel(R) Core(TM) i7-13700, NVIDIA GeForce RTX 4090, and Ubuntu 20.04.5 LTS.} 

In the first example, given the linear nature of the inverse problem, we utilize functional Householder flow and functional projected transformation flow as our flow models. For the second and third examples, which involves a nonlinear inverse problem, we employ functional planar flow and functional Sylvester flow. These flow models are specifically designed to capture intricate nonlinear relationships, enabling the transformation of a Gaussian distribution into a non-Gaussian distribution that more closely approximates the posterior distribution associated with the nonlinear inverse problem. 

\subsection{Simple Elliptic Equation}\label{subsec4.1}
Consider the inverse source problem associated with the elliptic equation:
\begin{align}\label{prob1}
\begin{split}
 -\beta \Delta w + w &= u, \quad \text{in}\ \Omega, \\ 
 \frac{\partial w}{\partial \bm{n}} &= 0, \quad \text{on}\ \partial \Omega,
\end{split}
\end{align}
where $\Omega = (0, 1) \subset \mathbb{R}$, $\beta > 0$ is a positive constant ($\beta = 0.01$ in our experiment), and $\bm{n}$ denotes the outward unit normal vector. 
The forward operator is defined as
\begin{equation}
\mathcal{S}\mathcal{G}u =(w(x^1), w(x^2), \ldots, w(x^{N_{\bm{d}}}) )^T,
\end{equation}
where $\mathcal{G}$ is the PDE solution operator from $\mathcal{H}_u:=L^2(\Omega)$ to $\mathcal{H}_{w}:=H^2(\Omega)$, $\mathcal{S}$ is the measurement operator from $\mathcal{H}_{w}$ to $\mathbb{R}^{N_{\bm{d}}}$, $u \in \mathcal{H}_u$ denotes the model parameter, $w\in \mathcal{H}_w$ denotes the solution of the elliptic equation $(\ref{prob1})$, and $x^i \in \Omega$ for $i = 1, 2, \ldots ,  N_{\bm{d}}$.
With these notations, the problem can be written abstractly as
\begin{align}\label{equ:linear inverse}
 \bm{d} = \mathcal{S}\mathcal{G}u + \bm{\epsilon}.
\end{align}
For clarity, we list the specific choices for some parameters introduced in this subsection as follows:
\begin{itemize}
	\item \textcolor{black}{Assume that 5$\%$ random Gaussian noise $\bm{\epsilon} \sim \mathcal{N}(0, \bm{\Gamma}_{\text{noise}})$ is added, where $\bm{\Gamma}_{\text{noise}} = \tau^{-1}\textbf{I}$, and $\tau^{-1} = (0.05 \lVert \mathcal{S}\mathcal{G}u\rVert_{\infty})^2$.}
	\item Let the domain $\Omega$ be an interval $(0, 1)$ with $\partial \Omega = \lbrace 0, 1 \rbrace$, and the measurement data are assumed to be $\lbrace w(x^i) | i = 0, 1, \ldots, 10 \rbrace$, where $x^{i}=i/10$. 
	\item The covariance operator \(\mathcal{C}_0\) associated with the prior measure \(\mu_0\) is defined as \(\mathcal{C}_0 = (\text{I} - \alpha \Delta)^{-2}\), where \(\alpha = 0.1\) is a predetermined constant. The Laplace operator is defined on \(\Omega\) with homogeneous Neumann boundary condition. Additionally, the mean of the prior measure \(\mu_0\) is set to zero.
	\item In order to avoid inverse crime \cite{kaipio2006statistical}, the data is generated on a fine mesh with the number of grid points equal to $10^4$. And we use different sizes of mesh $n = \lbrace 50, 75, 100, 200, 300 \rbrace$ in the inverse stage.
\end{itemize}

Since the inverse problem under consideration is linear, this subsection demonstrates the efficacy of two linear flows---the functional Householder flow and the functional projected transformation flow---as well as the discretization invariance of the proposed model. 
We assume that the data produced from the underlying ground truth
\begin{align}\label{equ:smoothtruth}
		u^{\dagger} = \exp(-50(x-0.3)^2)-\exp(-50(x-0.7)^2).
\end{align}
To demonstrate the effectiveness of our proposed algorithm, we conduct a comparative experiment with the preconditioned Crank-Nicolson (pCN) algorithm. It has been theoretically established that the pCN algorithm can generate samples from the true posterior distribution \cite{dashti2013bayesian}. Consequently, we treat the samples obtained from pCN as benchmark results for comparison. However, as a Markov chain Monte Carlo (MCMC) method, pCN samples often exhibit high autocorrelation. To acquire a sufficient number of independent samples, numerous iterations are typically required, and the computational cost of pCN scales directly with the desired number of independent samples. In contrast, the computational burden of NF-iVI primarily lies in the training process of the flow model. Once trained, the generated samples are mutually independent, which enables efficient sampling. {\color{black} When performing the pCN algorithm, we produce $3 \times 10^6$
samples for each chain, and discard the first $10^5$
samples as burn-in when calculating quantities of interest (similar with \cite{pinski2015algorithms, cotter2013}). It is notable that running the pCN algorithm required approximately 9 minutes, whereas training our functional normalizing flow model took only 25 seconds.}

\begin{figure}[ht!]
	\centering
	
	\subfloat[Householder flow]{
		\includegraphics[ keepaspectratio=true, width=0.32\textwidth, clip=true]{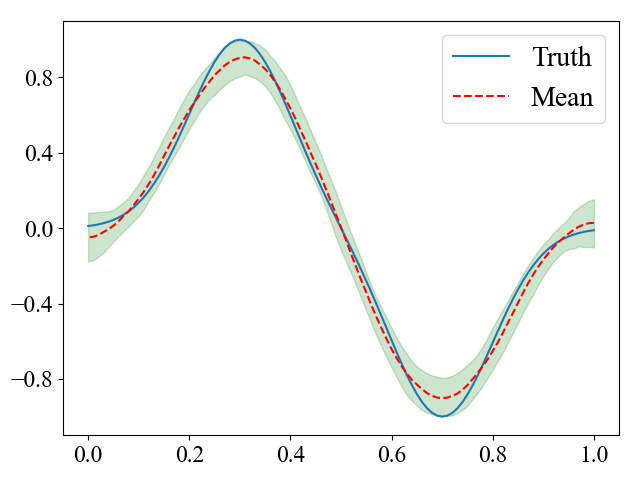}} 
    \subfloat[projected transformation flow]{
		\includegraphics[ keepaspectratio=true, width=0.32\textwidth, clip=true]{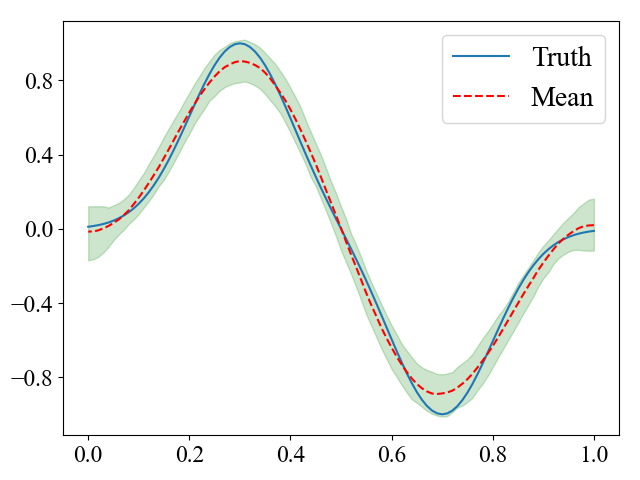}} 
	\subfloat[pCN algorithm]{
		\includegraphics[ keepaspectratio=true, width=0.32\textwidth, clip=true]{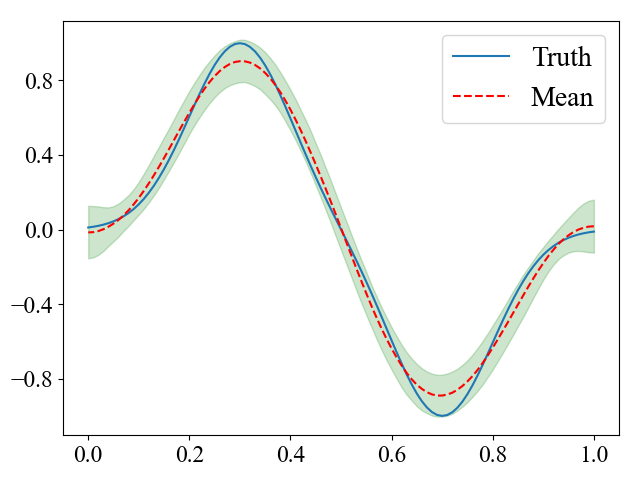}} 
	
	\caption{\emph{\small A comparison of the approximate posterior obtained by functional Householder flow with 24 layers Householder flow, 5 layers functional projected transformation flow and the posterior generated by the pCN algorithm. The green shade area represents the $95\%$ credibility region of estimated posterior. (a)(b): Result of the approximate posterior generated by functional Householder flow, functional projected transformation flow compared with the background truth of $u$. (c): Result of the posterior generated by pCN algorithm compared with the background truth of $u$.}}
	
       \label{fig:linearflow}
	
\end{figure}

Next, we will compare the posterior distribution obtained by the two algorithms in detail to verify the effectiveness of functional normalizing flow. Firstly, we provide some discussions of the mean function of approximate posterior. In subfigure (a)(b) of Figure \ref{fig:linearflow}, we compare the mean of the approximate posterior obtained by the NF-iVI method with the ground truth. The blue solid line represents the mean of NF-iVI estimate, while the red dashed line represents the ground truth. The shaded green region denotes the $95\%$ credibility interval of the approximate posterior. We observe that the $95\%$ credibility interval encompasses the ground truth, reflecting the inherent uncertainty in the parameter estimation. Subfigure (c) of Figure \ref{fig:linearflow} presents a similar comparison between the mean of the estimated posterior measure obtained by the pCN algorithm and the ground truth. As illustrated in the comparison between (a)(b) and (c), the mean of the approximate posterior generated by two linear flows closely resembles the mean of posterior obtained from the pCN algorithm.

To further support the conclusion, we provide numerical evidence. By comparing NF-iVI with the pCN sampling algorithm, we calculate the relative error between their estimates:
\begin{align}
    \textbf{Householder relative error:}= \frac{\lVert u_{HN}^{*}-u_{p}^{*}\rVert_{L^2}^2}{\lVert u_{p}^{*}\rVert_{L^2}^2}=0.00129, \label{equ:relative error smooth1}\\
    \textbf{projected transformation relative error:}= \frac{\lVert u_{PN}^{*}-u_{p}^{*}\rVert_{L^2}^2}{\lVert u_{p}^{*}\rVert_{L^2}^2}=0.00271, \label{equ:relative error smooth2}
\end{align}
where $u_{HN}^{*}$ denotes the estimated mean of functional Householder flow, $u_{PN}^{*}$ denotes the estimated mean of functional projected transformation flow, and $u_p^{*}$  represents the estimated mean of pCN algorithm. Obviously, the small relative errors suggest the effectiveness of our proposed algorithm.

As a result, combining visual (Figure \ref{fig:linearflow})
and quantitative (relative errors given in Equations (\ref{equ:relative error smooth1}), (\ref{equ:relative error smooth2})) evidences, the NF-iVI method provides a good approximation of the mean of posterior.

Secondly, we offer a discussion on the approximation of the posterior covariance function. The emphasis is placed on demonstrating that the NF-iVI method provides robust covariance estimates. To facilitate these comparisons, we need to introduce the definition of the variance and the covariance function.

Consider a random field $u$ on a domain $\Omega$ with mean $\bar{u}$ and the covariance function $c(x, y)$ describing the covariance between $u(x)$ and $u(y)$:
\begin{align*}
	c(x, y) = \mathbb{E}((u(x)-\bar{u}(x))(u(y)-\bar{u}(y))), \quad \text{for} \ x, y \in \Omega.
\end{align*}
The corresponding covariance operator $\mathcal{C}$ is 
\begin{align*}
	(\mathcal{C}\phi^{\prime})(x) = \int_{\Omega}c(x, y)\phi^{\prime}(y)dy,
\end{align*}
where the function $\phi^{\prime}$ is a sufficiently regular function defined on $\Omega$.

\begin{figure}[ht!]
	\centering
	
	\subfloat[variance (Householder)]{
		\includegraphics[ keepaspectratio=true, width=0.32\textwidth, ]{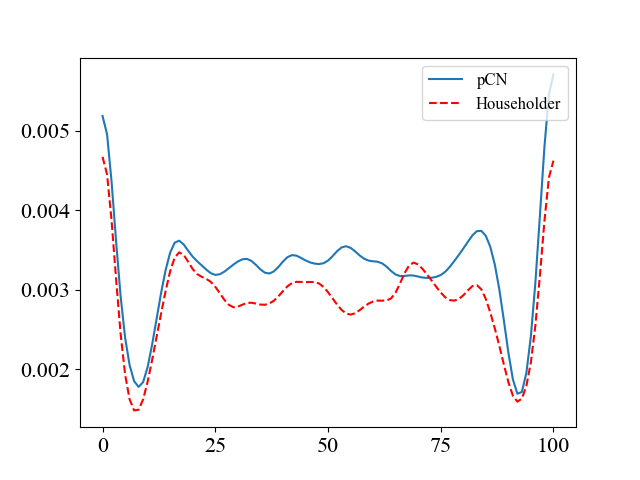}}
	\subfloat[covariance (Householder)]{
		\includegraphics[ keepaspectratio=true, width=0.32\textwidth, ]{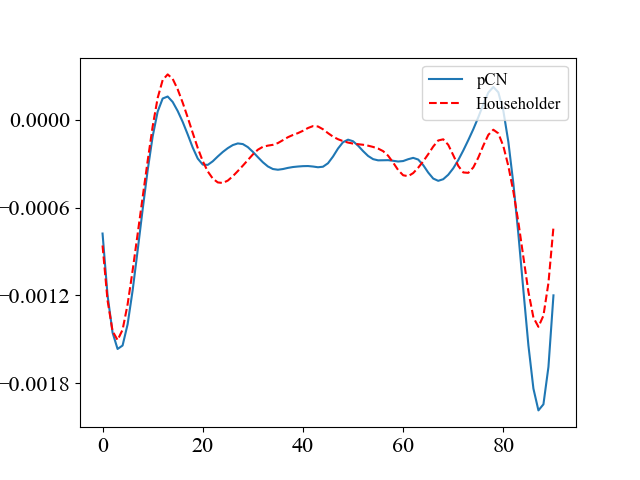}} 
	\subfloat[covariance (Householder)]{
		\includegraphics[ keepaspectratio=true, width=0.32\textwidth, ]{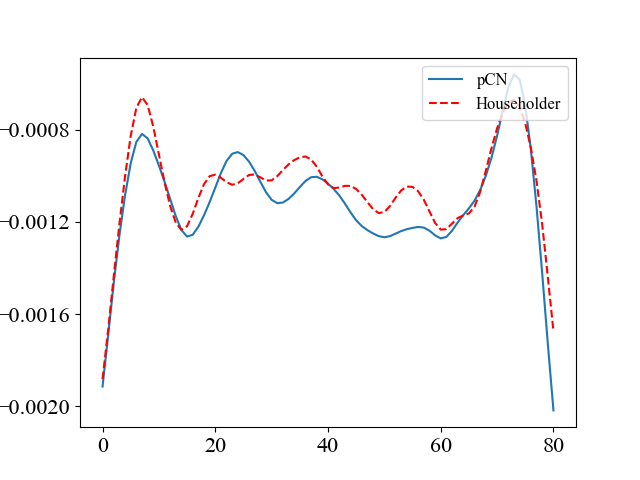}}

     \subfloat[variance (projected)]{
		\includegraphics[ keepaspectratio=true, width=0.32\textwidth, ]{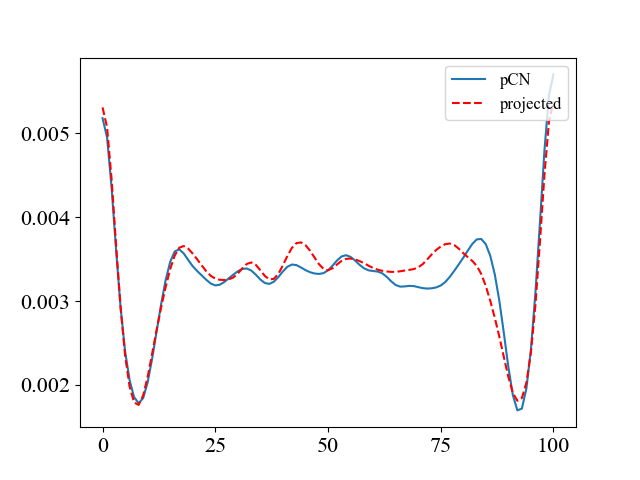}}
	\subfloat[covariance (projected)]{
		\includegraphics[ keepaspectratio=true, width=0.32\textwidth, ]{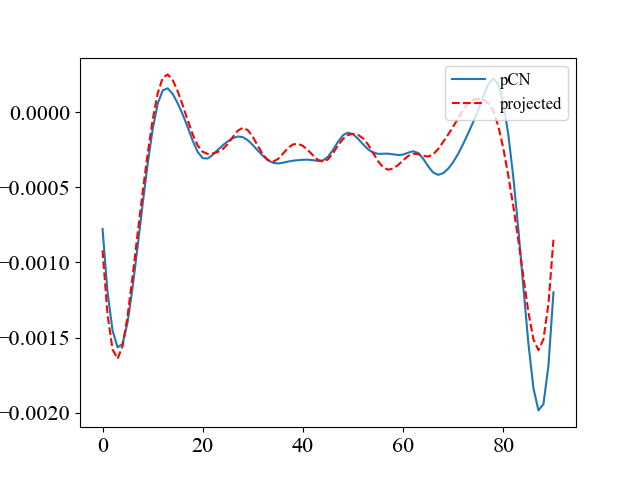}} 
	\subfloat[covariance (projected)]{
		\includegraphics[ keepaspectratio=true, width=0.32\textwidth, ]{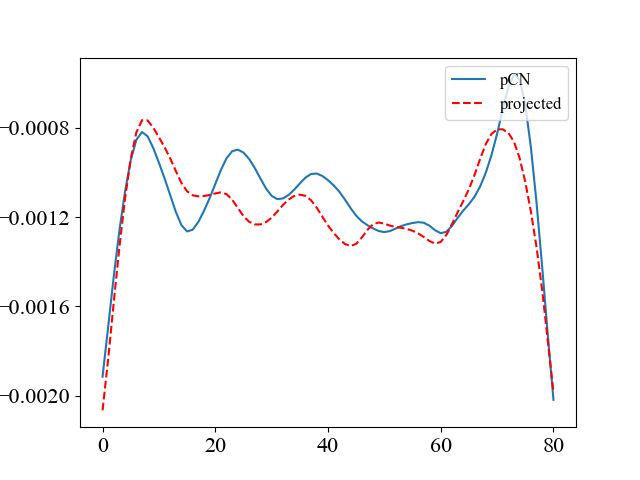}}

	\caption{\emph{\small {\color{black} The estimated variance and covariance functions obtained by the pCN algorithm (blue solid line), and NF-iVI algorithm (red dashed line).
				(a)(d): The covariance function $c(x,y)$ on all the mesh point pairs $\lbrace (x_i, x_{i})\rbrace^{N_g}_{i=1}$;
				(b)(e): The covariance function $c(x,y)$ on the mesh points $\lbrace (x_i, x_{i+10}) \rbrace^{N_g-10}_{i=1}$;
				(c)(f): The covariance function $c(x,y)$ on the mesh points $\lbrace (x_i, x_{i+20}) \rbrace^{N_g-20}_{i=1}$}}.}
\label{pic:cov linear}
	
\end{figure}

\begin{figure}[ht!]
	\centering
	
	\subfloat[covariance (Householder)]{
		\includegraphics[ keepaspectratio=true, width=0.32\textwidth, clip=true]{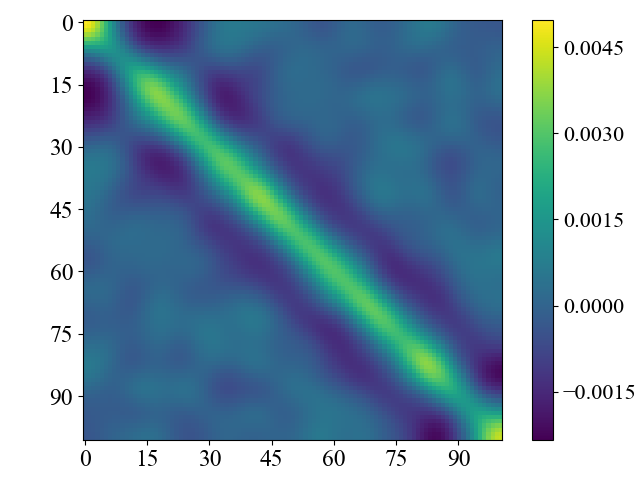}} 
    \subfloat[covariance (projected)]{
		\includegraphics[ keepaspectratio=true, width=0.32\textwidth, clip=true]{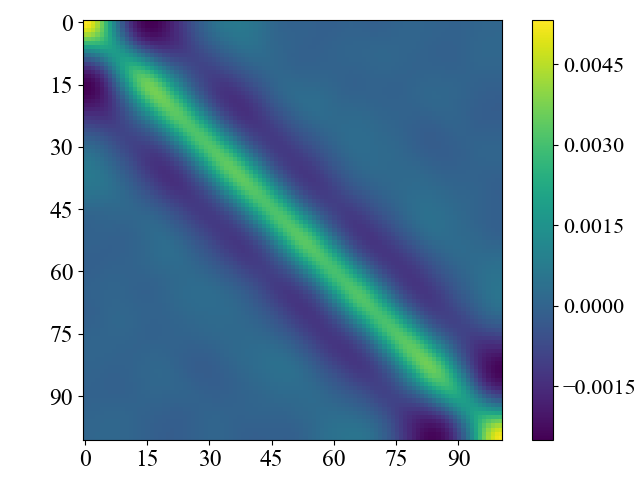}} 
	\subfloat[covariance (pCN)]{
		\includegraphics[ keepaspectratio=true, width=0.32\textwidth, clip=true]{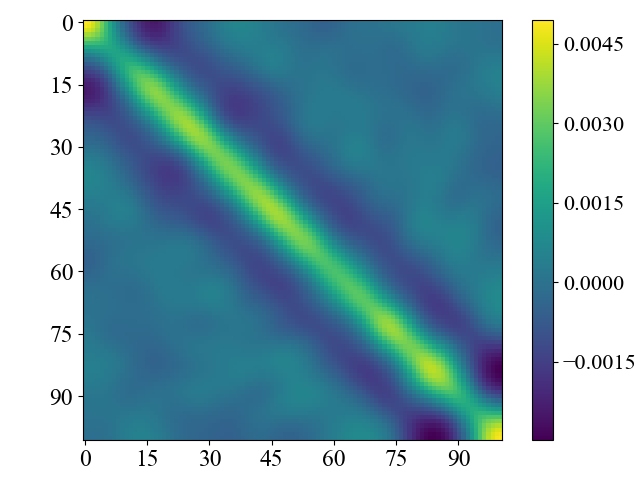}} \\

	\caption{\emph{\small {\color{black}The comparison of covariance. (a): The covariance given by Householder NF-iVI method; (b): The covariance given by Projected NF-iVI method; (c): The covariance given by the pCN method;}}}
\label{pic:var linear}
	
\end{figure}

The variance function, denoted by $\text{var}_u(x):=c(x,x)$, can be estimated as 
$$
\text{var}_u(x)\approx\frac{1}{m-1}\sum_{i=1}^m(u_i(x)-\bar{u}(x))^2,
$$
where $x\in\Omega$ is a point residing in the domain $\Omega,\bar{u}$ is the mean function, and $m$ is the number of samples. The covariance function can be estimated as  $$\text{cov}_u(x_1,x_2)\approx\frac1{m-1}\sum_{i=1}^m(u_i(x_1)-\bar{u}(x_1))(u_i(x_2)-\bar{u}(x_2)),$$ where $x_1,x_2\in\Omega$ and $u,\bar{u}$ are defined as in var$_u(x).$ For simplicity, we compute these quantities on the mesh points and exhibit the results in Figure \ref{pic:cov linear}. In all of the subfigures in Figure \ref{pic:cov linear}, the estimates obtained by the pCN and NF-iVI are drawn in blue solid line and red dashed line, respectively. In Figure \ref{pic:cov linear} (a)(d), we show the variance function calculated on all of the mesh points, i.e., $\{\text{var}_{u}(x_{i})\}_{i=1}^{N_{g}}$ ($N_{g}$ is the number of mesh points, which equals to $101$ in this problem). In Figure \ref{pic:cov linear} (b)(e) and (c)(f), we show the covariance function calculated on the pairs of points $\left\{(x_i,x_{i+10})\right\}_{i=1}^{N_g-10}$ and $\{(x_i,x_{i+20})\}_{i=1}^{N_g-20}$, respectively. Subfigures (a), (b), and (c) in Figure \ref{pic:var linear} depict the matrix representations of the covariance operators \(\mathcal{C}_{\text{HN}}\), \(\mathcal{C}_{\text{PN}}\) and \(\mathcal{C}_{\text{p}}\). Here $\mathcal{C}_{\text{HN}}$, \(\mathcal{C}_{\text{PN}}\) represents the covariance operator corresponding to the Householder flow and projected transformation flow, while $\mathcal{C}_{\text{p}}$ represents the covariance operator corresponding to the pCN algorithm. The results confirm that the NF-iVI method provides an accurate estimate of the covariance function. 


Then we make numerical comparisons of the variance and covariance functions by showing the relative errors:
\begin{align*}
    \text{relative error of variance} &= \frac{\sum^{N_{g}}_{i=1}\big (c_N(x_i, x_i) - c_p(x_i, x_i)\big )^2}{\sum^{N_{g}}_{i=1}\big (c_p(x_i, x_i)\big)^2}, \\
    \text{relative error of covariance} &= \frac{\sum^{N_{g}-k}_{i=1}\big(c_N(x_i, x_{i+k}) - c_p(x_i, x_{i+k})\big)^2}{\sum^{N_{g}-k}_{i=1}\big(c_p(x_i, x_{i+k})\big)^2},\\
    \text{total relative error $\bm{c}$} &= \frac{\sum^{N_{g}}_{i=1}\sum^{N_{g}}_{j=1}\big(c_N(x_i, x_{j}) - c_p(x_i, x_{j})\big)^2}{\sum^{N_{g}}_{i=1}\sum^{N_{g}}_{i=1}\big(c_p(x_i, x_{j})\big)^2},
\end{align*}
where $k$ is an integer, $c_p(x,y)$ is the covariance function generated from pCN algorithm, and $c_N(x,y)$ is the covariance function generated from NF-iVI method. The relative errors are shown in Table \ref{table:relative linear}. Here, the notation $\bm{c}$ means the total relative error, and the notation $\lbrace c(x_i, x_{i+k}) \rbrace^{N_g-k}_{i=1}$ means the covariance function values on the pair of mesh points $\lbrace (x_i, x_{i+k}) \rbrace^{N_g-k}_{i=1}$ ($k=0, 10, 20$). The numbers below this notation are the relative errors between the vectors obtained by pCN method and NF-iVI, respectively.

As depicted in Table \ref{table:relative linear}, the relative errors are small, indicating that the posterior covariance functions obtained by the NF-iVI and pCN sampling algorithms are quantitatively similar.

\begin{table}[ht!]
	\renewcommand{\arraystretch}{1.5}
	\centering
		\caption{\emph{\small The relative errors between variance function, and covariance functions.}} 
	\begin{tabular}{c|cccc}
		\hline $\text{Relative Error} $& $\bm{c}$  & $\lbrace c(x_i, x_i)\rbrace^{N_g}_{i=1}$& $\lbrace c(x_i, x_{i+10}) \rbrace^{N_g-10}_{i=1}$& $\lbrace c(x_i, x_{i+20}) \rbrace^{N_g-20}_{i=1}$\\
		\hline $\text{Householder Flow}$& $0.03513$& $0.01911$& $0.09921$& $0.00921$\\
        \hline $\text{Projected Transformation}$ & $0.03656$& $0.00686$& $0.11220$& $0.03715$\\
		\hline
	\end{tabular}
 \label{table:relative linear}
\end{table}

Additionally, we will leverage this problem to demonstrate the discretization invariance of our flow model. {\color{black} To maintain conciseness, we present the discretization-invariance numerical experiments only for the Householder flow. The analogous results for the other three proposed flow models, which exhibit similar properties, are omitted but can be fully reproduced and verified using the publicly available code.}

To assess the impact of discretization, we discretize the problem using 50, 75, 100, 200, and 300 grid points in a regular mesh, and compute the \(L^2\)-error between the mean functions of the approximate posterior and the background truth across various discrete levels. In Table \ref{table:L2error1}, the numbers above are different discrete level we used for NF-iVI, and the numbers below are the relative errors between the means function obtained by NF-iVI and the true function, respectively. As detailed in Table \ref{table:L2error1}, the \(L^2\)-errors between the mean functions derived from our algorithm and the truth are consistent across different discretizations. Furthermore, Table \ref{table:covariancesame1} presents a comparison of the covariance functions of the approximate posterior for several discrete levels. In the table, the numbers above are different discrete level we used for NF-iVI, and the numbers below are the exact values of the covariance function evaluated at the points $(x_1, x_1), (x_1,x_2), (x_1,x_3), (x_2,x_2)$, where $x_1=0.3, x_2=0.5, x_3=0.7$. Collectively, these findings substantiate the assertion that the approximate posteriors generated by our proposed algorithm are invariant to the different discrete levels.

\begin{table}[ht!]
	\renewcommand{\arraystretch}{1.5}
	\centering
		\caption{\emph{\small The $L^2$-errors between the mean of approximate posterior and true function.}} 
	\begin{tabular}{c|ccccc}
		\hline $\text{Discrete Level} $& 50 & 75 & 100& 200 & 300\\
		\hline $\text{$L^2$-Error}$& $0.00286$& $0.00294$& $0.00305$& $0.00311$& $0.00302$\\
		\hline
	\end{tabular}
 \label{table:L2error1}
\end{table}

\begin{table}[ht!]
	\renewcommand{\arraystretch}{1.5}
	\centering
		\caption{\emph{\small The values of covariance function, where $x_1=0.3$, $x_2=0.5$, $x_3=0.7$}.} 
	\begin{tabular}{c|ccccc}
		\hline $\text{Discrete Level} $& 50  & 75 & 100 & 200 & 300\\
		\hline $\text{$c(x_1,x_1)$}$& $0.00353$& $0.00340$& $0.00323$& $0.00326$&$0.00315$\\
  		\hline $\text{$c(x_1,x_2)$}$& $-0.00120$& $-0.00123$& $-0.00122$& $-0.00130$&$-0.00127$\\
    		\hline $\text{$c(x_1,x_3)$}$& $0.00028$& $0.00036$& $0.00016$& $0.00030$&$0.00023$\\
      		\hline $\text{$c(x_2,x_2)$}$& $0.00299$& $0.00272$& $0.00290$& $0.00329$&$0.00287$\\
		\hline
	\end{tabular}
 \label{table:covariancesame1}
\end{table}

By integrating the numerical findings presented in Tables \ref{table:L2error1} and \ref{table:covariancesame1} with the theoretical results established in Theorem \ref{thm:discrete invariance}, we have successfully demonstrated the discretization invariance of the flow model.

As discussed in Subsection \ref{subsec:linear transformation}, functional projected transformation flow are expected to exhibit superior approximation capabilities compared with functional Householder flow. To demonstrate this, we utilize 5 layers for the functional projected transformation flow and explore different configurations with 5, 9, 14, 18, and 21 layers for the functional Householder flow. Following Algorithm \ref{alg A}, all the functional normalizing flow models are trained for 5000 iterations ($K=5000$) using 30 samples ($N=30$) in each iteration. The initial learning rate is set to $\alpha_0=0.01$, with a step decay learning rate schedule of a multiplicative factor $\tau=0.8$ after every 500 iterations.

\begin{figure}[ht!]
	\centering
	
	\subfloat[projected transformation flow-5]{
		\includegraphics[ keepaspectratio=true, width=0.31\textwidth, clip=true]{PIC/linear/1D_pro.png}} 
	\subfloat[Householder flow-5]{
		\includegraphics[ keepaspectratio=true, width=0.31\textwidth, clip=true]{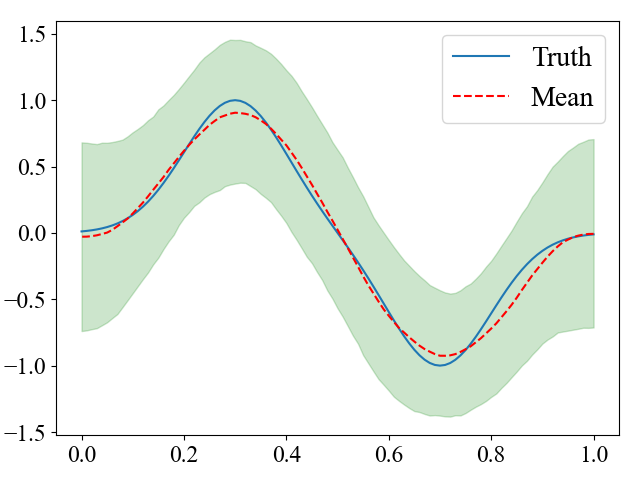}} 
        \subfloat[Householder flow-9]{
		\includegraphics[ keepaspectratio=true, width=0.31\textwidth, clip=true]{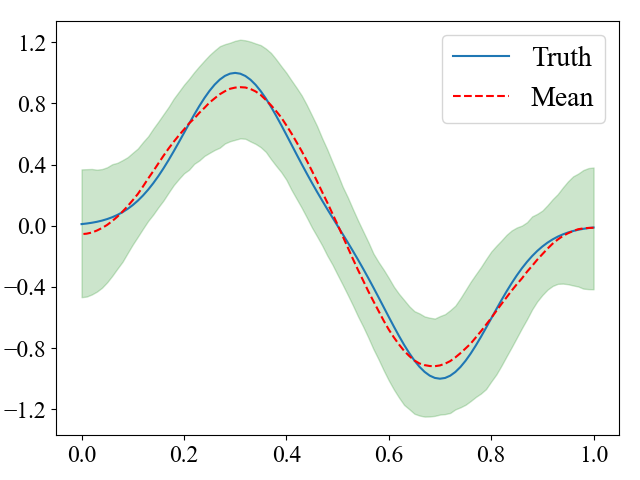}} 
	\\
	\subfloat[Householder flow-14]{
		\includegraphics[ keepaspectratio=true, width=0.31\textwidth, clip=true]{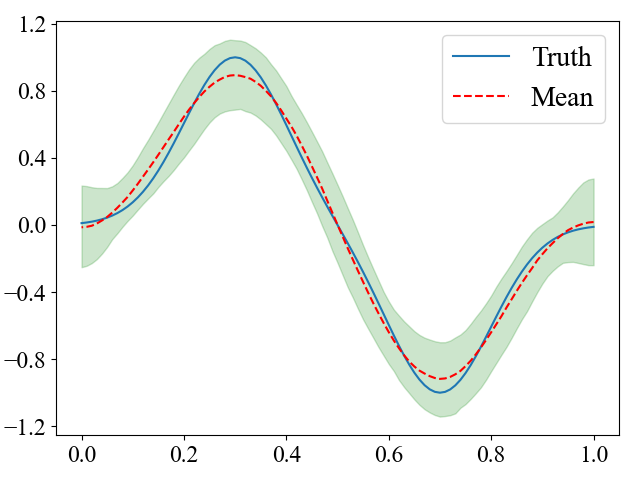}} 
	\subfloat[Householder flow-18]{
		\includegraphics[ keepaspectratio=true, width=0.31\textwidth, clip=true]{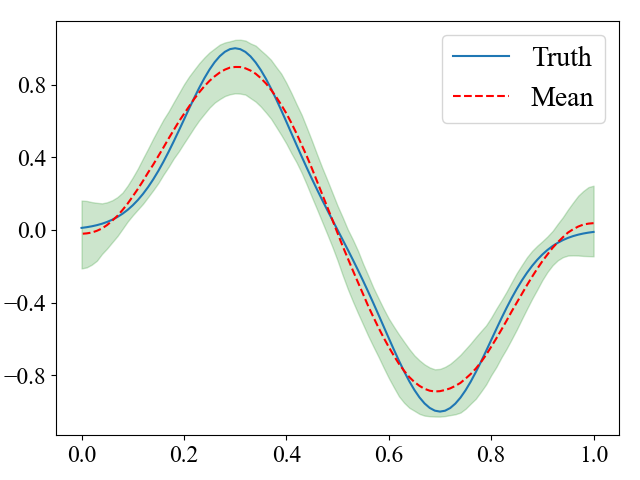}} 
     \subfloat[Householder flow-21]{
		\includegraphics[ keepaspectratio=true, width=0.31\textwidth, clip=true]{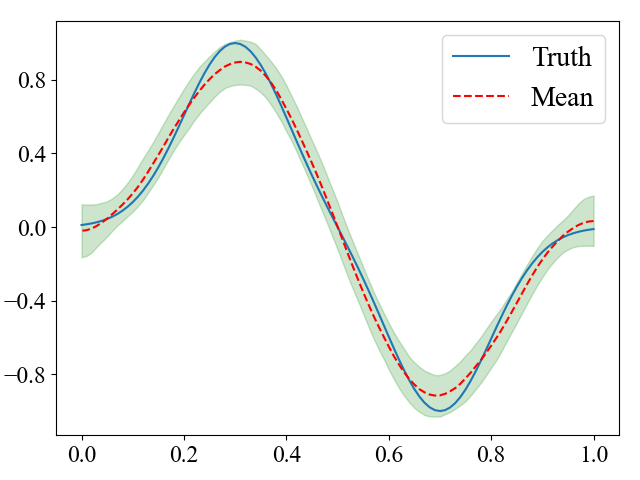}} 
	\caption{\emph{\small The comparison of the approximation capabilities of functional Householder flow and functional projected transformation flow. The green shade area represents the $95\%$ credibility region of approximate posterior. (a): Result of functional projected transformation flow with 5 layers. (b)-(f): Result of functional Householder flow with 5, 9, 14, 18, and 21 layers, respectively. }}
	
       \label{fig:projectflow}
	
\end{figure}

In subfigure (a) of Figure \ref{fig:projectflow}, we present the results obtained from the 5-layer functional projected transformation flow. We observe that these results are consistent with the posterior generated by the pCN algorithm, as shown in subfigure (c) of Figure \ref{fig:linearflow}, indicating the effectiveness of our proposed method. In subfigures (b)-(f) of Figure \ref{fig:projectflow}, we present the results obtained from the functional Householder flow with 5, 9, 14, 18, and 21 layers, respectively. We observe that when the number of layers is small, the functional Householder flow may not achieve satisfactory approximations.

The comparisons presented above highlight the superior approximation capabilities of the functional projected transformation flow, even when using a small number of layers. 
{\color{black}
Furthermore, to demonstrate the advantage of the Bayesian approach in solving inverse problems, we utilize a reduced number of measurement points to investigate the resulting posterior distribution under the constraint of sparse measurement data. Our algorithm employ a 5-layer projected transformation flow. The inversion is performed under three measurement scenarios, described below:

\begin{itemize}
\item \textbf{Full measurement}: The measurement data are assumed to be $\lbrace w(x^i) | x^{i}=i/10, i = 0, 1, 2, \ldots, 10 \rbrace$.  
\item \textbf{Left measurement}: The measurement data are assumed to be $\lbrace w(x^i) | x^{i}=i/10, i = 0, 1, \ldots, 5 \rbrace$.  
\item \textbf{Right measurement}: The measurement data are assumed to be $\lbrace w(x^i) | x^{i}=i/10, i = 5, 6, \ldots, 10 \rbrace$.  
\end{itemize}

Full measurement scenario involves measurement points spread across the domain $\Omega$. Left measurement scenario localizes measurement points on the left side of $\Omega$, and right measurement scenario localizes the measurement points on the right. All other settings for the inverse problem remain unchanged. The resulting posterior corresponding to these three distinct measurements are presented in Figure \ref{pic:sparse_measurement} below.

\begin{figure}[ht!]
	\centering
	\subfloat[full measurement]{
		\includegraphics[ keepaspectratio=true, width=0.32\textwidth, clip=true]{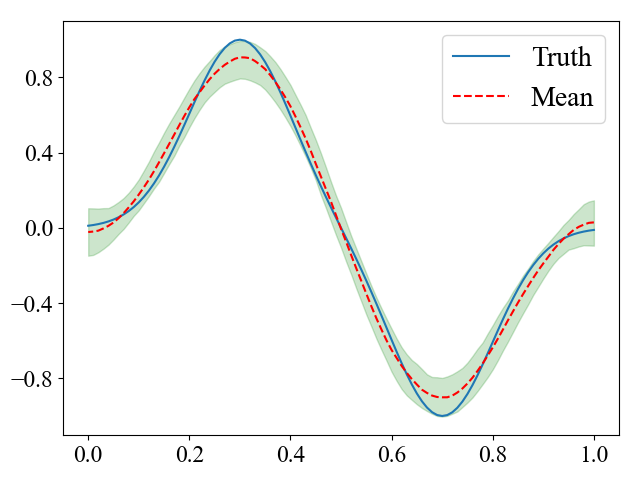}} 
    \subfloat[left measurement]{
		\includegraphics[ keepaspectratio=true, width=0.32\textwidth, clip=true]{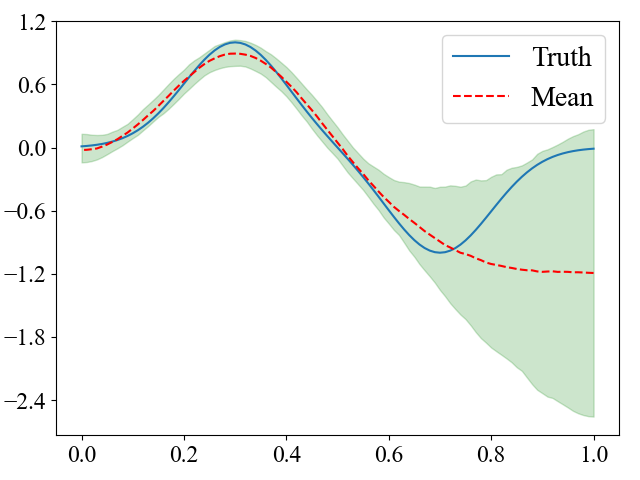}} 
	\subfloat[right measurement]{
		\includegraphics[ keepaspectratio=true, width=0.32\textwidth, clip=true]{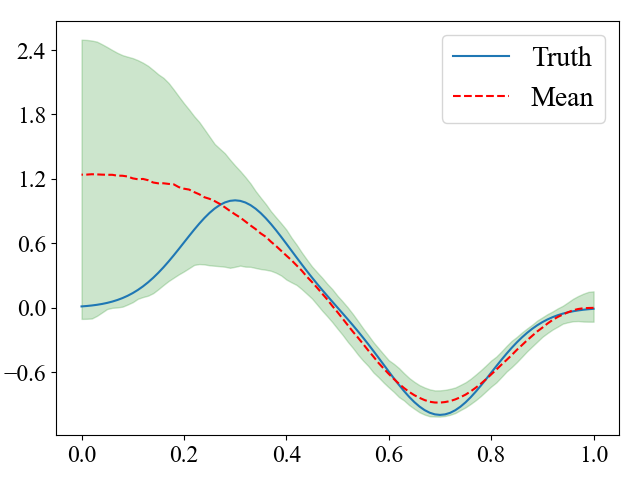}} \\
	
	\caption{\emph{\small A comparison of the approximate posterior with different measurement scenarios. The green shade area represents the $95\%$ credibility region of estimated posterior. (a): Posterior of full measurement. (b): Posterior of left measurement. (c): Posterior of right measurement.
}}
\label{pic:sparse_measurement}	
\end{figure}

Subfigures (a), (b), and (c) of Figure \ref{pic:sparse_measurement} present the results of the inverse problem with full measurement, left measurement, and right measurement, respectively. It is evident that for this linear inverse problem, the posterior distribution---when informed by full measurements---accurately reflects the true state of the model parameters. Conversely, when observations are sparse (specifically, when measurements are concentrated on the left side of the domain while the right side is neglected), the posterior distribution exhibits narrow credible intervals on the left, with the posterior mean closely matching the truth. In sharp contrast, the right side is characterized by significantly wider credible intervals, and the posterior mean deviates substantially from the ground truth. An analogous behavior is observed when measurements are concentrated on the right side of the domain.

These results do not reflect a deficiency in the model’s performance under sparse measurement settings; instead, they indicate that the measurement data are insufficient to accurately recover the true model parameters in both left-side and right-side measurement scenarios. When measurement data are inadequate, the credible regions provided by the Bayesian posterior are critical, as they delineate the regions where accurate inference is feasible and, conversely, where reliable inference cannot be achieved.}

\subsection{Steady-State Darcy Flow Equation}\label{subsec:darcy}
In this subsection, we focus on the inverse problem of estimating the permeability distribution in a porous medium from a discrete set of pressure measurements, as studied in \cite{calvetti2018iterative}. Consider the following steady-state Darcy flow equation:
\begin{align}\label{equ:darcy}
\begin{split}
 -\nabla \cdot (e^u\nabla w) &= f, \quad x \in \Omega, \\
 w &= 0, \quad x \in \partial \Omega,
\end{split}
\end{align}
where $f \in H^{-1}(\Omega)$ is the source function ($f=sin(\pi x_1)sin(\pi x_2)$ in our experiment, where $x=(x_1,x_2)$), and $u \in \mathcal{H}_u := L^{\infty}(\Omega)$ is called log-permeability for the computational area $\Omega = (0, 1)^2$. The forward operator has the following form:
\begin{align*}
 \mathcal{S}\mathcal{G}(u) = (w(x^1), w(x^2), \ldots, w(x^{N_d}))^T,
\end{align*}
where $\mathcal{G}$ is the PDE solution operator from $\mathcal{H}_u$ to $\mathcal{H}_{w}$, $\mathcal{S}$ is the measurement operator from $\mathcal{H}_{w}$ to $\mathbb{R}^{N_d}$, and $x^i \in \Omega$ for $i = 1, \ldots, N_{\bm{d}}$. With these notations, the problem can be written abstractly as:
\begin{align}\label{equ:nonlinear}
 \bm{d} = \mathcal{S}\mathcal{G}(u) + \bm{\epsilon},
\end{align}
where $\bm{\epsilon} \sim \mathcal{N}(0, \bm{\Gamma}_{\text{noise}})$ is the random Gaussian noise, and $\bm{d}$ is the measurement data.

In our experiments, the prior measure of $u$ is a Gaussian probability measure $\mu_0$ with mean zero and covariance $\mathcal{C}_0$. For clarity, we list the specific choices for some parameters introduced in this subsection as follows:
\begin{itemize}
	\item Assume that 5$\%$ random Gaussian noise $\bm{\epsilon} \sim \mathcal{N}(0, \bm{\Gamma}_{\text{noise}})$ is added, where $\bm{\Gamma}_{\text{noise}} = \tau^{-1}\textbf{I}$, and $\tau^{-1} = (0.05\lVert \mathcal{S}\mathcal{G}u\rVert_{\infty})^2$.
    \item We assume that the data is produced from the underlying log-permeability: 
    \begin{align*}
		u^{\dagger} = \ &\exp(-20(x_1-0.3)^2 - 20(x_2-0.3)^2) \\
		&+ \exp(-20(x_1-0.7)^2 - 20(x_2-0.7)^2).
	\end{align*}
    \item Let the domain $\Omega$ be a bounded area $(0, 1)^2$. The measurement points $\lbrace x^i \rbrace^{N_d}_{i=1}$ are taken at the coordinates $\{ i/21,j/21 \}_{i,j=1}^{20}$.
    \item The operator $\mathcal{C}_0$ is given by $\mathcal{C}_0 = (\text{I} - \alpha\Delta)^{-2}$, where $\alpha = 0.1$ is a fixed constant. Here, the Laplace operator is defined on $\Omega$ with zero Neumann boundary condition.
    \item To avoid the inverse crime, a fine mesh with the number of grid points equal to $500 \times 500$ is employed for generating the data. For the inversion, a mesh with a number of grid points equal to $20 \times 20$ is employed.
\end{itemize}
Given the nonlinear nature of the inverse problem, we employ two nonlinear flows, i.e., functional planar flow and functional Sylvester flow to approximate the posterior. Functional planar flow is constructed with 32 layers, while functional Sylvester flow is constructed with 5 layers. Following Algorithm \ref{alg A}, all the flow models are trained over 5000 iterations $(K=5000)$, with each iteration utilizing 30 samples $(N=30)$. We utilize a step decay learning rate schedule, starting with an initial learning rate of $\alpha_0=0.01$ and decreasing it by a multiplicative factor of $\tau=0.8$ after every 500 iterations. For comparison, we also present the posterior generated by the pCN algorithm with $3 \times 10^6$ samples. {\color{black} It is worth noting that running the pCN algorithm required approximately 472 minutes, whereas training the functional normalizing flow model took only about 28 minutes.}

\begin{figure}[ht!]
	\centering
	
	\subfloat[mean (pCN)]{
		\includegraphics[ keepaspectratio=true, width=0.31\textwidth, clip=true]{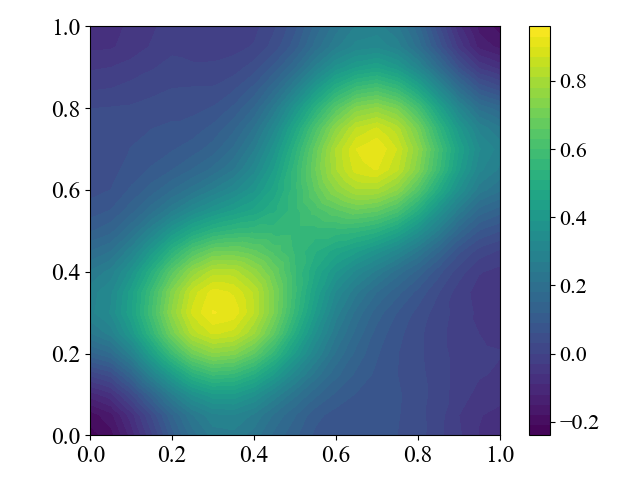}} 
	\subfloat[mean (planar)]{
		\includegraphics[ keepaspectratio=true, width=0.31\textwidth, clip=true]{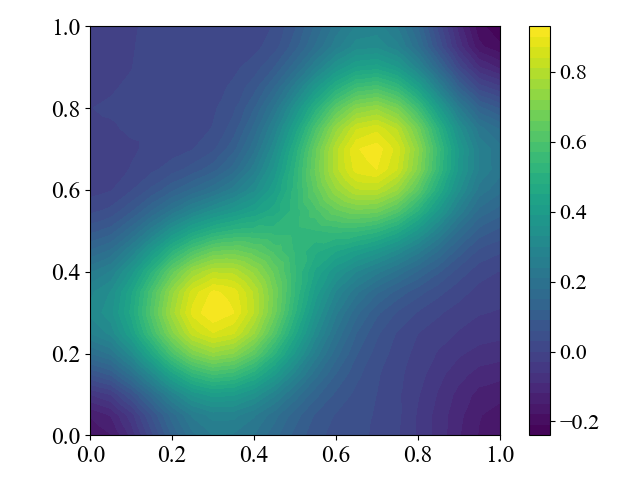}} 
        \subfloat[mean (Sylvester)]{
		\includegraphics[ keepaspectratio=true, width=0.31\textwidth, clip=true]{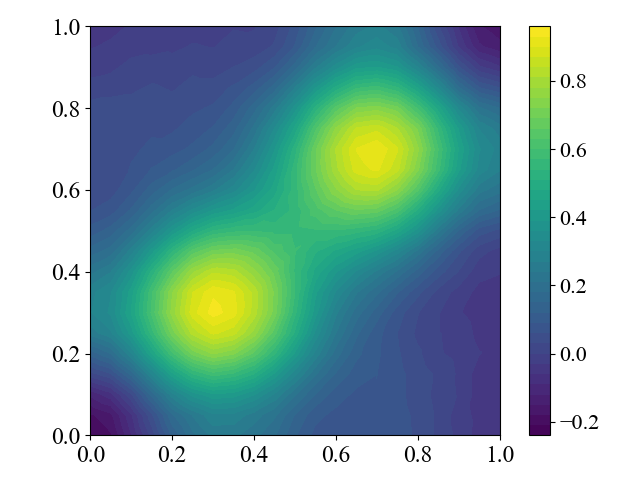}} 
	\\
	\subfloat[variance (pCN)]{
		\includegraphics[ keepaspectratio=true, width=0.31\textwidth, clip=true]{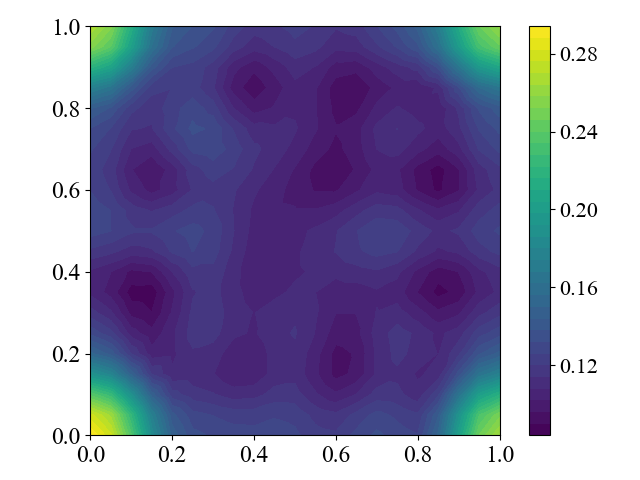}} 
	\subfloat[variance (planar)]{
		\includegraphics[ keepaspectratio=true, width=0.31\textwidth, clip=true]{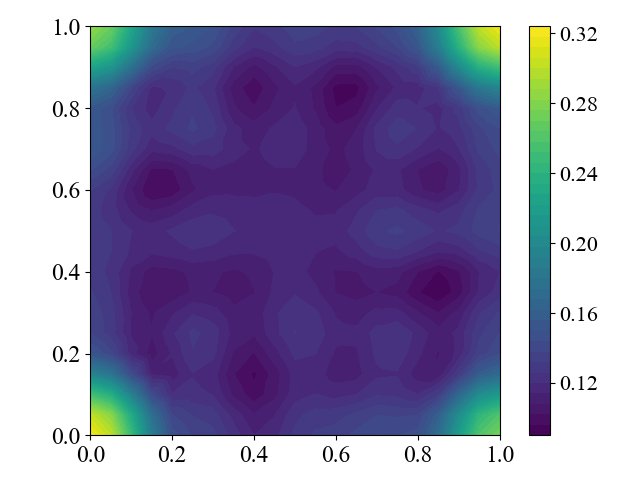}} 
        \subfloat[variance (Sylvester)]{
		\includegraphics[ keepaspectratio=true, width=0.31\textwidth, clip=true]{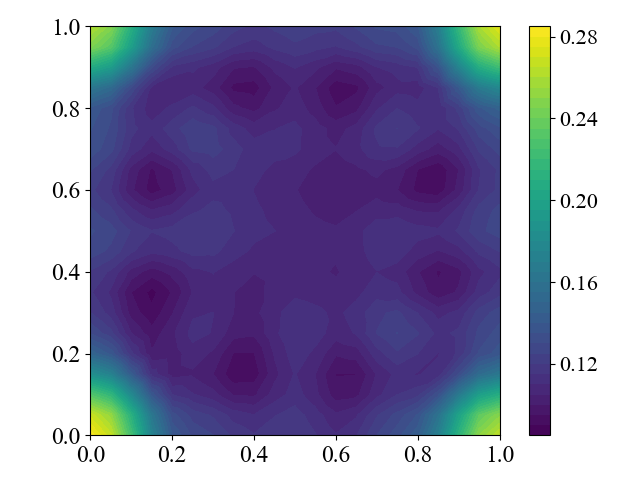}} 
	\caption{\emph{\small  The comparison of the measures obtained by flow models and pCN algorithm. (a)(b)(c): The mean function of posterior obtained by pCN algorithm, 32 layers planar flow, 5 layers Sylvester flow, respectively. (d)(e)(f): The variance function of posterior obtained by pCN algorithm, 32 layers planar flow, 5 layers Sylvester flow, respectively. }}
       \label{fig:nonlinearflow}
\end{figure}

Firstly, in Figure \ref{fig:nonlinearflow}, we present a comparative analysis of the mean function generated by the pCN algorithm and the approximate posterior obtained from the flow models. Subfigures (a), (b), and (c) of Figure \ref{fig:nonlinearflow} illustrate the mean function of the estimated posterior. The similarity between these results demonstrates the effectiveness of our proposed method in accurately estimating the mean function.

To further support the effectiveness of the NF-iVI method in estimating the mean function, we provide numerical evidence. The relative errors between the mean of the approximate posterior generated by 32 layers planar flow, 5 layers Sylvester flow and the mean of the posterior $u_p^*$ generated by the pCN algorithm are given by:
\begin{align}
    \textbf{planar relative error:}= \frac{\lVert u_{pN}^{*}-u_{p}^{*}\rVert^2}{\lVert u_{p}^{*}\rVert^2}=0.00142,\label{equ:relative error darcy1}\\
    \textbf{Sylvester relative error:}= \frac{\lVert u_{SN}^{*}-u_{p}^{*}\rVert^2}{\lVert u_{p}^{*}\rVert^2}=0.00188
    \label{equ:relative error darcy2},
\end{align}
where $u_{SN}^*$ is the mean of planar flow, and $u_{SN}^*$ is the mean of Sylvester flow. The estimated values are small, indicating the effectiveness of our algorithm.

Combining the visual evidence from Figure \ref{fig:nonlinearflow} and the numerical evidence provided by the relative error in Equations (\ref{equ:relative error darcy1}) and (\ref{equ:relative error darcy2}), we conclude that the NF-iVI method effectively approximates the mean of the posterior.

Secondly, to further explore the approximate posterior, we analyze the covariance functions. In subfigures (e)(f) of Figure \ref{fig:nonlinearflow}, we present the point-wise variance field of the posterior measure obtained from the NF-iVI method. For comparison, subfigure (d) of Figure \ref{fig:nonlinearflow} depicts the point-wise variance field estimated from pCN method. The similarity between these three plots indicates the effectiveness of the NF-iVI method in quantifying the uncertainties of the parameter $u$.

\begin{figure}[ht!]
	\centering
	\subfloat[variance (planar)]{
		\includegraphics[ keepaspectratio=true, width=0.31\textwidth, clip=true]{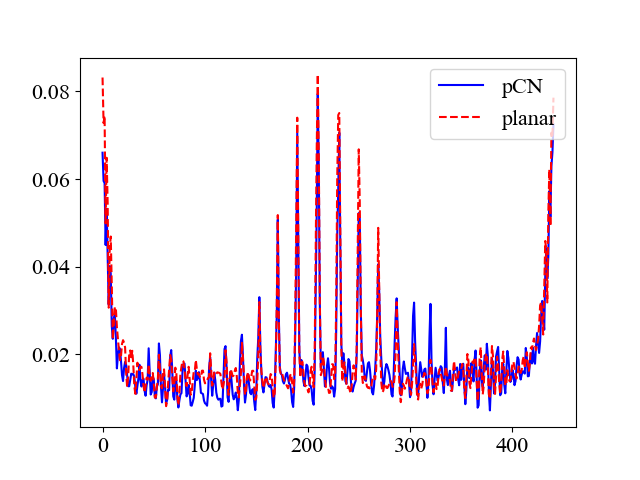}} 
	\subfloat[covariance (planar)]{
		\includegraphics[ keepaspectratio=true, width=0.31\textwidth, clip=true]{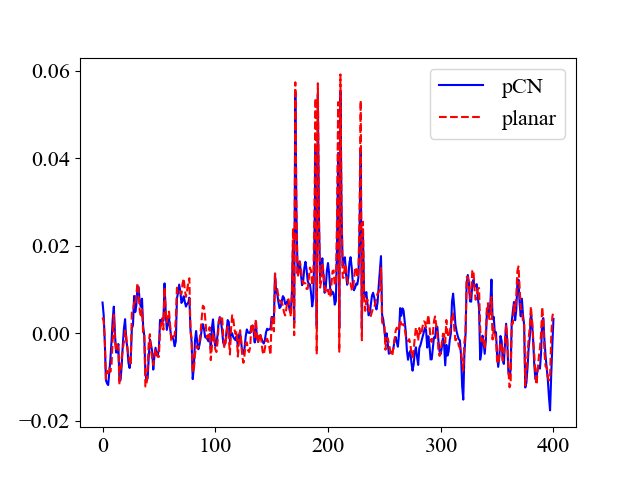}} 
        \subfloat[covariance (planar)]{
		\includegraphics[ keepaspectratio=true, width=0.31\textwidth, clip=true]{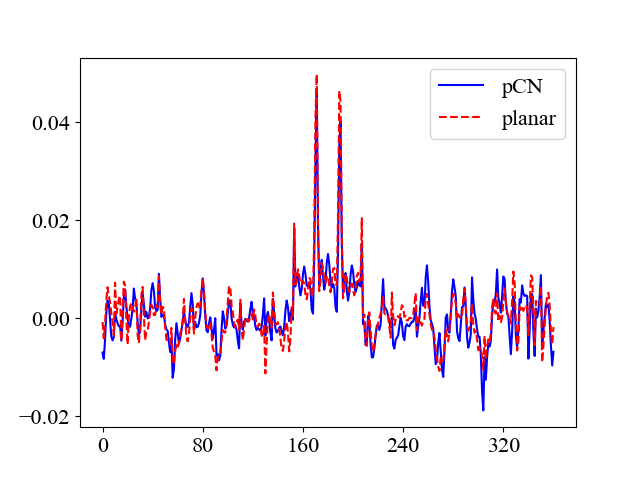}} 

        \subfloat[variance (Sylvester)]{
		\includegraphics[ keepaspectratio=true, width=0.31\textwidth, clip=true]{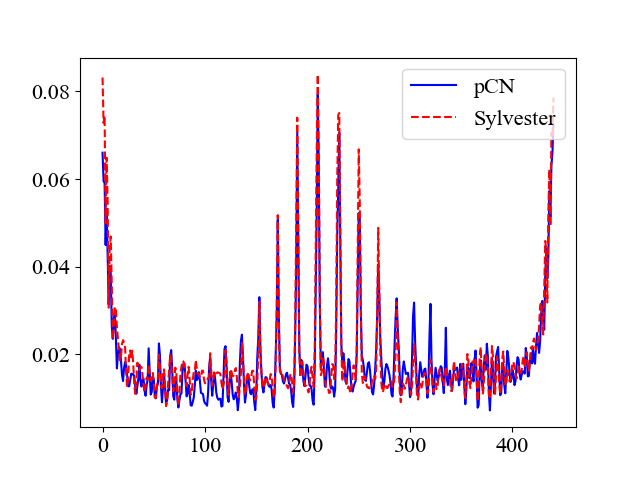}} 
	\subfloat[covariance (Sylvester)]{
		\includegraphics[ keepaspectratio=true, width=0.31\textwidth, clip=true]{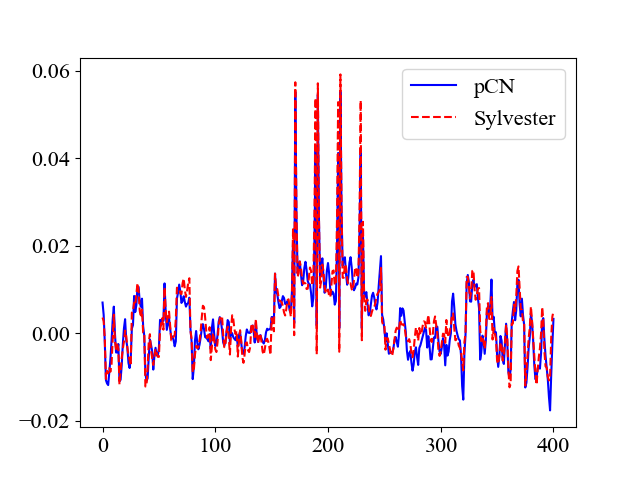}} 
        \subfloat[covariance (Sylvester)]{
		\includegraphics[ keepaspectratio=true, width=0.31\textwidth, clip=true]{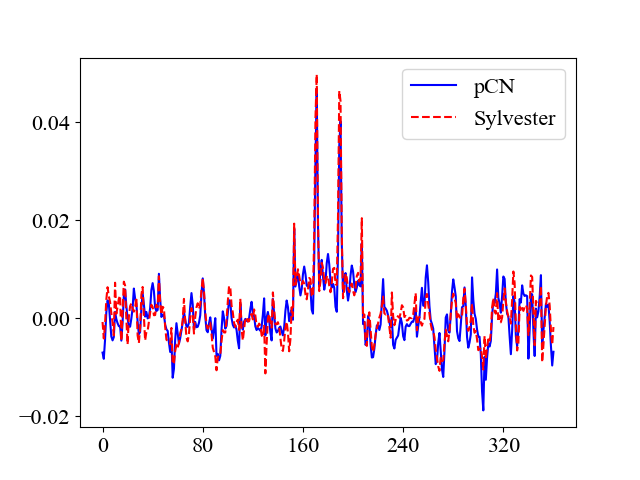}}

	\caption{\emph{\small {\color{black} The estimated variance and covariance functions obtained by the pCN algorithm (blue solid line), and NF-iVI algorithm (red dashed line).
				(a)(d): The covariance function $c(x,y)$ on all the mesh point pairs $\lbrace (x_i, x_{i})\rbrace^{N_g}_{i=1}$;
				(b)(e): The covariance function $c(x,y)$ on the mesh points $\lbrace (x_i, x_{i+40}) \rbrace^{N_g-40}_{i=1}$;
				(c)(f): The covariance function $c(x,y)$ on the mesh points $\lbrace (x_i, x_{i+80}) \rbrace^{N_g-80}_{i=1}$}}.}
       \label{pic:cov nonlinear}
\end{figure}

We also compare the covariance functions generated by these two methods. For simplicity, we compute these quantities at the mesh points and present the results in Figure \ref{pic:cov nonlinear}. In all subfigures of Figure \ref{pic:cov nonlinear}, the estimated covariance of posteriors obtained by the pCN and NF-iVI methods are represented by a blue solid line and a red dashed line, respectively. In subfigures (a)(d) of Figure \ref{pic:cov nonlinear}, we show the variance function calculated at all mesh points, i.e., $\{\mathrm{var}_{u}(x_{i})\}_{i=1}^{N_{g}}~(N_{g}$ is the number of mesh points, which is 441 in this problem). In subfigures (b)(e) and (c)(f) of Figure \ref{pic:cov nonlinear}, we show the covariance functions calculated on the pairs of points $\left\{(x_i,x_{i+40})\right\}_{i=1}^{N_g-40}$ and $\{(x_i,x_{i+80})\}_{i=1}^{N_g-80}$, respectively.

To numerically compare the variance and covariance functions generated by the two methods, we present the relative errors in Table \ref{table:relative nonlinear}. In the table, the notation $\bm{c}$ means the total relative error, and the notation $\lbrace c(x_i, x_{i+k}) \rbrace^{N_g-k}_{i=1}$ means the covariance function values on the pair of mesh points $\lbrace (x_i, x_{i+k}) \rbrace^{N_g-k}_{i=1}$ ($k=0, 40, 80$). The numbers below this notation are the relative errors between the vectors obtained by pCN method and NF-iVI, respectively. The small relative errors indicate that the posterior covariance functions obtained by the NF-iVI and pCN sampling methods are quantitatively similar. Combining the visual evidence from Figures \ref{fig:nonlinearflow} and \ref{pic:cov nonlinear}, and the numerical evidence from Table \ref{table:relative nonlinear}, we conclude that the NF-iVI method effectively approximates the covariance function of the posterior.

\begin{table}[ht!]
	\renewcommand{\arraystretch}{1.5}
	\centering
		\caption{\emph{\small The relative errors between the variance function and covariance functions.}} 
	\begin{tabular}{c|cccc}
		\hline $\text{Relative Error} $& $\bm{c}$  & $\lbrace c(x_i, x_i)\rbrace^{N_g}_{i=1}$& $\lbrace c(x_i, x_{i+40}) \rbrace^{N_g-40}_{i=1}$& $\lbrace c(x_i, x_{i+80}) \rbrace^{N_g-80}_{i=1}$\\
		\hline $\text{Planar Flow}$& $0.07252$& $0.00889$& $0.01254$& $0.03345$\\
        \hline $\text{Sylvester Flow}$ & $0.09725$& $0.01227$& $0.02312$& $0.05118$\\
		\hline
	\end{tabular}
 \label{table:relative nonlinear}
\end{table}

As discussed in Subsection \ref{subsec:nonleaner transformation}, the functional Sylvester flow is anticipated to exhibit superior approximation capabilities compared with the functional planar flow. To demonstrate this, we consider the same inverse problem defined in equation (\ref{equ:nonlinear}). In this experiment, we set the dimension of $\text{Im}(\mathcal{F}_n)$ in functional Sylvester flow to 20 ($M=20$ in Subsection \ref{subsec:nonleaner transformation}). By employing only 5 layers of functional Sylvester flow, we achieve satisfactory results in approximating the posterior. In contrast, using 5 layers of functional planar flow would be insufficient to achieve similar performance.

\begin{figure}[ht!]
	\centering
	
	\subfloat[background truth]{
		\includegraphics[ keepaspectratio=true, width=0.31\textwidth, clip=true]{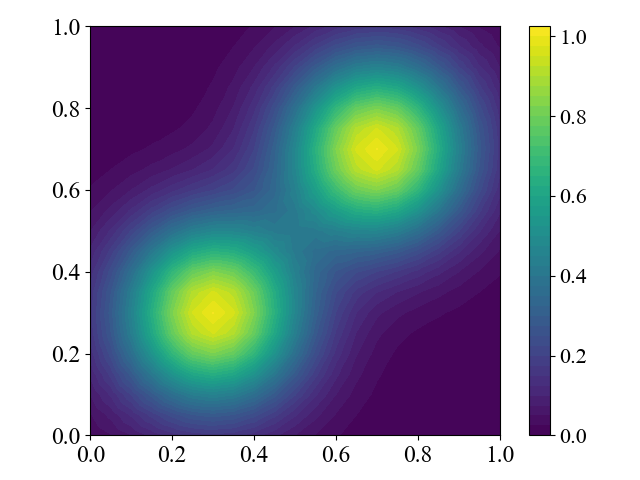}} 
	\subfloat[mean (Sylvester-5)]{
		\includegraphics[ keepaspectratio=true, width=0.31\textwidth, clip=true]{PIC/nonlinear/s_mean.png}} 
        \subfloat[variance (Sylvester-5)]{
		\includegraphics[ keepaspectratio=true, width=0.31\textwidth, clip=true]{PIC/nonlinear/s_pointwise.png}} 
	\\
	\subfloat[background truth]{
		\includegraphics[ keepaspectratio=true, width=0.31\textwidth, clip=true]{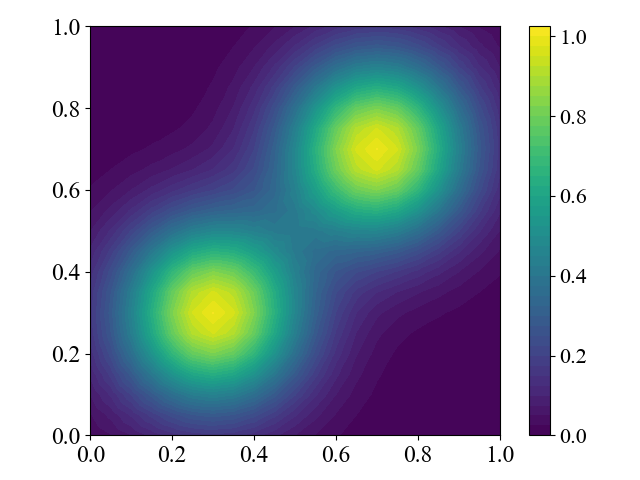}} 
	\subfloat[mean (planar-5)]{
		\includegraphics[ keepaspectratio=true, width=0.31\textwidth, clip=true]{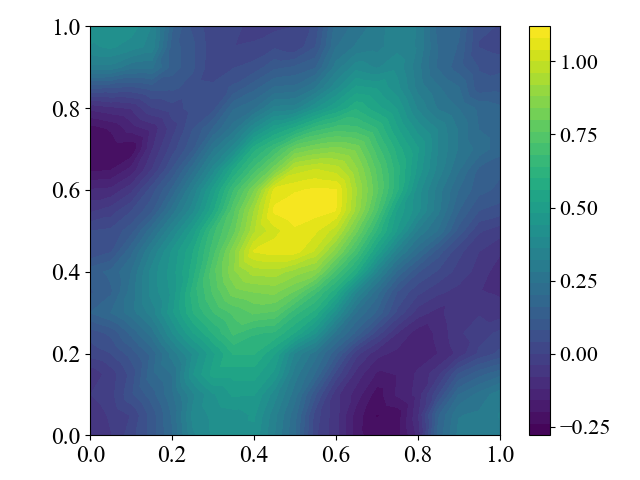}} 
        \subfloat[variance (planar-5)]{
		\includegraphics[ keepaspectratio=true, width=0.31\textwidth, clip=true]{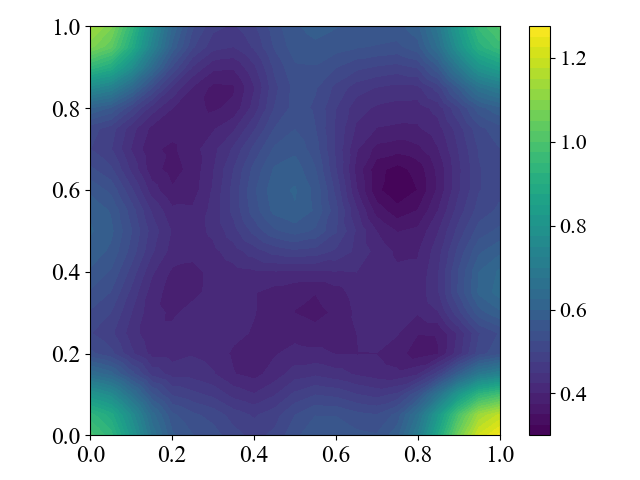}} 
	\caption{\emph{\small The comparison of approximation ability between functional planar flow and functional Sylvester flow. We use 5 layers for both flow models to approximate the posterior. (a)(d): The background truth of $u$. (b):  The mean of approximate posterior obtained by functional Sylvester flow. (c): The variance function of approximate posterior obtained by functional Sylvester flow. (e): The mean function of approximate posterior obtained by functional planar flow. (f):The variance function of approximate posterior obtained by functional planar flow.}}
	
       \label{fig:sflow}
	
\end{figure}

Subfigures (a) and (d) of Figure~\ref{fig:sflow} show the ground truth. Subfigures (b), (c), (e), and (f) depict the mean functions and the estimated pointwise variance fields of the approximate posterior obtained using the 5-layer functional Sylvester flow and the 5-layer functional planar flow. The comparisons highlight the limitations of functional planar flow in achieving accurate approximations with a small number of layers. In contrast, functional Sylvester flow exhibits superior approximation capabilities.

\subsection{Electrical Impedance Tomography}
{\color{black}
In this subsection, we focus on the electrical impedance tomography (EIT) problem. EIT is a technique to see inside an object by measuring electricity on its surface. It works by attaching electrodes, sending a small electrical current through them, and then reading the voltages that appear. We utilize the complete electrode
model (CEM), proposed in \cite{somersalo1992existence}. The strong form of the PDE governing this model is given as follows:
\begin{align*}
\begin{cases} 
-\nabla \cdot (\kappa \nabla v(x)) = 0, & x \in \Omega, \\ 
\int_{e_l} \kappa \frac{\partial v}{\partial n}  dS = I_l, & l = 1, \ldots, L, \\ 
\kappa \frac{\partial v}{\partial n}(x) = 0, & x \in \partial \Omega \setminus \bigcup_{l=1}^L e_l, \\ 
v(x) + z_l \kappa \frac{\partial v}{\partial n}(x) = V_l, & x \in e_l, l = 1, \ldots, L,
\end{cases}
\end{align*}
where $\Omega \subset \mathbb{R}^2$ is the domain and $\{e_l\}_{l=1}^L \subset \partial \Omega$ are electrodes on the boundary upon which currents $I_{\textbf{all}}=\{I_l\}_{l=1}^L$ are injected and voltages $\{V_l\}_{l=1}^L$ are read. The constants $\{z_l\}_{l=1}^L$ represent the contact impedances of the electrodes. The function $\kappa \in \mathcal{H}_{\kappa}$ represents the conductivity of the body and $v \in \mathcal{H}_v$ represents the potential within the body. 

The inverse problem we consider is the recovery of $\kappa \in \mathcal{H}_{\kappa}$ from a sequence of boundary voltage measurements on the electrodes. The mathematical formulation of the inverse problem is given as follows.

Following \cite{somersalo1992existence}, to improve the recovery of the model parameter $\kappa$, a standard strategy is to apply multiple distinct currents injection $I_{\mathrm{all}}^{(1)}, I_{\mathrm{all}}^{(2)}, \dots, I_{\mathrm{all}}^{(M)}$ and record the corresponding electrode voltages, thereby collecting richer measurement information. Under a single stimulation pattern \( I_{\mathrm{all}}^{(i)} = \{I^{(i)}_l\}_{l=1}^L \), the solution of the PDE consists of the electric potential field $v^{(i)}$ inside the conductor and the set of voltage data $\{V_l^{(i)}\}_{l=1}^{L}$ on the $L$ electrodes, respectively. Thus, the PDE solution operator $\mathcal{G}^{(i)}$: $\mathcal{H}_{\kappa} \rightarrow \mathcal{H}_{v} \times \mathbb{R}^L$ can be written as
$$
\mathcal{G}^{(i)}(\kappa)=\left(v^{(i)}, V^{(i)}_1, V^{(i)}_2, \ldots, V^{(i)}_L\right).
$$
In practice, direct measurement of the internal voltage $v^{(i)}$  is challenging, and only the electrode voltages \( \{V^{(i)}_l\}_{l=1}^L \) are accessible. Hence the measurement operator $\mathcal{S}^{(i)}$: $\mathcal{H}_{v} \times \mathbb{R}^L \rightarrow \mathbb{R}^L$ is defined as

$$
\mathcal{S}^{(i)}\left(v^{(i)}, V^{(i)}_1, V^{(i)}_2, \ldots, V^{i}_L\right)=\left(V^{(i)}_1, V^{(i)}_2, \ldots, V^{(i)}_L\right).
$$
We perform \( M \) independent current stimulations, and the complete set of measured electrode voltages is therefore  

\[
V = \left( V^{(1)}, V^{(2)}, \dots, V^{(M)} \right),
\]
where $V^{(m)}=\{V_l^{(m)}\}_{l=1}^L$ is the boundary voltage data of the $m$-th current stimulation, $m=1, ..., M$. Since $M$ current stimulations are performed, the inverse problem can be written as
\begin{align}\label{equ:eit}
    \bm{d}=\mathcal{S}\mathcal{G}(\kappa)+\bm{\epsilon},
\end{align}
where $\bm{\epsilon}$ is the noise, $\mathcal{G}$ is a collection of PDE solution operators mapping from $\kappa$ to
their corresponding solutions:
$$
\mathcal{G}(\kappa)=\left(\mathcal{G}^{(1)}(\kappa), \mathcal{G}^{(2)}(\kappa), \ldots, \mathcal{G}^{(M)}(\kappa) \right),
$$
and $\mathcal{S}$ is the measurement operator:
$$
\mathcal{S}\mathcal{G}(\kappa)=\left(\mathcal{S}^{(1)}\mathcal{G}^{(1)}(\kappa), \mathcal{S}^{(2)}\mathcal{G}^{(2)}(\kappa), \ldots, \mathcal{S}^{(M)}\mathcal{G}^{(M)}(\kappa) \right).
$$

\begin{figure}[ht!]
	\centering
	\subfloat[mesh]{
		\includegraphics[ keepaspectratio=true, width=0.32\textwidth, clip=true]{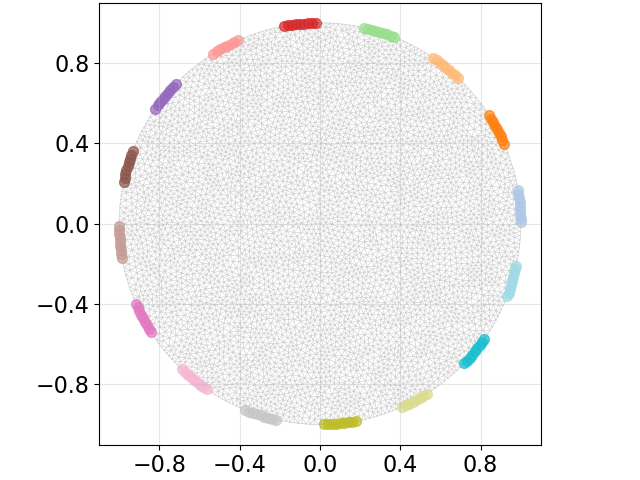}} 
    \subfloat[measurement data]{
		\includegraphics[ keepaspectratio=true, width=0.32\textwidth, clip=true]{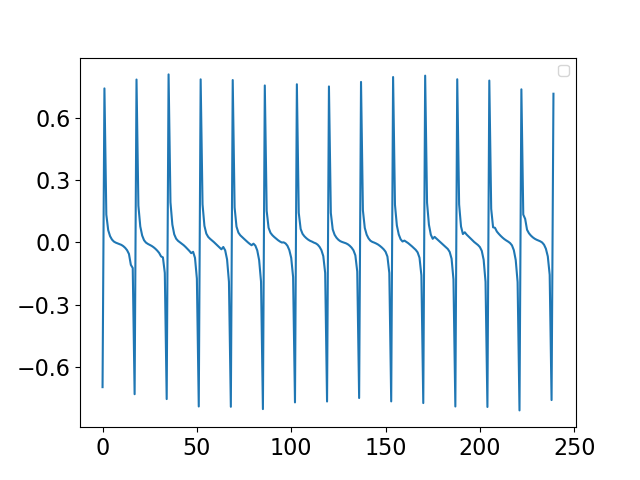}} 
	\subfloat[measurement data]{
		\includegraphics[ keepaspectratio=true, width=0.32\textwidth, clip=true]{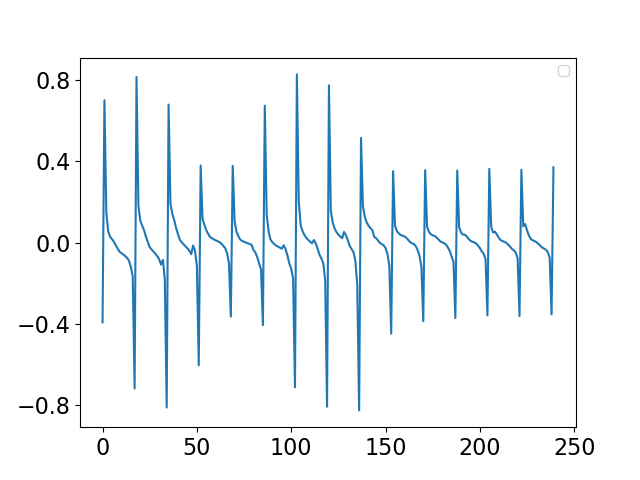}} \\
	
	\caption{\emph{\small (a) Mesh and electrodes. (b)(c) Measurement data $\bm{d}$ with noise $\bm{\epsilon}$.
}}
\label{pic:eit_data}	
\end{figure}

For the numerical experiments conducted in this study, we specify the computational domain as $ \Omega = \{x \in \mathbb{R}^2\mid \|x\| < 1\} $, where $ \|\cdot\| $ denotes the Euclidean norm in $ \mathbb{R}^2 $. Sixteen electrodes are equally spaced along the boundary $ \partial\Omega $, achieving a 50\% boundary coverage rate (see Subfigure (a) of Figure \ref{pic:eit_data}). All contact impedances are taken to be \( z_l = 0.025 \). Adjacent electrodes are stimulated with a current of 1, so that the matrix of stimulation patterns \( I_{\textbf{total}} \in \mathbb{R}^{16 \times 15} \) ($I_{\textbf{total}}(l,i)=I_l^{(i)}$) is given by
\[
I = 1 \times 
\begin{pmatrix}
+1 & 0 & \cdots & 0 \\
-1 & +1 & \cdots & 0 \\
0 & -1 & \ddots & 0 \\
\vdots & \vdots & \ddots & +1 \\
0 & 0 & 0 & -1
\end{pmatrix}.
\]
{\color{black}
In this numerical experiment, consistent with much of the classical literature on EIT inverse problems \cite{karimi2021optimal, dunlop2017hierarchical}, we assume that the background truth $\kappa^{\dagger}$ is a piecewise constant function. The direct application of the smooth Gaussian prior $\mu_0$, as employed in Subsections \ref{subsec4.1} and \ref{subsec:darcy}, is ill-suited to this setting owing to a regularity mismatch: the non-smooth character of $\kappa^{\dagger}$ is incompatible with the smooth prior assumed in the inversion. Although several alternatives exist (including Besov priors \cite{Dashti2012IPI} and TV-Gaussian priors \cite{tvgaussian}), we adopt the level set method \cite{burger2001level, iglesias2016bayesian, burger2005survey} for its intuitive appeal and simplicity.

The level set method is widely employed in PDE-constrained inverse problems where the model parameter is piecewise constant. In this experimrnt, for observation model \eqref{equ:eit}, we assume $\kappa$ takes only two distinct values, denoted as $a$ and $b$. A latent variable $\hat{\kappa}$ is introduced via the level set mapping
\begin{equation}
\kappa(x)=\mathcal{F}_{\text{LS}}\bigl(\hat\kappa(x)\bigr), \quad 
\mathcal{F}_{\text{LS}}(t)=
\begin{cases} 
a, & t \leq 0, \\[4pt]
b, & t>0,
\end{cases}
\end{equation}
where we set $a=1$ and $b=3$ for this study. This link function relates the discontinuous parameter $\kappa$ to the smooth latent variable $\hat{\kappa}$, thereby enabling the use of the smooth Gaussian prior $\mu_0$ for Bayesian inversion. Substituting this mapping into the inverse problem yields the revised observation model
\begin{equation}\label{equ:LS}
\bm{d}=\mathcal{S}\mathcal{G}\big(\mathcal{F}_{\text{LS}}(\hat{\kappa})\big)+\bm{\epsilon}.
\end{equation}
Inference is performed on the smooth latent variable $\hat{\kappa}$, and the final parameter estimate is obtained via $\kappa_{\text{mean}}=\mathcal{F}_{\text{LS}}(\hat{\kappa}_{\text{mean}})$.

The experiment is conducted on the unit disk $\Omega$. We consider two scenarios to generate the ground-truth conductivity. 
For the first scenario $\kappa_1^{\dagger}(x)$, we define the subdomain of high-conductivity anomalies $T$ as
\[
T=\left\{(x_1,x_2)\in\Omega \;\bigg|\;
\exp\left[-3\big(2(x_1-0.5)^2+x_2^2\big)\right]
+\exp\left[-3\big(2(x_1+0.4)^2+x_2^2\big)\right]>0.5\right\}.
\]
Let $\chi_T$ denote the indicator function of $T$, i.e., 
\[ 
\chi_T(x) = \begin{cases} 1, & \text{if } x \in T, \\ 0, & \text{if } x \in \Omega \setminus T. \end{cases} 
\]
 The ground truth is constructed as $\kappa_1^{\dagger}(x)=1+2\chi_T(x_1, x_2)$, which corresponds to the values $a=1$ and $b=3$ defined above.
For the second scenario, we generate the ground truth $\kappa_2^{\dagger}(x)$ by sampling a smooth latent function $\hat{\kappa}_2(x)$ directly from the Gaussian prior $\mu_0$ and mapping it using the level set mapping: $\kappa_2^{\dagger}(x)=\mathcal{F}_{\text{LS}}(\hat{\kappa}_2(x))$. 
We generate synthetic observation data $\bm{d}_1$ and $\bm{d}_2$ corresponding to $\kappa_1^{\dagger}$ and $\kappa_2^{\dagger}$ using the observation model \eqref{equ:eit}. The underlying ground truths are depicted in Subfigures (a) and (e) of Figure \ref{pic:eit}, and the corresponding measurement data are shown in Subfigures (b) and (c) of Figure \ref{pic:eit_data}. We then train the proposed functional normalizing flow model within the Bayesian inversion framework using $\bm{d}_1$ and $\bm{d}_2$ to recover the latent parameters $\hat{\kappa}_1$ and $\hat{\kappa}_2$, thereby approximating the true parameters $\kappa_1^{\dagger}$ and $\kappa_2^{\dagger}$.

 It is noteworthy that the construction of the level set mapping $\mathcal{F}_{\text{LS}}$ relies on prior knowledge of the two constant values $a$ and $b$. Consequently, level set method cannot be directly applied if these values are unknown. When only the piecewise-constant structure of $\kappa$ is available while $a$ and $b$ remain undetermined, standard level set inversion becomes infeasible. As a remedy, these constant values can be inferred directly from measurement data, as demonstrated in \cite{akyildiz2025efficient}. While there are methods exist to infer these values, this paper focuses on constructing the functional normalizing flow model and approximating Bayesian posteriors. Thus, following common practice in level set inversion \cite{burger2001level, iglesias2016bayesian, burger2005survey, dunlop2017hierarchical}, we assume $a$ and $b$ are known a priori and directly define the level set link function $\mathcal{F}_{\text{LS}}$.

Upon applying the level set method, the inverse problem is reformulated as \eqref{equ:LS}. For clarity, we list the specific choices for the parameters introduced in this problem as follows:
}
\begin{itemize}
\item Assume that 1$\%$ random Gaussian noise $\bm{\epsilon} \sim \mathcal{N}(0, \bm{\Gamma}_{\text{noise}})$ is added, where $\bm{\Gamma}_{\text{noise}} = \tau^{-1}\textbf{I}$, and $\tau^{-1} = (0.01 \lVert \mathcal{S}\mathcal{G}(\mathcal{F}_{\text{LS}}(\hat{\kappa}))\rVert_{\infty})^2$.
\item The covariance operator \(\mathcal{C}_0\) associated with the prior measure \(\mu_0\) is defined as \(\mathcal{C}_0 = (\text{I} - \alpha \Delta)^{-2}\), where \(\alpha = 0.1\) denotes a pre-specified constant. Here, the Laplace operator \(\Delta\) is defined on the computational domain \(\Omega\) subject to homogeneous Neumann boundary conditions. Furthermore, the prior measure \(\mu_0\) is taken to have a zero mean.
\item To avoid inverse crime, the synthetic data is generated on a fine mesh consisting of $20054$ grid points, while a mesh of distinct resolution ($n = 5124$ grid points) is adopted in the inversion stage.
\end{itemize}

\begin{figure}[ht!]
	\centering
	
	\subfloat[background truth]{
		\includegraphics[ keepaspectratio=true, width=0.25\textwidth, clip=true]{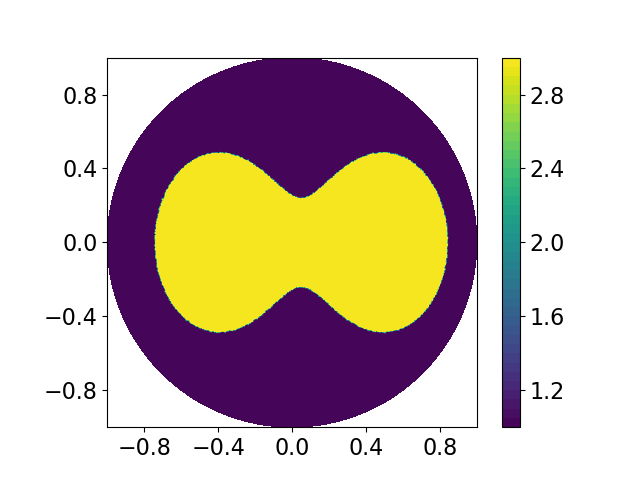}} 
	\subfloat[mean of latent $\hat{\kappa}_{\text{mean}}$]{
		\includegraphics[ keepaspectratio=true, width=0.25\textwidth, clip=true]{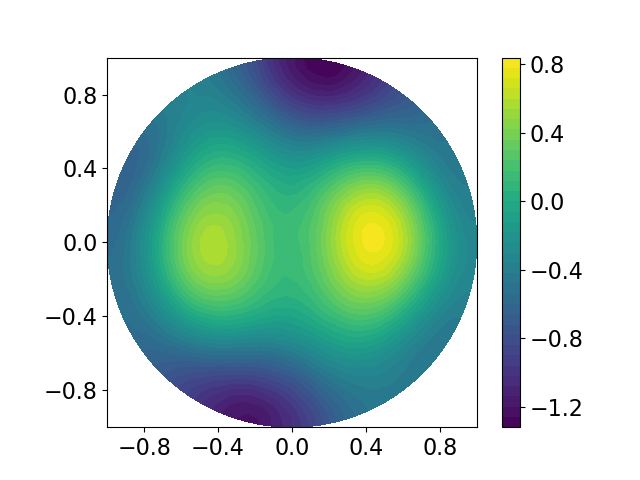}} 
        \subfloat[$\mathcal{F}_{\text{LS}}(\hat{\kappa}_{\text{mean}})$]{
		\includegraphics[ keepaspectratio=true, width=0.25\textwidth, clip=true]{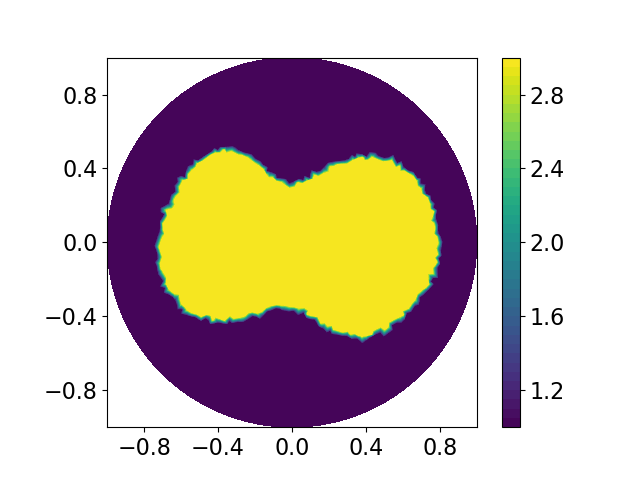}} 
	\subfloat[variance (Sylvester)]{
		\includegraphics[ keepaspectratio=true, width=0.25\textwidth, clip=true]{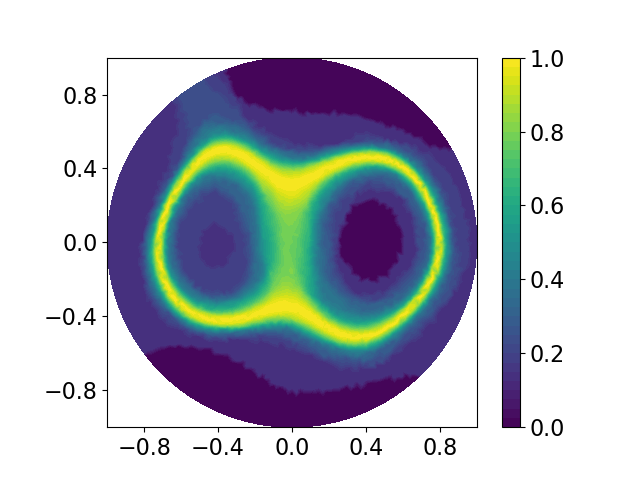}}
 \\
	\subfloat[background truth]{
		\includegraphics[ keepaspectratio=true, width=0.25\textwidth, clip=true]{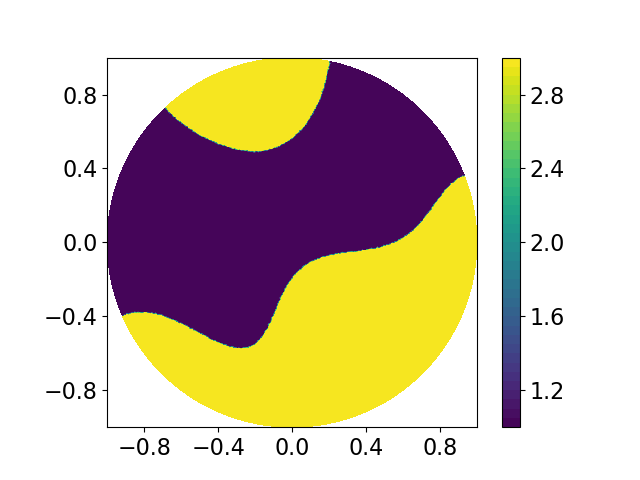}} 
        \subfloat[mean of latent $\hat{\kappa}_{\text{mean}}$]{
		\includegraphics[ keepaspectratio=true, width=0.25\textwidth, clip=true]{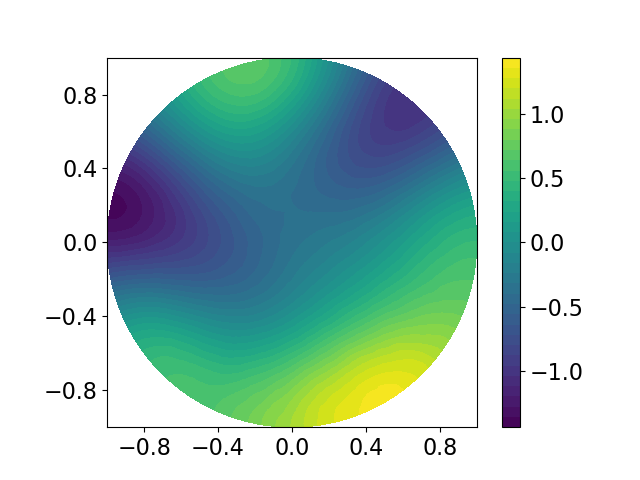}} 
        \subfloat[$\mathcal{F}_{\text{LS}}(\hat{\kappa}_{\text{mean}})$]{
		\includegraphics[ keepaspectratio=true, width=0.25\textwidth, clip=true]{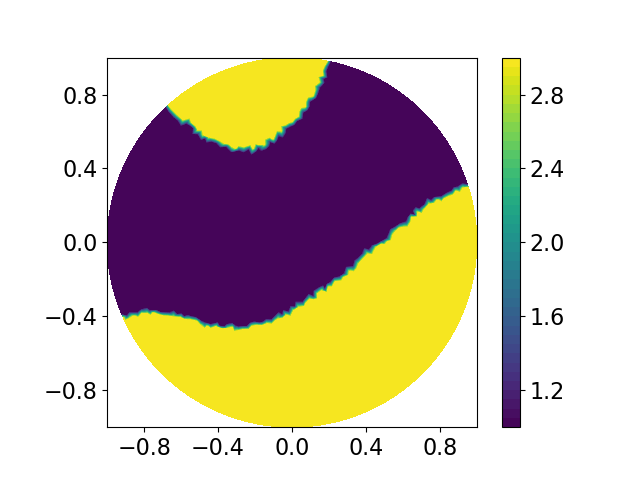}}
        \subfloat[variance (Sylvester)]{
		\includegraphics[ keepaspectratio=true, width=0.25\textwidth, clip=true]{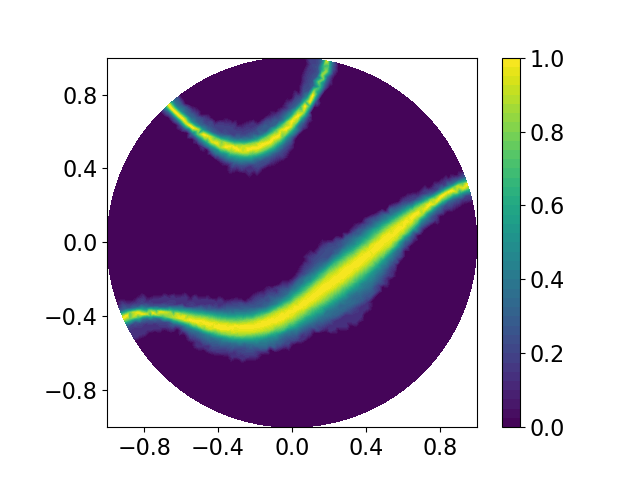}} 
	\caption{\emph{\small  (a)(e) Background truth. (b)(f) Mean of the approximate latent posterior $\hat{\kappa}_{\text{mean}}$. (c)(g) The result of $\mathcal{F}_{\text{LS}}(\hat{\kappa}_{\text{mean}})$, (d)(h) Pointwise variance of the approximate posterior. }}
       \label{pic:eit}
\end{figure}

Owing to the nonlinear nature of the inverse problem, we employ a five-layer functional Sylvester flow model to approximate the posterior distribution about $\hat{\kappa}$. {\color{black}It is worth noting that the level set mapping employed herein presents a numerical drawback that complicates gradient-based training \cite{akyildiz2025efficient}. The original operator $\mathcal{F}_{\text{LS}}$ yields vanishing gradients almost everywhere, thereby impeding gradient descent optimization during the training of the functional normalizing flow model.
 To circumvent this difficulty, we substitute $\mathcal{F}_{\text{LS}}$ with a smoothed surrogate during training, given by
\[
\mathcal{F}_{\text{LS}}^{\text{smooth}}(t)
= a + \frac{b-a}{1 + \exp\left(-10 t\right)}.
\]
In practice, $\mathcal{F}_{\text{LS}}^{\text{smooth}}$ is employed for all gradient computations throughout the functional normalizing flow training process, while the original discontinuous mapping $\mathcal{F}_{\text{LS}}$ is utilized for the final reconstruction and quantitative evaluation of $\kappa$. Such smoothing is a standard remedy for level set inverse problems \cite{akyildiz2025efficient}.}

Following Algorithm \ref{alg A}, the flow model is trained for \(K = 2000\) iterations, with each iteration drawing \(N = 10\) samples. A step-decay learning rate schedule is adopted: the initial learning rate is set to \(\alpha_0 = 0.01\), and it is multiplied by a decay factor \(\tau = 0.8\) every 200 iterations. The application of the pCN algorithm to the EIT problem resulted in prohibitive computational cost. Consequently, solving the EIT problem with the pCN algorithm was omitted in our experiment.

{\color{black}Subfigures (b)--(d) and (f)--(h) of Figure \ref{pic:eit} display the result generated by the trained functional Sylvester flow. Observing Subfigures (b) and (f), the reconstruction $\hat{\kappa}$ appears excessively smooth, whereas $\kappa=\mathcal{F}_{\text{LS}}(\hat{\kappa})$ recovers the true profile, as shown in Subfigures (a)(c) and (e)(g). Furthermore, Subfigures (d) and (h) depict the covariance functions of the estimated posteriors. Here, to obtain this covariance function, we first sample from the trained functional normalizing flow model to generate the sample set $\{\hat{\kappa}_i\}_{i=1}^N$. Subsequently, we apply the level set function $\mathcal{F}_{\text{LS}}$ to obtain the transformed sample set $\{\mathcal{F}_{\text{LS}}(\hat{\kappa}_i)\}_{i=1}^N$. Finally, the covariance is computed via these transformed samples according to the formulas presented in Subsection \ref{subsec4.1}. Collectively, these numerical experiments demonstrate the efficacy of our proposed algorithm for the EIT problem.
}

}
\vskip 0.2 cm

\section{Numerical Examples of Conditional Functional Normalizing Flow}\label{sec5}
In this section, we present the results of the conditional functional normalizing flow, as described in Section~\ref{sec3}. For a given measurement pair $(\bm{x}^*,\bm{d}^*)$, the trained network $\mathcal{N}_{\lambda}(\bm{v})$ can efficiently provides a rough estimate of the corresponding posterior distribution $\mu_{(\bm{x}^*,\bm{d}^*)}$. In addition, for a specific measurement pair, the results can be further refined through a retraining process.

\subsection{Conditional Functional Normalizing Flow}\label{sec5.1}
In this subsection, we detail the training process of the conditional network and the extraction of an initial posterior approximation directly from it. Considering the steady-state Darcy flow equation (\ref{equ:darcy}), the inverse problem aims to infer the posterior distribution of the parameter $u$, given measurement data $\bm{d} = (w(x^1), \ldots, w(x^{N_d}))$ at points $\bm{x} = (x^1, \ldots, x^{N_d})$. To approximate the posterior, we employ a 5-layer functional Sylvester flow similar with Subsection \ref{subsec:darcy}. 

\begin{figure}[ht!]
	\centering
	\subfloat[sample 1]{
		\includegraphics[ keepaspectratio=true, width=0.32\textwidth, clip=true]{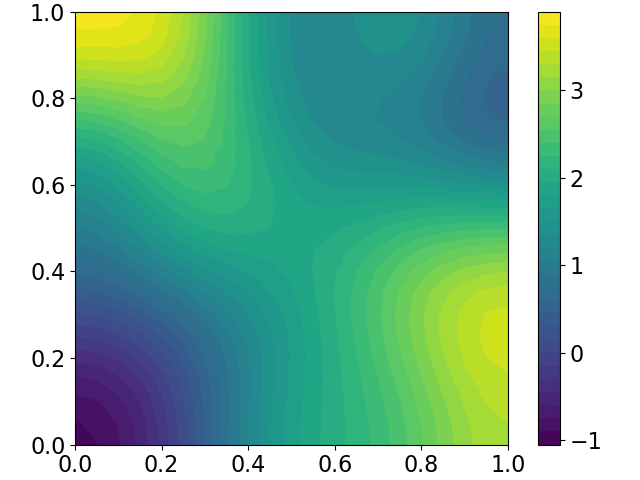}} 
    \subfloat[sample 2]{
		\includegraphics[ keepaspectratio=true, width=0.32\textwidth, clip=true]{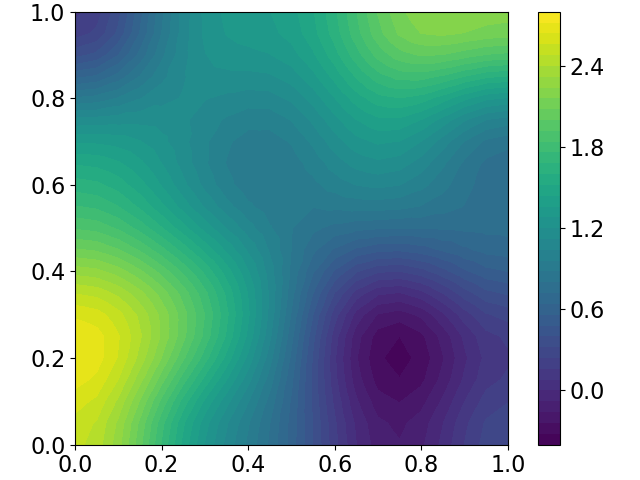}} 
	\subfloat[sample 3]{
		\includegraphics[ keepaspectratio=true, width=0.32\textwidth, clip=true]{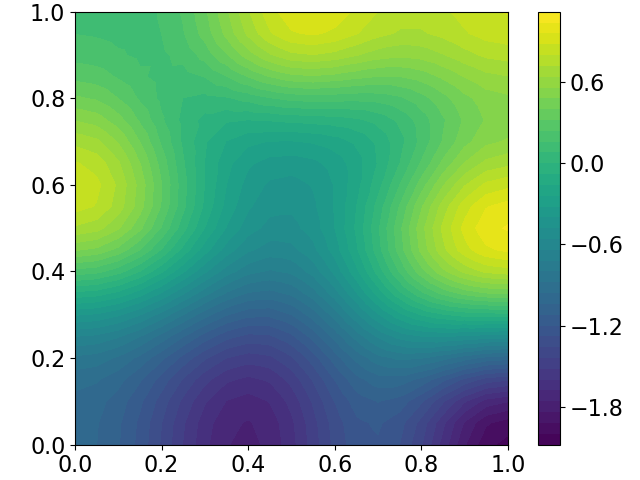}} \\
	
	\caption{\emph{\small Samples from $\mu_0$.
}}
\label{pic:samples}	
\end{figure}

To effectively train the conditional network, we generate the training dataset based on the following process:

\begin{itemize}
	\item {\color{black} Draw a large number of samples from the measure $\mu_0$, denoted as $\{u_1,u_2,\ldots,u_{N_{train}}\}$, where $\mu_0$ is the prior measure we defined in Subsection \ref{subsec:darcy}. Examples of the samples drawn from $\mu_0$ are shown in Figure \ref{pic:samples}.}
	\item Use these samples as parameters for the PDE equation (\ref{equ:darcy}), and generate their corresponding solutions $ \{ w_1, w_2,\ldots, w_{N_{train}}\}$.
	\item Determine the measurement points $\{\bm{x}_1,\bm{x}_2,\ldots,\bm{x}_{N_{train}}\}$ for each solution and calculate their corresponding measurement data $\{\bm{d}_1,\bm{d}_2,\ldots\bm{d}_{N_{train}}\}$.
        \item Combine the measurement points and data into pairs $(\bm{x}_i, \bm{d}_i)$ and form the training dataset $D_{train}=\{(\bm{x}_1,\bm{d}_1),(\bm{x}_2,\bm{d}_2),\ldots,(\bm{x}_{N_{train}},\bm{d}_{N_{train}})\}$.
\end{itemize}

Here, we briefly describe the symbols used in this context. The vector $\bm{x}=(x^1, x^2 , \ldots, x^{\bm{N}_d})$ represents the measurement points, where each $x^i$ corresponds to a distinct measurement point. The notation $\bm{x}_j$ refers to different measurement points associated with different measurement pairs $(\bm{x}_j, \bm{d}_j)$. For convenience, all the measurement points \(\bm{x}_1, \ldots, \bm{x}_{N_{\text{train}}}\) are selected uniformly across the domain. Their coordinates are given by  $\{(i/21,j/21)\}_{i,j=1}^{20}$.

\begin{figure}[ht!]
	\centering
	\subfloat[background truth]{
		\includegraphics[ keepaspectratio=true, width=0.33\textwidth, clip=true]{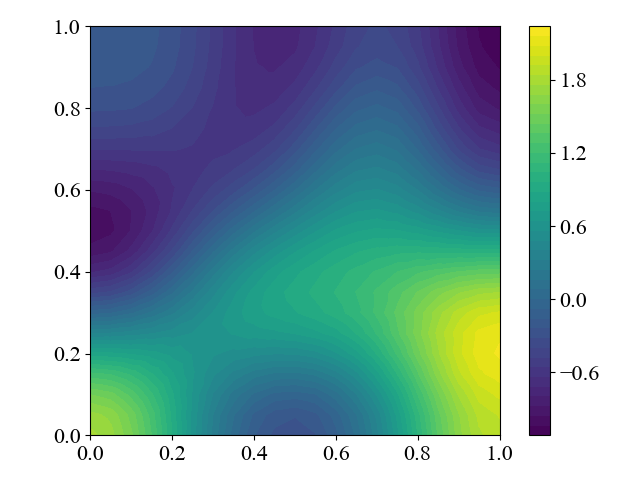}} 
	\subfloat[mean (CNF-iVI)]{
		\includegraphics[ keepaspectratio=true, width=0.33\textwidth, clip=true]{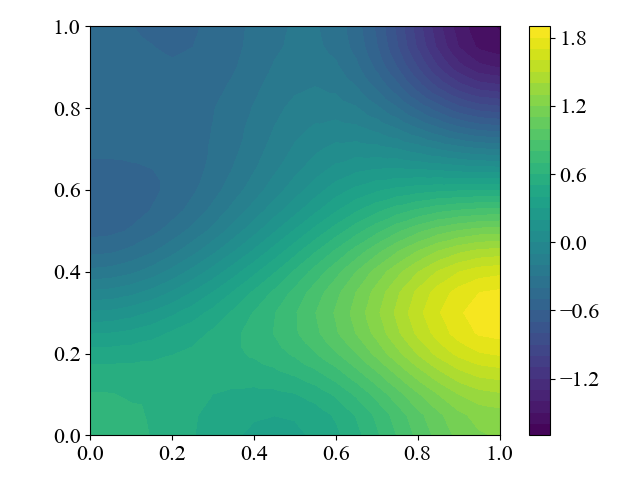}} 
        \subfloat[mean (pCN)]{
		\includegraphics[ keepaspectratio=true, width=0.33\textwidth, clip=true]{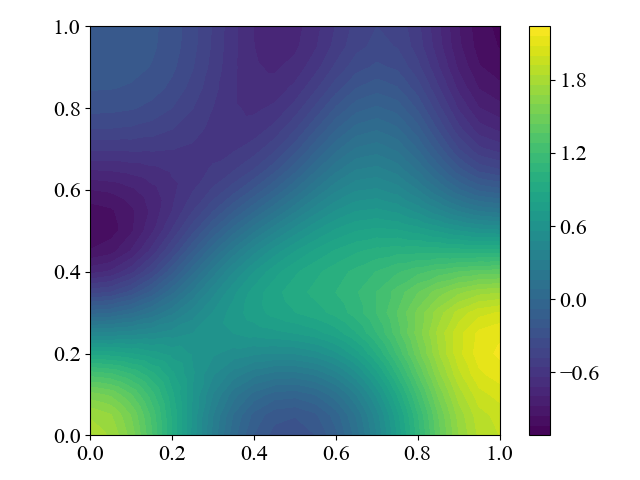}} \\
  \subfloat[variance]{
		\includegraphics[ keepaspectratio=true, width=0.31\textwidth, clip=true]{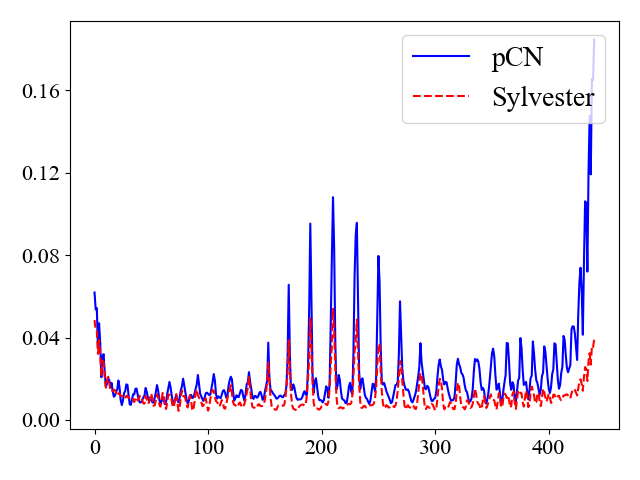}} 
	\subfloat[covariance]{
		\includegraphics[ keepaspectratio=true, width=0.31\textwidth, clip=true]{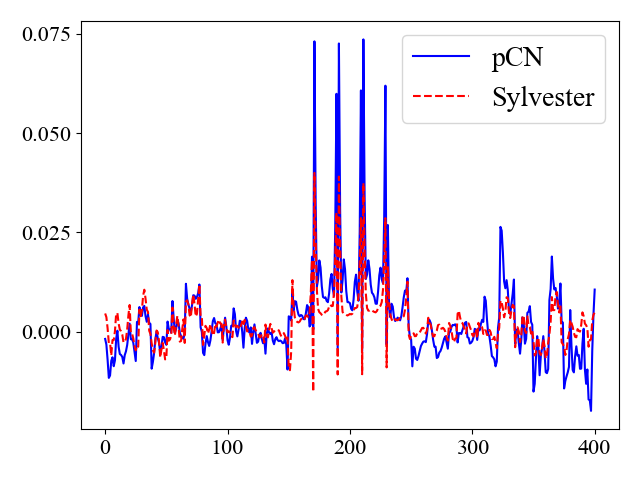}} 
        \subfloat[covariance]{
		\includegraphics[ keepaspectratio=true, width=0.31\textwidth, clip=true]{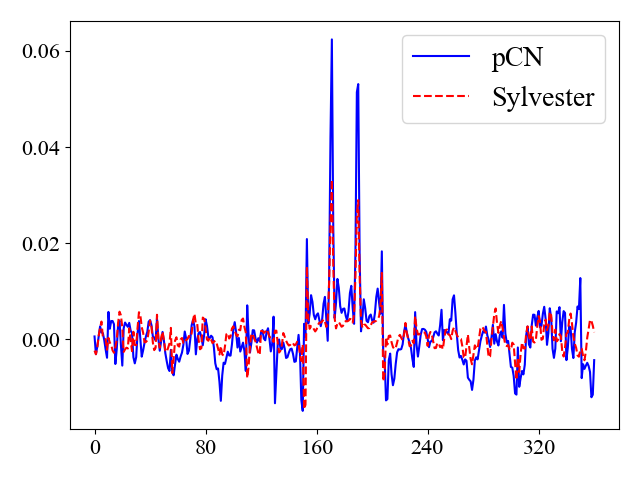}} 
  \caption{\emph{\small  The comparison of the measures obtained by CNF-iVI method and the pCN method. The estimated variance and covariance functions obtained by the pCN algorithm are drawn in blue solid line, and by CNF-iVI algorithm the red dashed line. (a): The background truth of $u$. (b): The mean of approximate posterior obtained by CNF-iVI method. (c): The mean of posterior obtained by pCN method. (d): The covariance function $c(x,y)$ on all the mesh point pairs $\lbrace (x_i, x_{i})\rbrace^{N_g}_{i=1}$;
				(e): The covariance function $c(x,y)$ on the mesh points $\lbrace (x_i, x_{i+40}) \rbrace^{N_g-40}_{i=1}$;
				(f): The covariance function $c(x,y)$ on the mesh points $\lbrace (x_i, x_{i+80}) \rbrace^{N_g-80}_{i=1}$.}}
\label{pic:examplebefore1}
\end{figure}

The model is trained on the dataset $D_{\text{train}}$ with Algorithm~\ref{alg B}. We set $M=10$ and $N_u=20$, with the number of training steps $K=50000$ in Algorithm \ref{alg B}. The initial learning rate is set to $\alpha_0=0.001$, following a step decay learning rate schedule with a multiplicative factor of $\tau = 0.95$ after every $1000$ iterations.
To evaluate the performance of the trained model, we generate a test dataset $D_{test} = \{(\bm{x}_1,\bm{d}_1),(\bm{x}_2,\bm{d}_2),\ldots,(\bm{x}_{N_{test}},\bm{d}_{N_{test}})\}$ with $N_{test}=100$ following the same data generation process with training dataset $D_{train}$. For each pair $(\bm{x}_i, \bm{d}_i)$ in the test dataset $D_{test}$, the corresponding approximate posterior $\nu_{\lambda}(\bm{\bm{v}_i})$ can be computed directly using the trained conditional network $\mathcal{N}_{\lambda}(\bm{v})$.

\begin{figure}[ht!]
\centering
\subfloat[background truth]{
\includegraphics[ keepaspectratio=true, width=0.33\textwidth, clip=true]{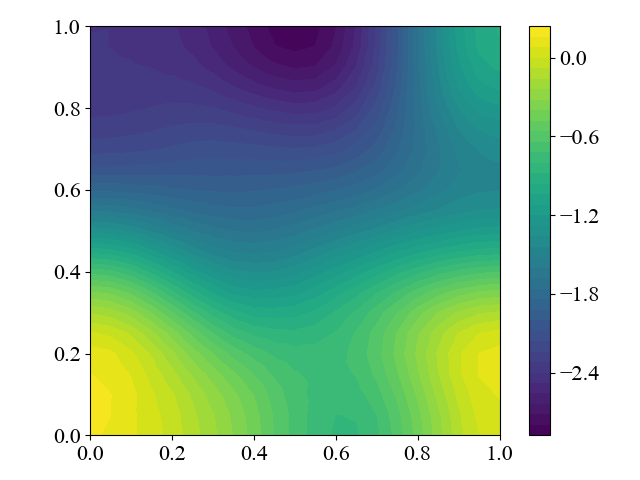}} 
\subfloat[mean (CNF-iVI)]{
\includegraphics[ keepaspectratio=true, width=0.33\textwidth, clip=true]{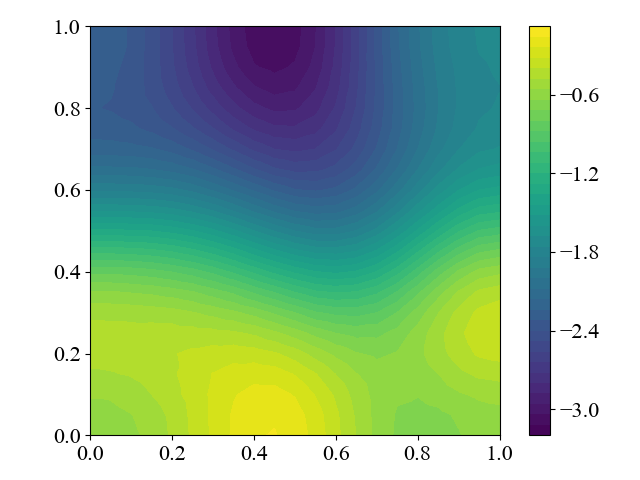}} 
\subfloat[mean (pCN)]{
\includegraphics[ keepaspectratio=true, width=0.33\textwidth, clip=true]{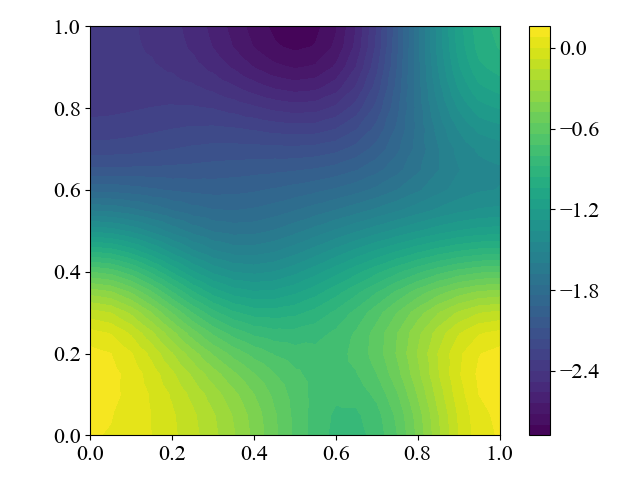}} \\
\subfloat[variance]{
\includegraphics[ keepaspectratio=true, width=0.31\textwidth, clip=true]{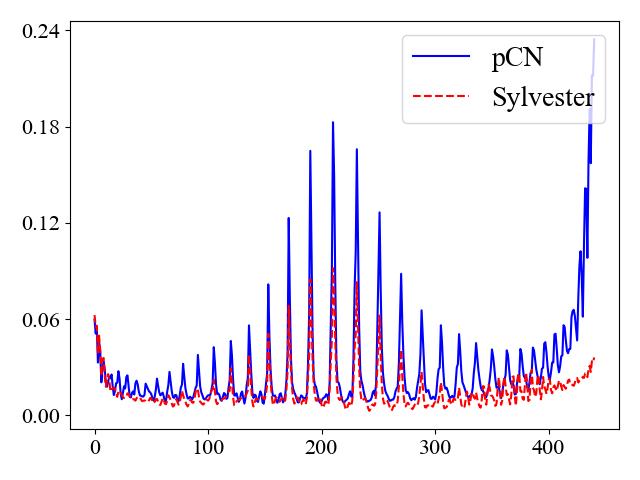}} 
\subfloat[covariance]{
\includegraphics[ keepaspectratio=true, width=0.31\textwidth, clip=true]{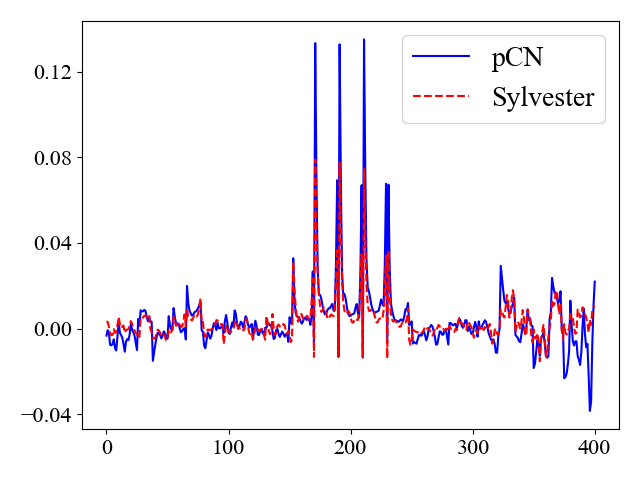}} 
\subfloat[covariance]{
\includegraphics[ keepaspectratio=true, width=0.31\textwidth, clip=true]{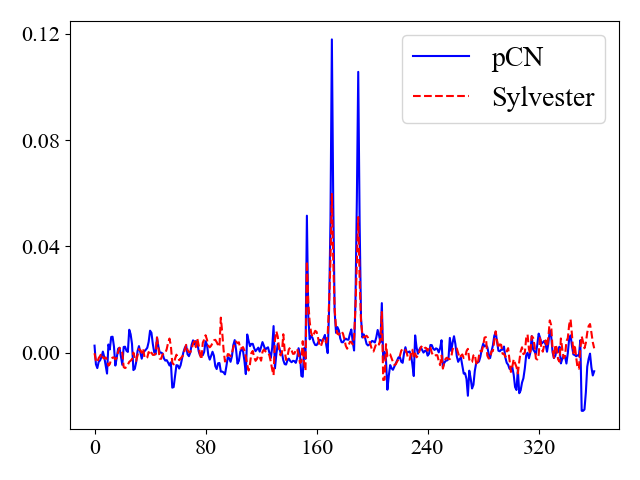}} 
\caption{\emph{\small  The comparison of the measures obtained by CNF-iVI method and the pCN method. The estimated variance and covariance functions obtained by the pCN algorithm are drawn in blue solid line, and by CNF-iVI algorithm the red dashed line. (a): The background truth of $u$. (b): The mean of approximate posterior obtained by CNF-iVI method. (c): The mean of posterior obtained by pCN method. (d): The covariance function $c(x,y)$ on all the mesh point pairs $\lbrace (x_i, x_{i})\rbrace^{N_g}_{i=1}$;
(e): The covariance function $c(x,y)$ on the mesh points $\lbrace (x_i, x_{i+40}) \rbrace^{N_g-40}_{i=1}$;
(f): The covariance function $c(x,y)$ on the mesh points $\lbrace (x_i, x_{i+80}) \rbrace^{N_g-80}_{i=1}$.}}
\label{pic:examplebefore2}
\end{figure}

We now compare the approximate posterior obtained from the conditional network with that produced by the pCN algorithm. For each pCN run, we perform $3\times 10^6$ sampling iterations. Owing to the high computational cost of the pCN algorithm, we restrict our comparison to the posteriors corresponding to two examples in the test dataset.

The conditional network yields only a rough estimate, resulting in a substantial discrepancy between the mean of the approximate posterior and that obtained from the pCN algorithm. Consequently, the covariance function produced by the conditional network also differs significantly from that generated by the pCN algorithm. Figures~\ref{pic:examplebefore1} and \ref{pic:examplebefore2} present qualitative results of the approximate posterior inferred using the trained conditional network $\mathcal{N}_{\lambda}(\bm{v})$ for two different measurement pairs. To further demonstrate the effectiveness of the CNF-iVI method, we compute the re-simulation error:
\begin{align}\label{equ:resimulation error}
\begin{split}
    \textbf{re-simulation error:}&= E_{(\bm{x},\bm{d})}\left(E_{\nu_{\lambda}(\bm{v}^{(\bm{x}, \bm{d})})}\lVert \mathcal{S}\mathcal{G}(u)-\mathcal{S}\mathcal{G}(u_{truth}^{(\bm{x}, \bm{d})})\rVert\right)\\
    &\approx \frac{1}{N_{test}N_{samples}}\sum\limits_{i=1}^{N_{test}}\sum\limits_{j=1}^{N_{samples}}\lVert \mathcal{S}\mathcal{G}(u_{ij})-\mathcal{S}\mathcal{G}(u_{truth}^{(\bm{x}_i,\bm{d}_i)})\rVert\\
    &=0.10834.
\end{split}
\end{align}
where $u_{truth}^{(\bm{x},\bm{d})}$ represents the ground truth and $\bm{v}^{(\bm{x},\bm{d})}$ denotes the vector corresponding to the measurement information $(\bm{x}, \bm{d})$, and $N_{samples}=500$ for this evaluation.

The visual evidence presented in Figures~\ref{pic:examplebefore1} and~\ref{pic:examplebefore2}, together with the quantitative results provided by the re-simulation error in (\ref{equ:resimulation error}), demonstrates the effectiveness of the proposed method in roughly capturing the properties of the posterior conditioned on the measurement information $(\bm{x},\bm{d})$. This further highlights the efficiency of the conditional functional normalizing flow in providing a rough inference of the posterior.

\begin{figure}[ht!]
	\centering
	\subfloat[background truth]{
		\includegraphics[ keepaspectratio=true, width=0.33\textwidth, clip=true]{PIC/conditional/truth1.png}} 
	\subfloat[mean (retrained model)]{
		\includegraphics[ keepaspectratio=true, width=0.33\textwidth, clip=true]{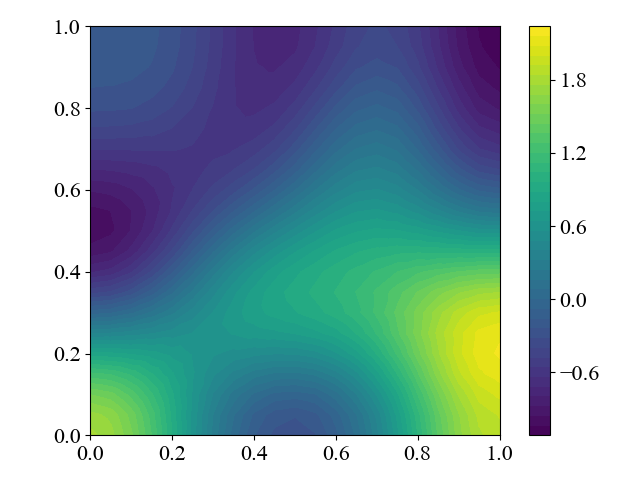}} 
        \subfloat[mean (pCN)]{
		\includegraphics[ keepaspectratio=true, width=0.33\textwidth, clip=true]{PIC/conditional/pcnmean1.png}} 
	\\
	\subfloat[variance]{
		\includegraphics[ keepaspectratio=true, width=0.31\textwidth, clip=true]{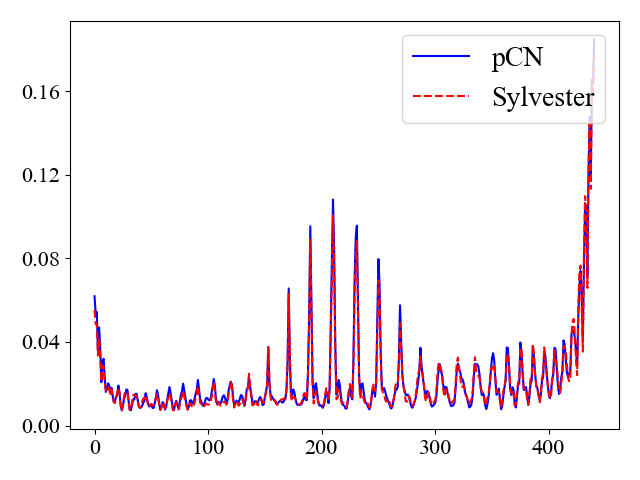}} 
	\subfloat[covariance]{
		\includegraphics[ keepaspectratio=true, width=0.31\textwidth, clip=true]{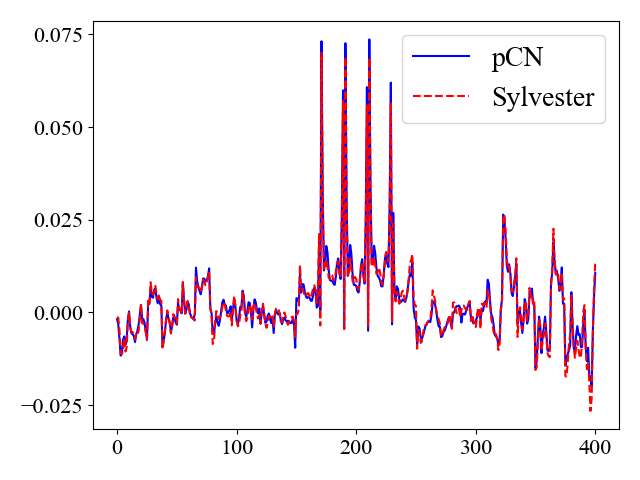}} 
        \subfloat[covariance]{
		\includegraphics[ keepaspectratio=true, width=0.31\textwidth, clip=true]{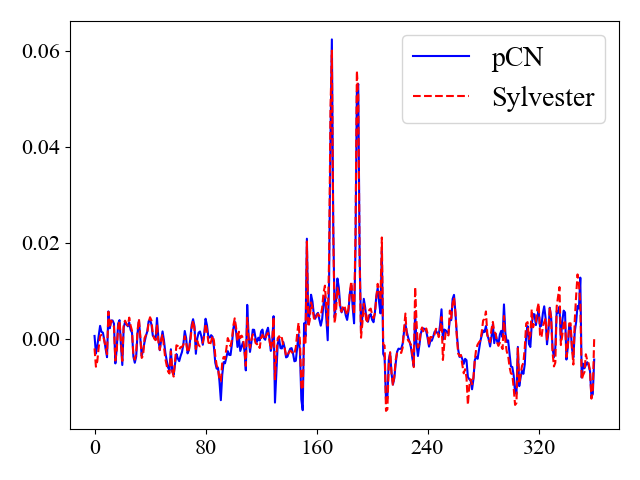}} 
	\caption{\emph{\small  The comparison of the measures obtained by retrained NF-iVI method and the pCN method. The estimated variance and covariance functions obtained by the pCN algorithm are drawn in blue solid line, and by retrained NF-iVI algorithm the red dashed line. (a): The background truth of $u$. (b): The mean of approximate posterior obtained by retrained NF-iVI method. (c):The mean of posterior obtained by pCN method. (d): The covariance function $c(x,y)$ on all the mesh point pairs $\lbrace (x_i, x_{i})\rbrace^{N_g}_{i=1}$;
				(e): The covariance function $c(x,y)$ on the mesh points $\lbrace (x_i, x_{i+40}) \rbrace^{N_g-40}_{i=1}$;
				(f): The covariance function $c(x,y)$ on the mesh points $\lbrace (x_i, x_{i+80}) \rbrace^{N_g-80}_{i=1}$}.}
       \label{pic:retrainexample1}
\end{figure}

\subsection{Further Training for Specific Data}\label{sec5.2}
For the inverse problem corresponding to a specific measurement pair, we can refine the approximate posterior obtained from the trained conditional network $\mathcal{N}_{\lambda}(\bm{v})$ using a retraining method. Consider the test dataset $\{(\bm{x}_1,\bm{d}_1),(\bm{x}_2,\bm{d}_2),\ldots,(\bm{x}_{N_{test}},\bm{d}_{N_{test}})\}$ used in Subsection \ref{sec5.1}, we first convert each pair $(\bm{x}_i,\bm{d}_i)$ into its corresponding vector $\bm{v}_i$, and then inputting $\bm{v}_i$ into the trained neural network $\mathcal{N}_{\lambda}(\bm{v})$ to obtain initial approximations $\{\nu_{\lambda}(\bm{v}_1),\nu_{\lambda}(\bm{v}_2),\ldots,\nu_{\lambda}(\bm{v}_{N_{test}})\}$ of the posteriors. Subsequently, we refine these approximate posteriors with Algorithm \ref{alg A} using $\theta_0^{(i)}=\mathcal{N}_{\lambda}(\bm{v}_i)$ as a start point, resulting in a new set of refined posteriors $\{\hat{\nu}_1,\hat{\nu}_2,\ldots,\hat{\nu}_{N_{test}}\}$.

\begin{figure}[ht!]
\centering
\subfloat[background truth]{
\includegraphics[ keepaspectratio=true, width=0.33\textwidth, clip=true]{PIC/conditional/truth2.png}} 
\subfloat[mean (retrained model)]{
\includegraphics[ keepaspectratio=true, width=0.33\textwidth, clip=true]{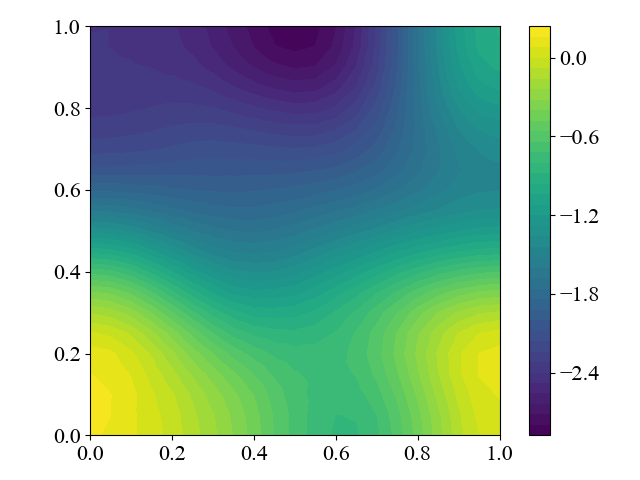}} 
\subfloat[mean (pCN)]{
\includegraphics[ keepaspectratio=true, width=0.33\textwidth, clip=true]{PIC/conditional/pcnmean2.png}} 
\\
\subfloat[variance]{
\includegraphics[ keepaspectratio=true, width=0.31\textwidth, clip=true]{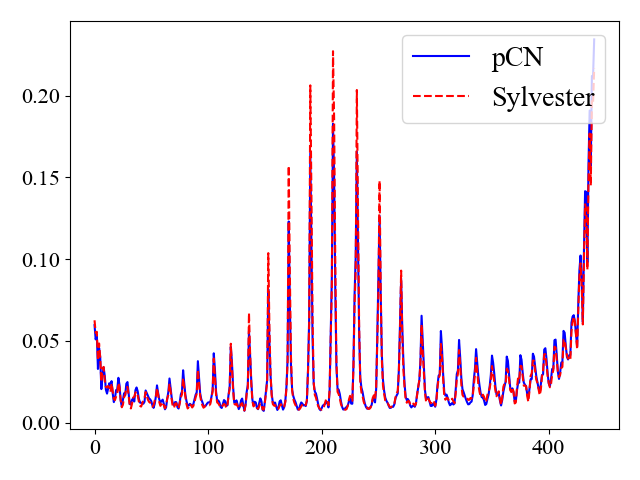}} 
\subfloat[covariance]{
\includegraphics[ keepaspectratio=true, width=0.31\textwidth, clip=true]{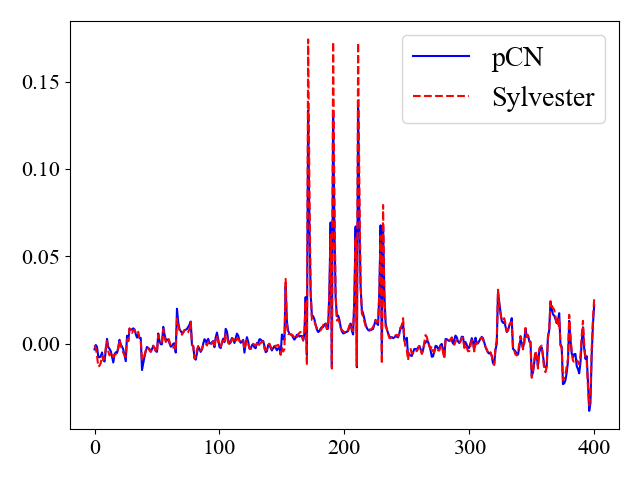}} 
\subfloat[covariance]{
\includegraphics[ keepaspectratio=true, width=0.31\textwidth, clip=true]{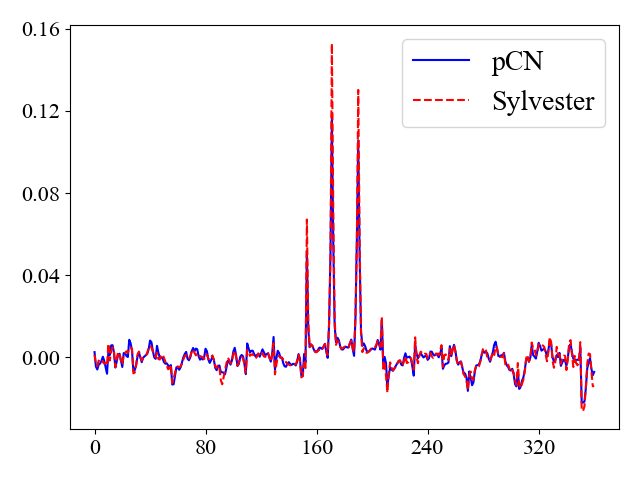}} 
\caption{\emph{\small  The comparison of the measures obtained by retrained NF-iVI method and the pCN method. The estimated variance and covariance functions obtained by the pCN algorithm are drawn in blue solid line, and by retrained NF-iVI algorithm the red dashed line. (a): The background truth of $u$. (b): The mean of approximate posterior obtained by retrained NF-iVI method. (c):The mean of posterior obtained by pCN method. (d): The covariance function $c(x,y)$ on all the mesh point pairs $\lbrace (x_i, x_{i})\rbrace^{N_g}_{i=1}$;
(e): The covariance function $c(x,y)$ on the mesh points $\lbrace (x_i, x_{i+40}) \rbrace^{N_g-40}_{i=1}$;
(f): The covariance function $c(x,y)$ on the mesh points $\lbrace (x_i, x_{i+80}) \rbrace^{N_g-80}_{i=1}$}.}
\label{pic:retrainexample2}
\end{figure}

We use Adam to retrain the approximate posterior. Following Algorithm \ref{alg A}, the functional Sylvester flow is retrained for 500 iterations $(K=500)$ using 30 samples $(N=30)$ in each iteration, with an initial learning rate of $\alpha_0=0.001$. The learning rate follows a step decay schedule with a multiplicative factor of $\tau=0.9$ after every 100 iterations.
Figures \ref{pic:retrainexample1} and \ref{pic:retrainexample2} present the refined results for the same scenarios shown in Figures \ref{pic:examplebefore1} and \ref{pic:examplebefore2}, obtained by retraining the initial output of the conditional flow model from Subsection \ref{sec5.1}.

To numerically compare the covariance operators, we also compute the relative errors of the variance functions and covariance functions between the estimated posteriors generated by the pCN algorithm and the retrained NF-iVI method, as reported in Table~\ref{Table:retrain}. In the table, the notation $\bm{c}$ denotes the total relative error, while $\lbrace c(x_i, x_{i+k}) \rbrace^{N_g-k}_{i=1}$ represents the values of the covariance function evaluated at the pairs of mesh points $\lbrace (x_i, x_{i+k}) \rbrace^{N_g-k}_{i=1}$ ($k=0, 40, 80$). The numerical values reported correspond to the relative errors between the vectors obtained using the pCN method and those obtained using the retrained NF-iVI method, for two different measurement pairs, respectively.

Furthermore, we provide numerical evidence to support the effectiveness of the retraining process. The re-simulation error of the retrained model is computed as
\begin{align}\label{equ:resimulation error2}
\begin{split}
    \textbf{re-simulation error:}&= E_{(\bm{x},\bm{d})}\left(E_{\hat{\nu}_{(\bm{x},\bm{d})}}\lVert \mathcal{S}\mathcal{G}(u)-\mathcal{S}\mathcal{G}(u_{truth}^{(\bm{x}, \bm{d})})\rVert\right)\\
    &\approx \frac{1}{N_{test}N_{samples}}\sum\limits_{i=1}^{N_{test}}\sum\limits_{j=1}^{N_{samples}}\lVert \mathcal{S}\mathcal{G}(u_{ij})-\mathcal{S}\mathcal{G}(u_{truth}^{(\bm{x}_i,\bm{d}_i)})\rVert\\
    &=0.02393,
\end{split}
\end{align}
where $u_{truth}^{(\bm{x},\bm{d})}$ represents the ground truth and $\hat{\nu}_{(\bm{x},\bm{d})}$ denotes the approximate posterior obtained from the retrained flow models corresponding to the measurement pair $(\bm{x}, \bm{d})$. We set $N_{samples}=500$ for the evaluation, consistent with (\ref{equ:resimulation error}). This result demonstrates a significant improvement over the re-simulation error \eqref{equ:resimulation error} obtained without retraining.

The visual evidence presented in Figures~\ref{pic:retrainexample1} and~\ref{pic:retrainexample2}, together with the re-simulation error in (\ref{equ:resimulation error2}) and the relative errors reported in Table~\ref{Table:retrain}, demonstrates the effectiveness of the retraining method.

\begin{table}[H]
	\renewcommand{\arraystretch}{1.5}
	\centering
		\caption{\emph{\small The relative errors between the covariance matrix, variance function, and covariance functions.}} \label{table:relative}
	\begin{tabular}{c|cccc}
		\hline $\text{Relative Error} $& $\bm{c}$  & $\lbrace c(x_i, x_i)\rbrace^{N_g}_{i=1}$& $\lbrace c(x_i, x_{i+40}) \rbrace^{N_g-40}_{i=1}$& $\lbrace c(x_i, x_{i+80}) \rbrace^{N_g-80}_{i=1}$\\
		\hline $\text{Example 1}$& $0.0672$& $0.0150$& $0.0277$& $0.0579$\\
  \hline $\text{Example 2}$& $0.0722$& $0.0195$& $0.0413$& $0.0633$\\
		\hline
	\end{tabular}
 \label{Table:retrain}
\end{table}

\section{Conclusion}\label{sec6}
In this paper, we introduce the NF-iVI and CNF-iVI methods in an infinite-dimensional setting, providing an efficient computational framework for applying variational inference to inverse problems in function spaces. The NF-iVI method constructs a transformation that guarantees measure equivalence while retaining sufficient flexibility in the transformed measure. Building on NF-iVI, the CNF-iVI method incorporates a conditional neural network to control the parameters of the functional normalizing flow, thereby significantly reducing computational costs.

The proposed NF-iVI and CNF-iVI approaches are applicable to PDE inverse problems with Gaussian priors and yield explicit representations of the approximate posterior measure. {\color{black}We have successfully applied these methods to three inverse problems: a simple smooth equation, the steady-state Darcy flow problem, and electrical impedance tomography.} Although CNF-iVI may initially provide a relatively rough approximation of the posterior distribution, the application of the retraining strategy can substantially improve its accuracy.

The current NF-iVI and CNF-iVI frameworks are based on transformations such as planar flow, Householder flow, Sylvester flow, and projected transformation flow. Future work may explore alternative classes of neural networks that satisfy the general theorem established in infinite-dimensional spaces. For example, in Euclidean settings, the non-centered parameterization formulation \cite{papaspiliopoulos2003non} has been employed as a transformation within normalizing flows. Under the non-centered parameterization framework, different types of neural networks can be introduced at each stage of the transformation, with carefully designed layers to enhance model expressiveness. Consequently, developing NF-iVI methods with measure-equivalence properties based on non-centered parameterization and integrating them with other neural architectures in function spaces, such as Fourier neural operators \cite{li2020fourier}, represents a promising direction for future research.

\section{Appendix}
Since the parameter $\theta$ plays no explicit role in the discussion, we omit explicit reference to it. Consequently, we introduce the following abbreviations:
\begin{itemize}
\item The operator $f_{\theta_n}^{(n)}$ is simplified as $f^{(n)}$.
\item The operator $\mathcal{F}_{\theta_n}^{(n)}$ is simplified as $\mathcal{F}^{(n)}$.
\item The operator $f_{\theta}$ is simplified as $f$.
\end{itemize}

\subsection{Functional Determinant}\label{subsec:functionaldeterminant}
To obtain a precise expression for the Radon-Nikodym derivative between $\mu_f$ and $\mu_0$ as discussed in Subsection \ref{subsec:normalizing flow in function space}, we will introduce the concept of the Fredholm-Carleman determinant. This determinant is applicable to operators of the form $I + \mathcal{K}$, where $\mathcal{K}$ is an operator with finite dimensional range mapping from a separable Hilbert space $\mathcal{H}_u$ to itself, and $I$ denotes the identity operator.
\begin{definition}
Let $\mathcal{H}_u$ denotes a separable Hilbert space, and let $\mathcal{K}$ be a finite-dimensional linear operator mapping from $\mathcal{H}_u$ to $\mathcal{H}_u$, with its range being $\mathcal{K}(\mathcal{H}_u)$. We introduce the following notations:
\begin{align*}
{\det}_{1}(I+\mathcal{K})=\det((I+\mathcal{K})|_{\mathcal{K}(\mathcal{H}_u)})
\end{align*}
and
\begin{align*}
{\det}_2(I+\mathcal{K})=\det((I+\mathcal{K})|_{\mathcal{K}(\mathcal{H}_u)})\exp(-\text{trace}(\mathcal{K})),
\end{align*}
where $(I+\mathcal{K})|_{\mathcal{K}(\mathcal{H}_u)}$ represents the matrix representation of the operator $I+\mathcal{K}$ restricted to the subspace $\mathcal{K}(\mathcal{H}_u)$.
\end{definition}
We refer to ${\det}_{1}(I+\mathcal{K})$ as the Fredholm-Carleman determinant associated with the operator $I+\mathcal{K}$, and ${\det}_2(I+\mathcal{K})$ as the regularized Fredholm-Carleman determinant of the same operator. For further details, the reader is directed to Chapter 6.4 in the reference \cite{gaussianmeasure}.

\subsection{Proof of Theorem \ref{thm:RN}}\label{subsec:equality}
$\quad$
\begin{proof}
We begin by proving the equivalence of measures. The proof can be divided in two steps. 

\textbf{Step 1: }Note that $f^{(n)}=I+\mathcal{F}^{(n)}$, we have
\begin{align*}
    f^{(N)} \circ \cdots \circ f^{(1)} =& (I+\mathcal{F}^{(N)}) \circ (I+\mathcal{F}^{(N-1)}) \circ \cdots \circ (I+\mathcal{F}^{(1)})\\
    =&(I+\mathcal{F}^{(N-1)})\circ \cdots \circ (I+\mathcal{F}^{(1)}) + \mathcal{F}^{(N)} \circ (I+\mathcal{F}^{(N-1)}) \circ \cdots \circ (I+\mathcal{F}^{(1)})\\
    =&I+\mathcal{F}^{(1)}+\mathcal{F}^{(2)} \circ (I+\mathcal{F}^{(1)}) + \mathcal{F}^{(3)} \circ (I+\mathcal{F}^{(2)}) \circ (I+\mathcal{F}^{(1)})+\cdots +\\
    &\mathcal{F}^{(N)} \circ (I+\mathcal{F}^{(N-1)}) \circ \cdots \circ (I+\mathcal{F}^{(1)})\\
    =&I+\mathcal{T},
\end{align*}
where $\mathcal{T}=\mathcal{F}^{(1)}+\mathcal{F}^{(2)} \circ (I+\mathcal{F}^{(1)}) + \mathcal{F}^{(3)} \circ (I+\mathcal{F}^{(2)}) \circ (I+\mathcal{F}^{(1)}) + \cdots +\mathcal{F}^{(N)} \circ (I+\mathcal{F}^{(N-1)}) \circ \cdots \circ (I+\mathcal{F}^{(1)})$.
Since $\text{Im}(\mathcal{F}^{(n)}) \subset \mathcal{H}$, we have $\text{Im}(\mathcal{T}) \subset \mathcal{H}$.

\textbf{Step 2: }We can find that $I+\mathcal{T}$ is Fr\'{e}chet differentiable and
\begin{align*}
    D(I+\mathcal{T})(u)&=D(I+\mathcal{F}^{(N)})\circ(I+\mathcal{F}^{(N-1)})\circ \cdots \circ(I+\mathcal{F}^{(1)})(u)\\
    &=(I+D\mathcal{F}^{(N)})(u^{N-1})\circ (I+D\mathcal{F}^{(N-1)})(u^{N-2})\circ \cdots \circ (I+D\mathcal{F}^{(1)})(u)\\
    &=\Gamma^{(N)} \circ \Gamma^{(N-1)}\circ \cdots \circ \Gamma^{(1)}. 
\end{align*}
where $\Gamma^{(n)}=(I+D\mathcal{F}^{(n)})(u^{n-1})$, and $u^{n-1}=(I+\mathcal{F}^{(n-1)})\circ\cdots\circ(I+\mathcal{F}^{(1)})(u)$ for $n=1, 2, \ldots, N$. Since the point spectrum of $D\mathcal{F}^{(n)}(u)$ does not lie within the interval $(-\infty,-1]$ for all $u \in \mathcal{H}_u$, it follows that the mapping $\Gamma^{(n)}$ is injective for $n = 1, 2, \ldots, N$ and $u \in \mathcal{H}_u$. Furthermore, since $D(I+\mathcal{T})(u)=\Gamma^{(N)} \circ \Gamma^{(N-1)}\circ\cdots\circ \Gamma^{(1)}$, we conclude that $D(I+\mathcal{T})(u)$ is injective for all $u \in \mathcal{H}_u$. 

Relying on Steps 1 and 2 and Example 10.27 of \cite{MR2663405}, we obtain that $\mu \sim \mu_0$ where $\mu=\mu_0 \circ (f^{(N)} \circ f^{(N-1)} \circ \cdots \circ f^{(1)})^{-1}$.

Next, considering the Radon--Nikodym derivative of \(\mu_{f_{\theta}}\) with respect to \(\mu_0\). We first examine the case \(N = 1\), identifying \(\mathcal{F}^{(1)}\) with \(\mathcal{F}\). Let $\{e_1,e_2,\ldots,e_M\}$ be a set of standard orthogonal bases of the space $\mathcal{F}(\mathcal{H}_u) \subset \mathcal{H}$. From Corollary 6.6.8 of \cite{gaussianmeasure}, we know that
\begin{align*}
\frac{d(\mu_0\circ (I+\mathcal{F})^{-1})}{d\mu_0}((I+\mathcal{F})(u))=\frac1{\Lambda_{\mathcal{F}}\left(u\right)},
\end{align*}
where
\begin{align*}
\Lambda_{\mathcal{F}}(u):=\left|{\det}_2\left(I+D\mathcal{F}(u)\right)\right|\exp\left[\delta \mathcal{F}(u)-\frac12\lVert\mathcal{F}(u)\rVert_{\mathcal{H}}^2\right],
\end{align*}
and
\begin{align*}
\delta \mathcal{F}(u)=\mathrm{trace}D \mathcal{F}(u)-\sum_{i=1}^M\widehat{e}_i(u)\mathcal{F}_i(u),
\end{align*}
where for any $i \in \{1,2,\ldots,M\}$, $\mathcal{F}_i(u)=\langle \mathcal{F}(u),e_i \rangle_{\mathcal{H}}$, $\widehat{e}_i(u)=\mathcal{C}_0^{-1} e_i (u)$ ($\widehat{h}$
is called the $\mu_0$-measurable linear functional generated by $h$ for any $h \in \mathcal{H}$, detail definition can be found in p.60 of \cite{gaussianmeasure}).
By simplifying the calculation, we denote that
\begin{align*}
\begin{split}
\left\langle u,\mathcal{F}(u) \right\rangle_{\mathcal{H}}:=\sum\limits_{i=1}^{M}\widehat{e}_i(u)\mathcal{F}_i(u)
=\left\langle \sum\limits_{i=1}^{M}\widehat{e}_i(u)e_i ,\sum\limits_{i=1}^{M}\mathcal{F}_i(u)e_i\right\rangle_{\mathcal{H}}.
\end{split}
\end{align*}
\begin{remark}
    For $u \sim \mathcal{N}(0, \mathcal{C}_0)$, the event $u \in \mathcal{H}$ is not guaranteed to occur, and thus the inner product $\left\langle u,\mathcal{F}(u) \right\rangle_{\mathcal{H}}$ is not always well-defined. However, when $u \in \mathcal{H}$, the identity $\left\langle u,\mathcal{F}(u) \right\rangle_{\mathcal{H}}=\sum\limits_{i=1}^{M}\widehat{e}_i(u)\mathcal{F}_i(u)$ holds. In our numerical experiments, we generate samples from the Gaussian distribution $\mathcal{N}(0,\mathcal{C}_0)$ by truncating its Karhunen-Loève (KL) expansion. This sampling method ensures that all generated samples $u$ satisfy $u \in \mathcal{H}$. Consequently, the use of $\left\langle u,\mathcal{F}(u) \right\rangle_{\mathcal{H}}=\sum\limits_{i=1}^{M}\widehat{e}_i(u)\mathcal{F}_i(u)$ is appropriate in the numerical context. To simplify the theorem's presentation, we uniformly employ the notation $\left\langle u,\mathcal{F}(u) \right\rangle_{\mathcal{H}}$ to represent this term.
\end{remark}
Hence, we obtain
\begin{align*}
\begin{split}
\Lambda_{\mathcal{F}}(u)=&\left|{\det}_2\left(I+D\mathcal{F}(u)\right)\right|\exp\left[\mathrm{trace}D \mathcal{F}(u)-\frac12\lVert\mathcal{F}(u)\rVert_{\mathcal{H}}^2-\sum_{i=1}^M\widehat{e}_i(u)\mathcal{F}_i(u)\right]\\
=&\left|{\det}_2\left(I+D\mathcal{F}(u)\right)\right|\exp\left[\mathrm{trace}D \mathcal{F}(u)-\frac12\lVert\mathcal{F}(u)\rVert_{\mathcal{H}}^2 - \langle u,\mathcal{F}(u) \rangle_{\mathcal{H}} \right]\\
=&\left|{\det}_{1}\left(I+D\mathcal{F}(u)\right)\right|\exp\left[-\frac12\lVert\mathcal{F}(u)\rVert_{\mathcal{H}}^2 - \langle u,\mathcal{F}(u) \rangle_{\mathcal{H}} \right].
\end{split}
\end{align*}

Next, considering the case where $N\neq 1$, we denote $f=I+\mathcal{T}$ as in Subsection \ref{subsec:equality}, and note that
\begin{align*}
\begin{split}
Df(u)&=D(f_N \circ \cdots \circ f_1)(u)\\
&=Df_N(f_{N-1}\circ \cdots \circ f_1(u)) \circ Df_{N-1}(f_{N-2}\circ \cdots \circ f_1(u)) \cdots \circ Df_1(u),
\end{split}
\end{align*}
we find that
\begin{align*}
{\det}_{1}(Df(u))={\det}_{1}(Df_N(u_{N-1}))\times {\det}_{1}(Df_{N-1}(u_{N-2})) \cdots \times {\det}_{1}(Df_1(u)),
\end{align*}
where $u_k=f_{k}\circ\cdots \circ f_1(u)$, $k=1, 2, \ldots, N-1$. Note that $\mathcal{T}(u)=f(u)-u$, from Corollary 6.6.8 of \cite{gaussianmeasure}, we have
\begin{align*}
\begin{split}
\frac{d\mu_0 \circ f^{-1}}{d\mu_0}(f(u))&=\frac{1}{\Lambda_{\mathcal{T}}(u)}=\frac{1}{\left|{\det}_{1}\left(I+D\mathcal{T}(u)\right)\right|}\exp\left[\frac12\lVert\mathcal{T}(u)\rVert_{\mathcal{H}}^2 + \langle u,\mathcal{T}(u) \rangle_{\mathcal{H}} \right]\\
&=\prod\limits_{k=1}^{N}\left | {\det}_{1}(Df_k(u_{k-1}))\right|^{-1}\exp\left(\frac{1}{2}\langle f(u)-u,f(u)-u\rangle_{\mathcal{H}}+\langle u,u-f(u)\rangle_{\mathcal{H}}\right).
\end{split}
\end{align*}
The proof of the theorem is completed.
\end{proof}

\subsection{Proof of Lemma \ref{lem:bijective}}
$\quad$
\begin{proof}
Since $\mathcal{F}$ is compact and continuous, we obtain that $f$ is a completely continuous field. Let us define
\begin{align*}
    \text{H}(u,t)=t\mathcal{F}(u),
\end{align*}
where $t\in [0,1]$. For any bounded closed set $\mathcal{M}$ in $\mathcal{H}_u$, we can conclude that
\begin{align*}
    \text{H}: \mathcal{M} \times [0,1] \rightarrow \mathcal{H}_u
\end{align*}
is a compact continuous operator. Let us denote
\begin{align*}
    h_t(u) = u + \text{H}(u,t),
\end{align*}
it is straightforward to see that $h_0(u)=u$ and  $h_1(u) = f(u)$.

For any $y \in \mathcal{H}_u$, let $\Omega = B(0, R) \subset \mathcal{H}_u$, where $R$ is a large enough constant. Thus for any $q \in \partial\Omega$ and $t \in [0,1]$, we have
\begin{align*}
    \left\Vert h_t(q) \right\Vert_{\mathcal{H}_u} \geq \left\Vert q \right\Vert_{\mathcal{H}_u} - \left\Vert \text{H}(q,t) \right\Vert_{\mathcal{H}_u} > 0,
\end{align*}
which ensures that $h_t(q) = q + \text{H}(q,t) \neq 0$. By the compact homotopy invariance of the Leray-Schauder degree \cite{MR1987179}, we obtain
\begin{align*}
    \text{deg}(f,\Omega,y) = \text{deg}(I , \Omega, y) = 1.
\end{align*}
Note that for any $u \in \mathcal{H}_u$, the eigenvalues of operator $Df(u)$ are strictly positive. By the definition of the Leray-Schauder degree, it is straightforward to know that
\begin{align*}
    \sum\limits_{z \in f^{-1}(y)}1 = 1.
\end{align*}
This implies that there is a unique point in \( f^{-1}(y) \), and hence, the solution to \( f(u) = y \) is unique within the ball \( B(0,R) \). 
Given the arbitrary nature of $R$, we can conclude that $f$ is a bijective mapping.
\end{proof}

\subsection{Proof of Theorems \ref{thm:planar} and \ref{thm:Sylvester}}\label{subsec:proof of nonlinear}
$\quad$
\begin{proof}
To prove the theorem, it suffices to verify that the proposed method satisfies the four conditions outlined in Theorem \ref{thm:RN}. As functional planar flow is a special case of functional Sylvester flow, we will concentrate on the latter. We need to prove the corresponding operator 
\begin{align*}
\mathcal{F}_{\theta_n}^{(n)}(u)=\mathcal{A}_nh(\mathcal{B}_nu+b_n)
\end{align*}
satisfies the four conditions of Theorem \ref{thm:RN}, i.e.:
\begin{itemize}
\item The space $\text{Im}(\mathcal{F}_{\theta_n}^{(n)}) \subset \mathcal{H}$, where $\text{Im}(\mathcal{F}_{\theta_n}^{(n)})$ denotes the image of $\mathcal{F}_{\theta_n}^{(n)}$.
\item The operator $\mathcal{F}_{\theta_n}^{(n)}$ has finite dimensional range.
\item The operator $f_{\theta_n}^{(n)}$ is bijective.
\item For any $u \in \mathcal{H}_u$,  all point spectrum of $D\mathcal{F}_{\theta_n}^{(n)}(u)$ are not in $(-\infty,-1]$.
\end{itemize}

Since $\mathcal{A}_n$ satisfies $\text{Im}(\mathcal{A}_n) \subset \mathcal{H}$, we know that $\mathcal{F}_{\theta_n}^{(n)} \subset \mathcal{H}$. At the same time, since the domain $D(\mathcal{A}_n)$ of operator $\mathcal{A}_n$ has finite dimensions and $\mathcal{A}_n$ is a linear operator, we can conclude that $\text{Im}(\mathcal{A}_n)$ has finite dimensions, which means that $\mathcal{F}_{\theta_{n}}^{(n)}$ has finite dimensional range. The above discussion verifies the first two conditions required by Theorem \ref{thm:RN}.

Since $\mathcal{A}_n$ is a linear operator with finite dimensional range, it must be a bounded operator. Note that $h(x)=\tanh(x) \subset (-1,1)$, so  $h(\mathcal{B}_nu+b_n)$ is bounded in $\mathbb{R}^M$ for any $u \in \mathcal{H}_{u}$. Consequently, $\mathcal{A}_nh(\mathcal{B}_nu+b_n)$ is bounded in $\mathcal{H}_{u}$, implies that the range of $\mathcal{F}_n$ is bounded in $\mathcal{H}_{u}$. By Lemma \ref{lem:bijective}, we know that $I+\mathcal{F}_{\theta_n}^{(n)}$ is bijective. This verifies the third condition required by Theorem \ref{thm:RN}.

Finally, we are able to determine that
\begin{align*}
D\mathcal{F}_{\theta_n}^{(n)}(u)&=\mathcal{A}_nDh(\mathcal{B}_nu+b_n)\mathcal{B}_n 
=\mathcal{A}_n\mathrm{diag}(h'(\mathcal{B}_nu+b_n))\mathcal{B}_n,
\end{align*}
and $\mathcal{A}_n\text{diag}(h'(\mathcal{B}_nu+b_n))\mathcal{B}_n$ has the same point spectrum with $\text{diag}(h'(\mathcal{B}_nu+b_n))\mathcal{B}_n\mathcal{A}_n$ (see Theorem \ref{thm:change} for details). Furthermore, we recognize that $\text{diag}(h'(\mathcal{B}_nu+b_n))$ is a diagonal matrix with diagonal elements within the interval $(0,1)$. Therefore, combined with the conditions that the eigenvalues of $\mathcal{B}_n\mathcal{A}_n$ are not in $(-\infty,1]$, we find that the eigenvalues of $\mathrm{diag}(h'(\mathcal{B}_nu+b_n))\mathcal{B}_n\mathcal{A}_n$ are not in $(-\infty,1]$, so the point spectrum of $\mathcal{A}_n\text{diag}(h'(\mathcal{B}_nu+b_n))\mathcal{B}_n$ is not in $(-\infty,1]$. Thus, we have verified the final condition required by Theorem \ref{thm:RN}.

In conclusion, $\mathcal{F}_{\theta_n}^{(n)}(u)=\mathcal{A}_nh(\mathcal{B}_nu+b_n)$ satisfies all the conditions of Theorem \ref{thm:RN}, which completes the proof. 
\end{proof}

\subsection{Proof of Theorems \ref{thm:Householder} and \ref{thm:pro}}
$\quad$
\begin{proof}
To prove the theorem, it suffices to verify that the proposed method satisfies the four conditions of Theorem \ref{thm:RN}. As functional Householder flow is a special case of functional projected transformation flow, we will focus on the latter. We need to prove that the corresponding operator
\begin{align*}
\mathcal{F}_{\theta_n}^{(n)}(u)=\mathcal{Q}R_n(\mathcal{P}u+b_n)
\end{align*}
satisfies the four conditions of Theorem \ref{thm:RN}. Similar to the approach used in Subsection \ref{subsec:proof of nonlinear}, we are able to prove the first two conditions for functional projected transformation flow. The details are omitted here for brevity. 

For the third and fourth conditions, note that 
\begin{align*}
D\mathcal{F}_{\theta_n}^{(n)}(u)=\mathcal{Q}R_n\mathcal{P},
\end{align*}
and the operator $\mathcal{Q}R_n\mathcal{P}$ has the same eigenvalues with the matrix $R_n$, the point spectrum of $D\mathcal{F}_{\theta_n}^{(n)}(u)$ is not in $(-\infty,1]$. Thus, Lemma \ref{lem:libijective} implies that $I + \mathcal{F}_{\theta_n}^{(n)}$ is invertible. 
Based on the above discussion, we can conclude that $\mathcal{F}_{\theta_n}^{(n)}(u)=\mathcal{Q}R_n(\mathcal{P}u+b_n)$ satisfies the conditions of Theorem \ref{thm:RN}, which completes the proof of the theorem.
\end{proof} 

\subsection{Proof of Theorem \ref{thm:change}}
$\quad$
\begin{proof}
Note that both $\mathcal{A}$ and $\mathcal{B}$ are linear operators, and
\begin{align*}
\mathcal{B}: \mathcal{H}_u \rightarrow \mathbb{R}^M, \quad 
\mathcal{A}: \mathbb{R}^M \rightarrow \mathcal{H}_u.
\end{align*}
Without losing generality, $\mathcal{B}$ and $\mathcal{A}$ can be written as
\begin{align*}
\mathcal{B}v=(\langle \phi_1,v \rangle_{\mathcal{H}_u},\langle \phi_2,v \rangle_{\mathcal{H}_u},\ldots,\langle \phi_M,v \rangle_{\mathcal{H}_u})^{T},
\end{align*}
where $\phi_1, \phi_2, \ldots , \phi_M \in \mathcal{H}_u$, and
\begin{align*}
\mathcal{A}d=d_1\psi_1+d_2\psi_2+ \cdots +d_M\psi_M,
\end{align*}
where $\psi_1, \psi_2, \ldots , \psi_M \in \mathcal{H}_u$ are mutually orthogonal.

Let $\lambda$ be the point spectrum of $\mathcal{A}\mathcal{B}$, and $v_{\lambda}$ be its corresponding eigenfunction.
We have $\mathcal{A}\mathcal{B}v_{\lambda}=\lambda v_{\lambda}$. Obviously, since $\text{Im}(\mathcal{A}_n)\subset \text{span}\{\psi_1, \psi_2,\ldots,\psi_M \}$, let $v_{\lambda}=\sum\limits_{i=1}^{M}\alpha_i\psi_i$, then for any eigenpair $(\lambda, v_{\lambda})$, we have
\begin{align*}
\sum\limits_{k=1}^{M}\sum\limits_{i=1}^{M}\langle \phi_{k},\psi_{i} \rangle \alpha_i \psi_k = \lambda \sum\limits_{i=1}^{M}\alpha_i \psi_i.
\end{align*}
Hence, for any $r = 1, 2,\ldots, M$, we obtain
\begin{align}\label{Thm2.9-1}
\sum\limits_{i=1}^{M}\langle \phi_r,\psi_i \rangle \alpha_i = \lambda \alpha_r.
\end{align}
Equality (\ref{Thm2.9-1}) can be expressed as
\begin{align*}
\begin{pmatrix}
\langle \phi_1,\psi_1 \rangle & \cdots &\langle \phi_M,\psi_1 \rangle\\
\vdots & \ddots & \vdots\\
\langle \phi_1,\psi_M \rangle & \cdots &\langle \phi_M,\psi_M \rangle\\
\end{pmatrix}
\begin{pmatrix}
\alpha_1\\
\vdots\\
\alpha_M
\end{pmatrix}
=\lambda 
\begin{pmatrix}
\alpha_1\\
\vdots\\
\alpha_M
\end{pmatrix},
\end{align*}
which is just the equation $\mathcal{B}\mathcal{A} \alpha = \lambda \alpha$ with $\alpha = (\alpha_1, \alpha_2, \ldots , \alpha_M)^{T}$.
Therefore, the point spectrum of $\mathcal{A}\mathcal{B}$ is one-to-one correspondence with the eigenvalues of $\mathcal{B}\mathcal{A}$, which finishes the proof.
\end{proof}

\subsection{Proof of Theorem \ref{thm:discrete invariance}}
$\quad$
\begin{proof}
For clarity, we show the discretization invariant when $\mathcal{H}_u = L^2(D)$. The proof can be easily adapted to other settings, e.g., $\mathcal{H}_u = H^s(D)$ with $s\geq 0$. We have
\begin{align*}
\mathcal{F}_{\theta_n}^{(n)}: \mathcal{H}_u \rightarrow \mathcal{H}_u,
\end{align*}
and $\mathcal{F}_{\theta_{n}}^{(n)}$ has four different styles:
\begin{itemize}
\item Functional planar flow
\begin{align*}
    \mathcal{F}_{\theta_{n}}^{(n)}(u)=u_nh(\langle w_n,u_n\rangle_{\mathcal{H}_u}+b_n).
\end{align*}
\item Functional Householder flow
\begin{align*}
    \mathcal{F}_{\theta_{n}}^{(n)}(u)=-0.5v_n(\langle v_n,u\rangle_{\mathcal{H}_u}+b_n).
\end{align*}
\item Functional Sylvester flow
\begin{align*}
\mathcal{F}_{\theta_{n}}^{(n)}(u)=\mathcal{A}_nh(\mathcal{B}_nu+b_n).
\end{align*}
\item Functional projected transformation flow
\begin{align*}
\mathcal{F}_{\theta_{n}}^{(n)}(u)=\mathcal{Q}R_n(\mathcal{P}u+b_n).
\end{align*}
\end{itemize}
It is worth noting that functional planar flow and functional Householder flow are special cases of functional Sylvester flow and functional projected transformation flow, respectively. So we only need to prove the discrete invariance of only the latter two flows.

Let $K$ be a compact subset of $C(D)$, and $\{D_j\}_{j=1}^{\infty}$ a sequence of discrete refinements of $D$. To each discretization $D_j$ associate partition $P_{j}^{(1)},\ldots,P_{j}^{(j)} \subset D$, each contains a single, unique point of $D_j$, each has positive Lebesgue measure, and
\begin{align*}
\bigcup_{k=1}^jP_j^{(k)}=D.
\end{align*}

To establish the discretization invariance of the functional Sylvester flow, we aim to demonstrate the following: for any $\epsilon > 0$ and any $a \in K$, there exists an integer $L > 0$ (independent of $a$) such that for all $m > L$, the inequality
\begin{align*}
\lVert\mathcal{A}_nh(\mathcal{B}_na+b_n)-\mathcal{A}_nh(\mathcal{B}_n^{(m)}a_{(m)}(x)+b_n)\rVert_{\mathcal{H}_u}<\epsilon
\end{align*}
holds. Here, $a_{(m)}(x) = (a(x_1),\ldots,a(x_m))$, and $\mathcal{B}_{n}^{(m)}a_{(m)}(x)$ is defined as
\begin{align*}
\mathcal{B}_{n}^{(m)}a_{(m)}(x)=R_{\mathcal{B}}^{n}
\begin{pmatrix}
& \sum\limits_{i=1}^{m}a(x_i)\phi_1(x_i)\left|P_j^{(i)}\right| \\
& \vdots \\
& \sum\limits_{i=1}^{m}a(x_i)\phi_M(x_i)\left|P_j^{(i)}\right| 
\end{pmatrix}.
\end{align*}
Note that $\mathcal{A}_n$ is a bounded linear operator and $h(x)$ is a continuous function, the the only thing we need to prove is that for any $a\in K$ and $m>L$, we have
\begin{align*}
\lVert \mathcal{B}_n a-\mathcal{B}_n^{(m)}a_{(m)}\rVert_{\mathbb{R}^M} < \epsilon,
\end{align*} 
which equals to
\begin{align*}
\left|\sum\limits_{i=1}^{m}a(x_i)\phi_r(x_i)\left|P_j^{(i)}\right|-\langle a,\phi_r \rangle_{\mathcal{H}_u}\right|<\epsilon,
\end{align*}
for all $r = 1, \ldots , M$. We will provide a proof of this claim at the end of the proof.

Similarly, to establish the discretization invariance of the functional projected transformation flow, we aim to demonstrate the following: for any $\epsilon > 0$ and any $a \in K$, there exists an integer $L > 0$ (independent of $a$) such that for all $m > L$, the inequality
\begin{align*}
\lVert \mathcal{Q}R_n\mathcal{P}a-\mathcal{Q}R_n\mathcal{P}_{(m)}a_{(m)}\rVert_{\mathcal{H}_u}<\epsilon
\end{align*}
holds. Here, $a_{(m)}(x) = (a(x_1),\ldots,a(x_m))$, and $\mathcal{P}_{(m)}a_{(m)}(x)$ is defined as
\begin{align*}
\mathcal{P}_{(m)}a_{(m)}(x)=
\begin{pmatrix}
& \sum\limits_{i=1}^{m}a(x_i)\phi_1(x_i)\left|P_j^{(i)}\right| \\
& \vdots \\
& \sum\limits_{i=1}^{m}a(x_i)\phi_M(x_i)\left|P_j^{(i)}\right| 
\end{pmatrix}.
\end{align*}
Note that $\mathcal{Q}$ and $R_n$ are both bounded linear operators, what we need to prove is that
\begin{align*}
\lVert\mathcal{P}a-\mathcal{P}_{(m)}a_{(m)}\rVert_{\mathbb{R}^N}<\epsilon.
\end{align*}
To prove the desired result, it suffices to show that
\begin{align*}
\left|\sum\limits_{i=1}^{m}a(x_i)\phi_r(x_i)\left|P_j^{(i)}\right|-\langle a,\phi_r \rangle_{\mathcal{H}_u}\right|<\epsilon
\end{align*}
for all $r = 1, \ldots , M$.

Thus, combining the proofs for the linear and nonlinear cases, the key of the theorem is to prove that for any $a \in K$, there exists an integer $L>0$ independent of $a$, such that when $m>L$, we have
\begin{align*}
\left|\sum\limits_{i=1}^{m}a(x_i)\phi_r(x_i)\left|P_j^{(i)}\right|-\langle a,\phi_r \rangle_{\mathcal{H}_u}\right|<\epsilon,
\end{align*}
for all $r = 1, \ldots , M$.

Note that $K\subset C(D)$ is compact, we can find $a_1, \ldots , a_W \in K$ such that for any $a \in K$, there exists $w \in \{1, \ldots, W\}$ satisfying
\begin{align*}
\lVert a-a_{w}\rVert_{C(D)} < \frac{\epsilon}{3 |D| \sup_{r,i}|\phi_r(x_i)|}.
\end{align*}
Since $D_j$ is a discrete refinement, by convergence of the Riemann sum, for each $w = 1, \ldots , W$, there exists a positive constant $p_w>0$ such that when $t_w > p_w$, we have
\begin{align*}
\left|\sum\limits_{i=1}^{t_w}a_w(x_i)\phi_r(x_i)|P_j^{(i)}|-\langle a_w,\phi_r \rangle_{\mathcal{H}_u}\right|<\frac{\epsilon}{3},
\end{align*}
for all $r = 1, \ldots , M$.
Note that $a, a_w \in K$ are all continuous function, we have
\begin{align*}
\left|\sum\limits_{i=1}^{t_w}a_w(x_i)\phi_r(x_i)\left|P_j^{(i)}\right|-\sum\limits_{i=1}^{t_w}a(x_i)\phi_r(x_i)\left|P_j^{(i)}\right|\right|
< \|a-a_w\|_{C(D)} |D| \sup_{r,i}|\phi_r(x_i)| < \frac{\epsilon}{3}.
\end{align*}
Now, let $L \geq \max\{p_1,\ldots,p_W,p\}$. Then, for any $m>L$, we obtain
\begin{align*}
\begin{split}
\left|\sum\limits_{i=1}^{m}a(x_i)\phi_r(x_i)\left|P_j^{(i)}\right|-\langle a,\phi_r \rangle_{\mathcal{H}_u}\right|
\leq &\left|\sum\limits_{i=1}^{m}a(x_i)\phi_r(x_i)\left|P_j^{(i)}\right|-\sum\limits_{i=1}^{m}a_w(x_i)\phi_r(x_i)\left|P_j^{(i)}\right|\right|\\
& \, + \left|\sum\limits_{i=1}^{m}a_w(x_i)\phi_r(x_i)\left|P_j^{(i)}\right|-\langle a_w(x),\phi_r(x) \rangle_{\mathcal{H}_u}\right|\\
& \, + \left|\langle a,\phi_r\rangle_{\mathcal{H}_u}-\langle a_w,\phi_r \rangle_{\mathcal{H}_u}\right|\\
<& \epsilon,
\end{split}
\end{align*}
which completes the proof.
\end{proof}

\section*{Acknowledgments}
The authors would like to thank the anonymous referees for their comments and suggestions,
which helped to improve the paper significantly. 
The third author was supported in part by the National Natural Science Foundation of China (Grant Nos. 12322116, 12271428, and 12326606).
The fourth author was supported in part by the National Natural Science Foundation of China (Grant Nos. 12288201 and 12461160275) and by the Science Challenge Project (No. TZ2025006).

\bibliographystyle{amsplain}
\bibliography{references}
\end{document}